\documentclass[a4wide]{article}

\usepackage{a4wide}
\usepackage{authblk}
\usepackage{cite}
\usepackage[british]{babel}
\usepackage[utf8]{inputenc}
\usepackage{xcolor}
\usepackage{graphicx}
\usepackage{epstopdf}
\usepackage[colorlinks=true,linkcolor=blue,citecolor=blue]{hyperref}
\usepackage{xcolor}
\usepackage{amsfonts,amsmath,amsthm}
\usepackage{hyperref}
\usepackage{algpseudocode}
\usepackage{algorithm}
\usepackage{multirow}

\newtheorem{theorem}{Theorem}[section]

\newtheorem{lemma}[theorem]{Lemma}

\theoremstyle{remark}\newtheorem{remark}[theorem]{Remark}

\newcommand{\be}{\begin{equation}}
	\newcommand{\ee}{\end{equation}}

\newcommand{\F}{\mathcal{F}}

\allowdisplaybreaks

\begin{document}
	\title{Kinetic simulated annealing optimization with entropy-based cooling rate}
	
	\author[1]{Michael Herty \thanks{\texttt{herty@igpm.rwth-aachen.de}}}
	\author[2]{Mattia Zanella \thanks{\texttt{mattia.zanella@unipv.it}}}
	\affil[1]{RWTH Aachen University, Institut f\"ur Geometrie und Praktische Mathematik, Templergraben 55, 52062 Aachen, Germany, and,  Extraordinary Professor, Department of Mathematics and Applied Mathematics, University of Pretoria, Private Bag X20, Hatfield 0028, South Africa}
	\affil[2]{University of Pavia, Department of Mathematics "F. Casorati", Via A. Ferrata 5, 27100 Pavia, Italy}

	\date{}
	
	\maketitle
	\abstract
	We present a modified simulated annealing method with a dynamical choice of the cooling temperature. The latter is determined via a closed-loop control and is proven to yield exponential decay of the entropy of the particle system. The analysis is carried out through kinetic  equations for interacting particle systems describing the simulated annealing method in an extended phase space. Decay estimates are derived under the quasi-invariant scaling of the resulting system of Boltzmann-type equations to assess the consistency with their mean-field limit. Numerical results are provided to illustrate and support the theoretical findings.
	\medskip	
		
	
	\section{Introduction}
	
	The subject of this paper is the  solution of 
	\begin{equation}
		\label{eq:optim}
		x^* \in \textrm{argmin} \;  \mathcal  F(x) , \qquad x\in \mathbb R^d, 
	\end{equation}
	where $\mathcal F(\cdot):\mathbb R^d \to \mathbb R$ is a given non-convex cost by iterative, gradient--free particle methods. The literature on numerical methods for the solution of problem \eqref{eq:optim} is very rich, see {for example the textbooks \cite{grippo2023introduction,MR1678201,MR2244940}  and references therein.}
	We focus here on metaheuristic algorithms   \cite{KGV} {and \cite[Section 3.2]{blum2003metaheuristics}} that recently  have gained renewed interest, see \cite{P24,borghi2024kineticmodelsoptimizationunified}. This is due to the understanding of those methods as interacting particle systems where now tools of statistical physics and kinetic theory can be applied. These mathematical tools have, in particular, for so--called consensus-based methods, led to a different  point of view of such methods and they have proven fruitful in the subsequent analysis, see e.g. \cite{CCTT,MR4793478,PTTM,MR3629153,Ha_CBO20,Albi_23,TW_2020,CJLZ,PFZ_25} and references therein. In this paper, we also plan to exploit methods from kinetic theory focusing on a class of particle methods that is a continuous variant of the simulated annealing (SA) method \cite{GH}.  Since their introduction by Kirkpatrick, Gelatt, and Vecchi in 1983 \cite{KGV}, { the original method has gained interest in the mathematical community, see e.g. \cite{MR983115,MR1188544,nitanda2022convex} and the references therein.}  The SA  method is based on the evolution of an initial set of particles whose dynamics is obtained through an acceptance-rejection strategy {that} depends on a (control) parameter, {called temperature, which is decreasing over the number of iterations. In the last decades, a significant effort has been devoted to understanding a temperature schedule for efficient searching, see e.g. the pioneering works on adaptive simulated annealing (ASA) in \cite{ingber89,ingber93} and on simulated tempering \cite{parisi92} where the temperature becomes a dynamic variable. }    Our viewpoint is inspired by the recent work \cite{P24} where the links between SA, a Boltzmann--type description, and the corresponding Fokker--Planck equation (after suitable scaling) have been explored. { Following \cite{P24}, the main differences compared to the original work, is the introduction of a continuous time-scale in the SA method as well as a probabilistic interpretation of the particle dynamics. Those modifications lead to the dynamics described in Section \ref{sect:kinet_ann}.  } The relation between the continuous SA variant, see equation \eqref{eq:Xipart}, and the Langevin dynamics has been extensively  studied, see in this direction \cite{chizat2022mean,chak2023generalized}. Also, the dynamics of SA have been understood as an interacting particle system where the states correspond to approximations to the global minimiser $x^*.$  {The evolution of the temperature in the SA algorithm is of paramount importance for convergence \cite{MR942621} and is typically prescribed through a decay $O(1/\log(t))$, as discussed in \cite{P24}, in such a way that the the overall dynamics converge asymptotically toward an atomic measure centred in $x^*$. We summarise these results in Section~\ref{sect:kinet_ann}. }
	\par 
	{ In this work, we aim to improve these results by optimising the temperature dynamics to guide exponential convergence toward a quasi-equilibrium state, see Lemma~3.1. In a certain regime specified below, we are able to guide the equilibration rate of the Shannon entropy with an upper bound on the expected rate. To this end, we exploited recent results on Fokker-Planck-type equations for the formation of Gamma-type distributions and their equilibration rate, see \cite{ATZ_23,FPTZ22,Tosc}.
	A new entropic SA is proposed and consists of a coupled system of kinetic equations for the evolution of particles and temperature. At the microscopic level we consider the dynamics  \eqref{eq:trans}. Here, temperature is understood as an additional state space variable of the interacting particle system with its own parameterised dynamics. This interpretation allows to exploit the temperature to obtain improved performance of the SA method. The use of temperature as additional state variable allows to apply also tools from kinetic theory. Using this viewpoint, } a system of Boltzmann--type equations on the extended phase space can be derived, and, in the quasi-invariant scaling, we obtained a coupled system of mean-field equations \cite{MR1165528}. Those equations allow subsequently decay estimates on the distance to the solution $x^*$. Furthermore, the analysis of the arising Fokker--Planck equations allows to define  suitable choices for the parameter of temperature dynamics, such that { faster} decay rates can be established both theoretically and numerically, see Section \ref{sec:control}. In Section \ref{sec:numerics}, numerical results are presented that show the improved convergence rates.   
	
	\section{Kinetic simulated annealing}\label{sect:kinet_ann}

	For the solution of problem \eqref{eq:optim} we consider variants of the  SA algorithm \cite{KGV}. More precisely, we consider the  continuous version of the SA algorithm as proposed e.g. in \cite{GH}, where, starting from an initial sample $X(0)$ and an initial temperature $T(0)$ the SA process is described by a stochastic differential equation of the form 
	\begin{equation}
		\label{eq:Xipart}
		dX_i = -\nabla_x \mathcal F(X_i)dt + \sqrt{2T(t)}dW_i^t, 
	\end{equation}
	being $\{W_i\}_{i = 1,\dots,N}$ a set of independent Wiener processes. We now follow \cite{P24} and using classical methods of stochastic analysis, define the limit $N\to +\infty$ from \eqref{eq:Xipart}. This leads to the meanfield equation describing the evolution of the particle density $f(x,t)$, where $f(x,t)dx$ represents the proportion of particles in the volume $[x,x+dx)$ at time $t\ge0$. In detail, we get 
	\begin{equation}
		\label{eq:fmf}
		\partial_t f(x,t) = \nabla_x \cdot \left[ \nabla_x \mathcal F(x)f(x,t) + T(t)\nabla_x f(x,t) \right]. 
	\end{equation}
	It is worth to remark that the quasi-stationary state, i.e., the unique solution to the differential equation 
	\[
	\nabla_x \mathcal F(x)f(x,t) + T(t)\nabla_x f(x,t) = 0
	\] 
	corresponds to the density annihilating the flux of \eqref{eq:fmf}. It is given by the Boltzmann-Gibbs measure
	\begin{equation}
		\label{eq:gibbs}
		f^q_{\F}(x,t) = C(t) e^{-\F(x)/T(t)}, 
	\end{equation}
	being $C(t)>0$ a normalisation factor such that {$\int  f^q_{\F}(x,t) dx = 1$}. The factor $C(t)$ depends on time due to the dependence on the temperature $T.$  If $T(t)\ll 1$ we may expect that the $f^q_{\F}(x,t)$ is concentrated around the global minima of the cost functional $\mathcal F(x)$. { The index $q$ indicates that $f^q_{\F}$ is the quasi-stationary distribution. } 
	
	\begin{remark}
		We observe that \eqref{eq:fmf} can be equivalently rewritten in Landau form (if $f^q_\F>0$)
		\begin{equation}
			\label{eq:fmf2}
			\partial_t f(x,t) = \nabla_x\cdot \left[ T(t)f(x,t) \nabla_x  \log \dfrac{f(x,t)}{f^q_\F(x,t)}\right], 
		\end{equation}
		since the Gibb's measure is the unique solution of the following differential equation
		\[
		-\dfrac{\nabla_x f^q_\F(x,t)}{f^q_\F(x,t)} = \dfrac{\nabla_x \mathcal F(x)}{T(t)}. 
		\]
		{The equivalence of equation \eqref{eq:fmf2} to equation \eqref{eq:fmf} follows by direct computation, see e.g. \cite{FPTZ_17,PZ18}. 
		}
	\end{remark}
	
	In \cite{P24} a kinetic version of the simulated annealing algorithm has been presented.  We follow this presentation. Therein, a distribution  $f: \mathbb R^d \times \mathbb R_+ \to \mathbb R_+$ of particles with the position in $[x,x+dx)$ at time $t\ge0$ has been introduced. The evolution of the  distribution $f=f(x,t)$ follows the integro-differential equation
	\begin{equation}
		\label{eq:P24_boltzmann}
		\partial_t f(x,t) = \int_{\mathbb R^d} (B_\F(x^\prime \to x)f(x^\prime,t) - B_\F(x \to x^\prime)f(x,t))p(\xi)d\xi, 
	\end{equation}
	where the particles' updates are defined as follows
	\begin{equation}\label{001} 
	x^\prime = x + \sigma(t)\xi, \qquad \xi \sim p(\xi),
	\end{equation}
	and the following cross-section has been considered
	\begin{equation}
		\label{eq:kernel}
		B_\F(x \to x^\prime) = \min\left\{ 1, \dfrac{f^q_{\F}(x^\prime)}{f^q_\F(x)} \right\}. 
	\end{equation}
	The introduced dynamics may be written in weak form. Let $\varphi(\cdot)$ be a smooth test function, then the evolution of observable quantities are described by the
	\begin{equation}
		\label{eq:P24_weak}
		\begin{split}
			&\dfrac{d}{dt} \int_{\mathbb R^d} f(x,t)dx =\\
			&\quad \dfrac{1}{2} \int_{\mathbb R^d}\int_{\mathbb R^d} (\varphi(x^\prime) - \varphi(x)) (B_\F(x \to x^\prime)f(x,t) - B_\F(x^\prime\to x) f(x^\prime,t))p(\xi)dxd\xi. 
		\end{split}
	\end{equation}
	In strong form, \eqref{eq:P24_weak} reads
	\begin{equation}
		\label{eq:P24_strong}
		\partial_t f(x,t) = \mathbb E_{\xi}\left[ B_\F(x^\prime\to x)f(x^\prime,t) - B_\F(x\to x^\prime)f(x,t)\right]:= \mathcal L_\F(f)
	\end{equation}
	This kinetic model is coherent with the mean-field equation \eqref{eq:fmf} under a suitable quasi-invariant scaling of the parameters, see \cite{PT}. Thanks to this analogy, the temperature dynamics associated to the Gibbs measure $f^q_\F(x,t)$ and defining the cooling process can be investigated. Through direct computation of the  Shannon entropy, it has been shown that the decay rate of the temperature should be $o(T^2(t))$.  A typical examples is as introduced above by $T(t) \approx 1/\log(t)$.

	\subsection{Non-homogeneous temperature dynamics}\label{sec:dyn}
	
	To overcome the limitation provided by a slow temperature decay, we consider 
	non-homogeneous and non-autonomous temperature dynamics, whose decay rate is influenced by the particles' dynamics. We consider the distribution of particles $f(x,t)$ at position  $[x,x+dx) \subset R^d$ and at time $t\ge0$. Further, we introduce  the    distribution of temperature $g=g(T,t): \mathbb R_+ \times \mathbb R_+ \to \mathbb R_+$ describing the fraction of particles characterized by temperature in the interval $[T,T+dT) \subset \mathbb R_+$ at time $t\ge0$.  The evolution of such system of particles is defined in terms of the following microscopic scheme
	\begin{equation}
		\label{eq:trans}
		\begin{split}
			x^\prime  &= x + \mathcal D[g](t)\xi, \\
			T^\prime  &= T -\lambda[f] T + \kappa(T)\eta,
		\end{split}
	\end{equation}
	where  $\xi \sim p_1(\xi)$ is a symmetric random variable such that $\left\langle\xi \right\rangle = 0$ with identity covariance matrix $\Sigma =  I_d$.  Furthermore,  $\eta \sim p_2(\eta)$ is a  random variable independent on $\xi$ such that $\mathbb E[\eta] = 0$, $\mathbb E_\eta[\eta^2] = \sigma^2 <+\infty$.  We also assume that the third order moments of the introduced random variables are bounded, i.e.  $\left\langle\xi^3\right\rangle <+\infty$, $\left\langle \eta^3 \right\rangle<+\infty$. 
	
	In \eqref{eq:trans} the  strength of the diffusion is tuned by the positive operator $\mathcal D[g](t)\ge0$ depending on the (local) temperature distribution $g(T,t)$ and such that {$\mathcal D[\delta(T-0)] = 0$, where $\delta(T-\tau)$ is the Dirac delta distribution centred in $\tau$}. A prototype example, used also in the later analysis, is the $k$th moment of $g,$ i.e.,  
	\begin{equation}\label{operator D}
		\mathcal{D}[g](t) =  \int_{\mathbb R_+}T^k g(T,t) dT. 
	\end{equation}
	
	Finally,  the operator $\lambda = \lambda[f] \geq 0$ is a (feedback) control  that may depend on the distribution of the particles' positions  $f(x,t)$ and whose form, which will be presented later on, guarantee temperature cooling. It will be chosen in Section \ref{sec:control} to ensure exponential convergence. 
	
	From \eqref{eq:trans} the following relations are obtained 
	\[
	\mathbb E_{\xi} \left[x^\prime - x \right] = 0, \qquad \mathbb E_{\eta} \left[ T^\prime - T\right] = -\lambda T \le 0. 
	\]
	Therefore, the introduced microscopic transitions  leave, on average, the particle positions unchanged, whereas the expected  temperature of the system decays with a rate $\lambda \in [0,1]$. The microscopic scheme determining the temperature update $T^\prime$ should guarantee the positivity of the post-transition temperature. Depending on the choice of $\kappa(\cdot)$ we need to define suitable bounds on the support of $\eta$. For example,   in the case $\kappa(T) = T$ we obtain 
	\[
	T^\prime = T(1-\lambda + \eta)\ge0.
	\]
	and the positivity of $T^\prime$ is guaranteed {for any $\eta$ with bounded support in the interval} 
	\begin{equation}
		\label{eq:eta_bound}
		| \eta| \le 1-\lambda,
	\end{equation}
	being $\eta$ a zero mean random variable.
	We observe that, provided $\lambda \in [0,1]$, the condition  \eqref{eq:eta_bound} { is fulfilled. In case of a more general choice for the weight $\kappa(T)=T^p$, $ 0 \leq p<1$,  the positivity of $T'$ follows by considering a modified transition }  
	\begin{equation}
	\label{eq:T_dyn_mod}
	T^\prime = (1-\lambda) T  + T^p \chi\left(T \ge (1-p) \theta\right) \eta, 
	\end{equation}
	where $0<\theta<1$ and $\chi(A)$ is the characteristic function of the set $A \subseteq \mathbb R_+$, {see also the derivation proposed in \cite{BMTZ_25}}.  { The characteristic function is added to reduce the noise to zero for $T$ sufficiently small. } 	We may observe that, if $\kappa \not = 0,$  the { non-negativity of $T'$ is guaranteed provided that $\eta$ fulfils the uniform} bound
		\begin{equation}\label{eq:eta_bound_mod}
	|\eta| \le (1-\lambda)((1-p)\theta)^{1-p}.
	\end{equation}

	The  microscopic scheme \eqref{eq:trans} allows to investigate the coupled evolution of the distributions $(f,g)$ given in weak form by the equations 
	\begin{equation}
		\label{eq:kinetic}
		\begin{split}
			&\dfrac{d}{dt}\int_{\mathbb R^d}f(x,t) \varphi(x)dx \\
			&\quad = \dfrac{1}{2}\mathbb E_\xi\left[\int_{\mathbb R^d} (\varphi(x^\prime)-\varphi(x)) (B^g_\mathcal{F}(x\to x^\prime)f(x,t) - B^g_\mathcal{F}(x^\prime\to x)f(x^\prime,t))dx\right] \\
			&\dfrac{d}{dt}\int_{\mathbb R_+} g(T,t) \psi(T)dT = \dfrac{1}{\nu}\mathbb E_\eta\left[ \int_{\mathbb R_+}(\psi(T^\prime)-\psi(T))g(T,t) dT\right]. 
		\end{split}
	\end{equation}
	Here, $\varphi(\cdot)$, $\psi(\cdot)$ are smooth test functions and $\nu>0$ is related to the frequency of interactions characterising the temperature dynamics. In \eqref{eq:kinetic} we have introduced the cross--section  $B^g_{\mathcal F}(x\to x^\prime)$ which shares the structure defined in \eqref{eq:kernel} where now the Boltzmann-Gibbs distribution $f^q_\F(x,t)$ depends on the temperature dynamics through $\mathcal D[g]>0$ as follows
	\begin{equation} \label{eq:BG}
		f_\F^q(x,t) = C(t) e^{-\frac{\F(x)}{\mathcal D[g](t)}}, 
	\end{equation}
	being $C(t) >0$ a normalisation factor, {such that $\int_{\mathbb R^d}f_\F^q(x,t)dx = 1$.  }
	We remark that, at variance with the classical approach, the quasi-equilibrium distribution $f^q_\F$ depends now on $g(T,t)$ through the operator $\mathcal{D}$ given by equation \eqref{operator D}. 
	
	{The system \eqref{eq:kinetic} can be equivalently formulated in the following strong form} 
	\begin{equation}
		\label{eq:BG_strong}
		\begin{split}
			\partial f(x,t) &= \mathbb E_{\xi} \left[ B^g_\F(x \to x^\prime)f(x,t) - B_\F^q(x^\prime\to x)f(x^\prime,t)\right]:= \mathcal L_\F(f,g) \\
			\partial_t g(T,t) &= \dfrac{1}{\nu} \left(\mathbb E_{\eta} \left[ \int_{\mathbb R_+} \dfrac{1}{|{}^\prime J_f|}g({}^\prime T,t)dT\right] - g(T,t)\right):= \mathcal J(g,f).
		\end{split}
	\end{equation}
	being $|{}^\prime J_f|$ the determinant of the Jacobian matrix of the transformation ${}^\prime T \to T$ for a given $\lambda[f]>0$ in \eqref{eq:trans}, {being ${}^\prime T$ the pre-interaction temperature}. {We refer to \cite{pareschi2013interacting} for further details of the derivation of Fokker-Planck-type equations from kinetic dynamics \eqref{eq:BG_strong}.   }
	
	\subsection{The mean-field limit}
	In this section, we show that for a certain range
of the parameters characterizing the kinetic system of equations \eqref{eq:BG_strong}, usually referred to as the grazing limit or quasi-invariant collision regime, it is possible to derive a model of Fokker-Planck-type \cite{PT,villani98,bobylev00,toscani98} which facilitates the analytical study of the emerging equilibria. In the following we are interested in computing the quasi-stationary profiles emerging from the coupled system of  kinetic equations  \eqref{eq:kinetic} in suitable regimes of parameters. To this end, we consider the quasi-invariant scaling of the parameters as
	\begin{equation}
		\label{eq:scaling}
		t\to t/\epsilon, \qquad \sigma\to \sqrt{\epsilon}, \qquad \lambda \to \epsilon \lambda, \qquad \theta \to \epsilon
	\end{equation}
	with $\epsilon>0$. Following \cite{PT}, if $\epsilon \ll 1$, the transition schemes \eqref{eq:trans} is quasi-invariant, and we may expand the differences for  $\epsilon \to 0^+$
	\begin{equation}
		\label{eq:expansion}
		\begin{split}
			&\varphi(x^\prime) - \varphi(x) = (x^\prime-x) \cdot \nabla_x \varphi(x) + \dfrac{1}{2} \sum_{i,j=1}^d(x_i^\prime-x_i)\cdot(x_j^\prime-x_j) \dfrac{\partial^2 \varphi(x)}{\partial x_i\partial x_j} + O(\epsilon^{3/2}) \\
			&\psi(T^\prime) - \psi(T) = (T^\prime - T) \psi'(T) + \dfrac{1}{2}(T^\prime-T)^2 \psi''(T) + \dfrac{1}{6}(T^\prime-T)^3 \psi'''(\bar T).
		\end{split}
	\end{equation}
	Here, $\bar T  \in (\min \{T,T^\prime\}, \max\{T,T^\prime\})$ and the quasi-equilibrium density $f^q_\F(x^\prime)$ fulfills 
	\[
	\begin{split}
		f^q_\F(x^\prime) - f^q_\F(x) =& (x^\prime-x)\cdot \nabla_x f_\F^q(x) + O(\epsilon) =\\
		& - \dfrac{1}{D[g](t)}(x^\prime-x)\cdot\nabla_x \F(x) f_\F^q(x) + O(\epsilon),
	\end{split}\]
	As shown in \cite{P24} the  Boltzmann-type equation in \eqref{eq:kinetic} converges in a weak sense to the solution of the meanfield equation
	\begin{equation}\label{eq:syst}
		\begin{split}
			\partial_t f(x,t) &= \nabla_x \cdot \left[ \nabla_x \mathcal F(x) f(x,t) +\mathcal D[g](t)\nabla_x f(x,t)\right]. 
		\end{split}
	\end{equation}
	The quasi equilibrium density of the obtained Fokker-Planck equation corresponds to the Boltzmann-Gibbs distribution with temperature defined by the operator $\mathcal D[g]\ge0$. More precisely, we have 
	\begin{equation}\label{eq:equi}
		f^q_\F(x,t) = C(t) e^{-\F(x)/\mathcal D[g](t)}. 
	\end{equation}
	Note that since $g$ depends on time $t,$ also the normalization factor $C(t)$ is time dependent.  Hence, the  state $f^q_\F$ is only a quasi- equilibrium state. 
	
	To determine the equation for the evolution of the temperature $g$ we formally compute the expansion \eqref{eq:expansion} in the second equation of \eqref{eq:kinetic}. Under the scaling \eqref{eq:scaling},  the evolution of the density $g^\epsilon(T,\tau) = g(T,t/\epsilon)$ reads
	\[
	\begin{split}
		\nu\dfrac{d}{d\tau} \int_{\mathbb R_+}\psi(T)g^\epsilon(T,\tau)dT =&\dfrac{1}{\epsilon} \int_{\mathbb R_+} \mathbb (-\epsilon\lambda T)\psi'(T)g^\epsilon(T,\tau)dT +\\
		&  \dfrac{1}{2\epsilon }\int_{\mathbb R_+}(\epsilon^2\lambda^2[f]T^2 +\epsilon T^{2p} \sigma^2 ) \psi''(T)g^\epsilon(T,\tau)dT+ \\
		&R_\psi(g^\epsilon) + \bar R_\psi(g^\epsilon), 
	\end{split}
	\]
	being 
	\[
	\begin{split}
		R_\psi(g^\epsilon) =& \dfrac{1}{6\epsilon} \int_{\mathbb R_+} \mathbb E_\eta[-\lambda T + T^p\eta]^3 \psi'''(T)g^\epsilon(T,\tau)dT
	\end{split}
	\]
	and
	\[
	\begin{split}
		\bar R_\psi(g^\epsilon) =&\dfrac{1}{\epsilon} \int_{\mathbb R_+}\chi(0\le T \le (1-p)\epsilon)T^{2p} \mathbb E_\eta[\eta^2]\psi''(\bar T)g^\epsilon(T,\tau) dT +\\
		&\dfrac{1}{\epsilon} \int_{\mathbb R_+} \chi(0\le T \le (1-p)\epsilon) T^{3p} \mathbb E_\eta[3 \eta^2 + \eta^3]\psi'''(\bar T)g^\epsilon(T,\tau)dT. 
	\end{split}
	\]
	Assuming that the random variable $\eta$ has bounded third order moments, $\mathbb E_\eta[\eta^3]<+\infty$, we  obtain
	\[
	|R_\psi(g^\epsilon)|\approx \epsilon   + \epsilon^2 + \sqrt{\epsilon} \qquad |\bar R_\psi(g^\epsilon)|\approx \epsilon^{1+2p}. 
	\]
	Hence, in the limit $\epsilon\to 0^+$, all the reminder terms vanish and the distribution $g^\epsilon(T,\tau)$ converges to a limiting distribution $g(T,\tau)$ solution to 
	\[
	\nu\dfrac{d}{d\tau} \int_{\mathbb R_+}\psi(T)g(T,\tau)d\tau = \int_{\mathbb R_+} (-\lambda T) \psi'(T)g(T,\tau)dT + \dfrac{\sigma^2}{2} \int_{\mathbb R_+}T^{2p}\psi''(T)g(T,\tau)dT. 
	\]
	The strong form is  the Fokker-Planck equation for the evolution of the temperature 
	\begin{equation}\label{eq:FP_temp}
		\begin{split}
			\partial_t g(T,t) &= \dfrac{1}{\nu}\partial_T \left[ \lambda[f](t)Tg(T,t) +\dfrac{\sigma^2 }{2}\partial_T (T^{2p}g(T,t)) \right], 
		\end{split}
	\end{equation}
	with no-flux boundary conditions for the  equation of the temperature
	\begin{equation}
	\label{eq:boundary_terms}\begin{split}
		\lambda[f](t)Tg(T,t) +\dfrac{\sigma^2}{2} \partial_T (T^{2p}g(T,t)) \Big|_{T=0} = 0\\
		T^{2p}g(T,t)\Big|_{T=0 }= 0.  
	\end{split}
	\end{equation}
	
	Therefore, we are interested in the  evolution of the coupled  system of Fokker-Planck equations 
	\begin{equation}
	\label{eq:system_FP}
	\begin{split}
	\partial_t f(x,t) &= \nabla_x \cdot \left[\nabla_x \F(x)f(x,t) + \mathcal D[g](t)\nabla_x f(x,t) \right] \\
	\partial_t g(T,t) &= \dfrac{1}{\nu}\partial_T \left[ \lambda[f](t)Tg(T,t) +\dfrac{\sigma^2 }{2}\partial_T (T^{2p}g(T,t)) \right], 
	\end{split}
	\end{equation}
	with no-flux boundary conditions \eqref{eq:boundary_terms}. 
	
		We summarize the computations of this section as follows: under the quasi-invariant scaling, the evolution of the particle probability density $f(x,t)$ follows the mean-field equation \eqref{eq:syst}, where the temperature evolves in time according to a functional of the temperature distribution $g(T,t)$. Consequently, the dynamics of the temperature distribution is governed by the Fokker–Planck equation \eqref{eq:FP_temp}. The quasi-equilibrium distribution $f^q_\F$ depending on $g(T,t)$ and $\F$ is given by equation \eqref{eq:equi}. In the forthcoming section, we study the convergence of $f$ towards the equilibrium distribution $f^q_\F.$

	\begin{remark}
	By equation \eqref{eq:FP_temp} and in the limit $\nu\to 0^+$ we  compute the quasi-equilibrium state as the unique solution to the following differential equation 
	\[
	(\lambda[f](t)T +p\sigma^2 T^{2p-1})g^q(T,t) + \dfrac{\sigma^2}{2}T^{2p}\partial_Tg^q(T,t) = 0, 
	\]
	or equivalently
	\[
	\dfrac{\partial_T g^q(T,t)}{g^q(T,t)} = -\dfrac{2}{\sigma^2}\dfrac{\lambda[f]T + p\sigma^2 T^{2p-1}}{T^{2p}}. 
	\]
	From which we obtain  the quasi--equilibrium state 
	\begin{equation}\label{eq:equilibrium}
		g^q(T,t) = \dfrac{2(1-p)}{\Gamma\left(\frac{1-2p}{2(1-p)}\right)\left(\frac{\sigma^2(1-p)}{\lambda[f](t)}\right)^{\frac{1-2p}{2(1-p)}}}\exp\left\{ - \dfrac{\lambda}{\sigma^2(1-p)}T^{2(1-p)} \right\}T^{-2p} 
	\end{equation}
	with  $\Gamma(\cdot)$ being the Gamma function,  { i.e., for $z\in \mathbb{C}$ and $Re(z)>0,$ it is given by $\Gamma(z)=\int_0^\infty t^{z-1} \exp(-t)  dt.$ }	 The  quasi-equilibrium distribution is a generalized Gamma density \cite{stacy62} which is integrable for any $0<p<1/2$. 
	\end{remark}
	\par

	\subsection{Convergence to equilibrium}\label{sec:2.3}
	Similar to \cite{MR3352763}, we define the relative Shannon entropy of the distribution of $f$ with respect to a probability distribution $h$ as follows
	\begin{equation}
		\label{eq:shannon}
		H(f|h)(t) = \int_{\mathbb R^d} f(x)\log \dfrac{f(x)}{h(x)}dx. 
	\end{equation}
	The study of entropy decay has its roots in the classical kinetic theory of rarefied gases. In traditional modeling approaches, it is used to identify the equilibrium density of the system.. In the case of the classical Fokker-Planck equation, the proof of exponential convergence exploits the logarithmic Sobolev inequality \cite{Toscani99}. On the other hand, the evolution of the relative entropy of the solution density of the Fokker-Planck equation can be used to obtain a dynamical proof of the logarithmic Sobolev inequality under suitable assumptions \cite{OTTO2000361}. The definition of sharp rates strongly depends  on the structure of the equilibrium distribution as highlighted in \cite{FPTZ19,FPTZ22}. More recently, the study of the convergence to a quasi-equilibrium distribution has been considered as a follow-up question in the modelling of many-agent systems undergoing nonlocal interactions \cite{ATZ_23}.  
	
	In more detail, we are interested in the convergence of the particle density toward the quasi-equilibrium distribution given by the Boltzmann-Gibbs distribution $f^q_\F(x,t)$ defined in equation \eqref{eq:equi}. We may compute the evolution of the Shannon entropy \eqref{eq:shannon} along the solution of the Fokker-Planck equation \eqref{eq:system_FP}. The decay of the relative entropy is computed as follows
	\begin{equation}
		\label{eq:ddtH1}
		\dfrac{d}{dt} H(f|f^q_\F)(t) = \underbrace{\int_{\mathbb R^d}\partial_t f(x,t) \left( \log \dfrac{f(x,t)}{f^q_\F(x,t)}+1 \right)dx}_{I} - \underbrace{\int_{\mathbb R^d}f(x,t) \partial_t \log f^q_\F(x,t)dx}_{II}.
	\end{equation}
	We consider the terms  $(I)$ and $(II)$, separately. Thanks to the  conservation of mass we obtain for $(I):$
	\begin{equation}
		\label{eq:fisher}
		\begin{split}
			(I) &= \int_{\mathbb R^d}\partial_t f(x,t) \left( \log \dfrac{f(x,t)}{f^q_\F(x,t)}+1 \right)dx 
			= \int_{\mathbb R^d}\partial_t f(x,t)  \log \dfrac{f(x,t)}{f^q_\F(x,t)}dx \\
			&\quad= -\int_{\mathbb R^d}\mathcal D[g](t)f(x,t) \nabla_x \log \dfrac{f(x,t)}{f^q_\F(x,t)} \cdot  \nabla_x\log \dfrac{f(x,t)}{f^q_\F(x,t)} dx\\
			&\quad =: -I_H(f|f^q_\F)(t).
		\end{split}
	\end{equation}
	The term $I_H(f|f^q_\F)\ge0$ is the classical Fisher information \cite{CT07,Toscani99} where now the Maxwellian is the Boltzmann-Gibbs distribution. Furthermore,  the term $(II)$ in \eqref{eq:ddtH1}  may be rewritten   using 
	\[
	\log f^q_\F(x,t) = \log C(t) - \dfrac{\F(x)}{\mathcal D[g](t)},
	\]
	and 
	\[
	\partial_t \log f^q_\F(x,t) = \dfrac{\dot{C}(t)}{C(t)} + \dfrac{\F(x)\dot{\mathcal D}[g](t)}{\mathcal D^2[g](t)}. 
	\]
	Hence,  \eqref{eq:ddtH1} reads
	\begin{equation}\begin{split}
			\label{eq:ddtH2}
			\dfrac{d}{dt}H(f|f^q_\F) &= -I_H(f|f^q_\F) - \int_{\mathbb R^d}f(x,t) \left(\dfrac{\dot{C}}{C} + \dfrac{\F(x) \dot{\mathcal D}[g](t)}{\mathcal D^2[g](t)} \right)dx \\
			&= -I_H(f|f^q_\F) - \dfrac{\dot{C}}{C} - \dfrac{\dot{\mathcal D}[g](t)}{\mathcal D^2[g](t)} \int_{\mathbb R^d}\F(x) f(x,t)dx. 
		\end{split}
	\end{equation}
	The differential {equation} \eqref{eq:ddtH2} can be further simplified. Since the Gibb's measure $f^q_\F(x,t) = C(t) e^{-\F(x)/\mathcal D[g](t)}$ has unitary mass, we get
	\[
	0 = \int_{\mathbb R^d}\dot{C} e^{\F(x)/\mathcal D[g](t)} + \dfrac{C \dot{\mathcal D}[g]\F(x)}{\mathcal D^2[g](t)}e^{-\F(x)/\mathcal D[g](t)}dx, 
	\]
	that yields
	\[
	\dfrac{\dot{C}}{C} = -\dfrac{\dot{\mathcal D}[g](t)}{\mathcal D^2[g](t)}\int_{\mathbb R^d}\F(x) f^q_\F(x,t)dx. 
	\]
	Therefore, for any $f(x,t)$ solution to \eqref{eq:syst}, the evolution of the Shannon entropy is given by 
	\begin{equation}
		\label{eq:evo_H_part}
		\dfrac{d}{dt}H(f|f^q_\F)(t) = -I_H(f|f^q_\F) - \dfrac{\dot{\mathcal D}[g](t)}{\mathcal D^2[g](t)} \int_{\mathbb R^d}\F(x) (f(x,t)-f^q_\F(x,t))dx. 
	\end{equation}
	We  observe that, since the sign of the integrand  is unknown, we cannot conclude that the entropy  is dissipated in time for an arbitrary  distribution $f$. 
	
	\begin{remark}\label{rem1}
		In \cite{P24}  the case $\mathcal D[g](t) = T(t)$ has been considered.  If $ T = 1/\log(t)$ we obtain that $\dot T/T^2\approx 1/t$ and the second term annihilates for $t\to +\infty$. Therefore, we can expect to converge asymptotically towards the correct equilibrium. 
	\end{remark}

	%%%%%%%%%%%%%%%%
	
	\section{Choice of the  feedback control  $\lambda$ }\label{sec:control}
	
	This section is devoted to the choice of the parameter $\lambda$ in order to possibly obtain  a speed-up in the convergence to equilibrium of the continuous SA algorithm compared to the one observed in Remark \ref{rem1}. Since we aim to evaluate $\lambda=\lambda[f]$ within the particle dynamics, we assume that $\lambda$ is a closed--loop feedback on the state of the system, i.e., $f=f(x,t)$.  An open loop control might not be feasible from a computational point of view. The final form is presented in Lemma \ref{lem1}.  
	\par 
	In order to define a suitable closed--loop control for $\lambda$ we consider the relative entropy $H$ as a Lyapunov function. Equation \eqref{eq:evo_H_part} yields an equality for $H$ where the feedback enters in the operator $\mathcal{D}.$ First, the dependence of $\mathcal{D}$ and $\lambda$ is made explicit using the linearity of the operator $\mathcal D[\cdot]$
	\[
	\partial_t \mathcal D[g](t)=\mathcal D[\partial_t g]. 
	\]
	and  the particular form of the operator  \eqref{operator D}, i.e., for $k>0$ we define  
	\[
	\mathcal D[g](t)  = \int_{\mathbb R_+}T^k g(T,t)dT =: m_k(t) \geq 0.
	\]
	Using this choice, the following dependence of $\mathcal{D}$, resp. $m_k,$  is explicit in $\lambda:$ 
	\begin{equation}
		\label{eq:Devo}
		\begin{split}
			\dfrac{d}{dt}\mathcal D[g](t)  = \dfrac{d}{dt} m_k(t) =& -k \lambda[f]m_k(t)  + \dfrac{k(k-1)\sigma^2}{2}m_{k-2(p-1)}(t). 
		\end{split}
	\end{equation}
	Following the discussions in Section \ref{sec:dyn},   the operator $\mathcal D$ guarantees the dissipation of the temperature if 
	\begin{equation}
		\label{eq:lambda_cond_q}
		\lambda[f](t) > \dfrac{\sigma^2(k-1)}{2} \dfrac{m_{k-2(p-1)}(t)}{m_k(t)}
	\end{equation}
	and we obtain the following equality for the Lyapunov function $H:$ 
	\begin{equation}\label{eq:HF}
	\begin{split}
		\dfrac{d}{dt}H(f|f^q_\F)(t) =& - I_H(f|f^q_\F) \\
		&  -\dfrac{ k(k-1)\frac{\sigma^2}{2}m_{k-2(p-1)}(t) -k\lambda[f]m_k(t)  }{m_k^2(t)} \int_{\mathbb R^d}\F(x) (f(x,t) - f^q_\F(x,t))dx.
		\end{split}
	\end{equation}
	Note that,  in general, there is no information on the dependency of the sequence of moments  $\{m_k\}_{k>0}$ on each other.
	We  derive a  feedback control  in the relevant case $k=1$, where the condition  \eqref{eq:lambda_cond_q} reduces to $\lambda[f](t)>0.$ 
	
	\begin{lemma}[Case $k=1$]\label{lem1}
		Consider the relative entropy $H$ defined in equation \eqref{eq:shannon} for a regular solution $f$ to equation \eqref{eq:syst} with arbitrary initial distributions $f_0 \geq 0$, { $\int_{\mathbb R^d} f_0(x) dx = 1.$} Furthermore, let $\mathcal{D}$ be given by equation \eqref{operator D} with $k=1$. Assume  $g$ to be a regular solution to equation \eqref{eq:FP_temp} with arbitrary initial data $g_0 \geq 0 $ with { $\int_{\mathbb R^+} g_0(T) dT = 1.$ }  Furthermore, denote by $f^q_\F$ the quasi-equilibrium solution (depending on $\F$ and $g$) given by equation \eqref{eq:BG}. Then for any  $0<\alpha<\frac{ \sqrt{2} \| \F \|_\infty }{ m_1(0) \sqrt{H(0)} } < \infty$ for {$H(0)=H(f_0, f^q_\F(\cdot,0) )$} the following statement holds:  
		
		\noindent If for $t\geq0,$  we have 
		$\int_{\mathbb R^d}\F(x) ( f^q_\F(x,t) - f(x,t))dx > 0$, then for 
		\begin{align}\label{lem:lambda}
			\lambda[f](t) = \alpha \frac{ m_1(t) \sqrt{H(0)} }{ \sqrt{2} \| \F \|_\infty }
		\end{align}
		we have
		\begin{align*}
			\frac{d}{dt} H(f | f^q_\F) \leq - \alpha H(f | f^q_\F).
		\end{align*}
	\end{lemma}
	Some remarks are in order.
	\begin{remark}\label{rem:lem1}
		The condition on $\alpha$ ensures that $\lambda$ is bounded by one which is necessary to obtain a well-posed particle dynamics, see Section \ref{sect:kinet_ann} and equation \eqref{eq:trans}. The decay of $H$ ensures at least exponential convergence on all time intervals where the condition  
		$\int_{\mathbb R^d}\F(x) ( f^q_\F(x,t) - f(x,t))dx > 0$ holds true. We will later investigate also numerically and present examples where this condition is fulfilled. The feedback also depends on $g,$ but only through the moment $m_1(t).$ Since the moment decays \eqref{eq:Devo} (at least for $k=1$), this condition can be weakened.  The rate $\alpha$ depends in particular on $H(0)$ which measures the difference between the initial $f_0$ and the Gibbs distribution (at initial temperature). The closer the initial distribution $f_0$ is towards $f^q_\F,$ the faster the convergence to equilibrium. {Since condition $\int_{\mathbb R^d}\F(x) ( f^q_\F(x,t) - f(x,t))dx $ does not hold true} for arbitrary choices of $f_0$, the closed loop feedback $\lambda[f]$ needs to be defined also in the case  $\int_{\mathbb R^d}\F(x) ( f^q_\F(x,t) - f(x,t))dx \leq 0.$ 
		\par 
		{ For $k=1,$ the equation for $m_1$ reads in all cases
			\begin{align*}
				\dot m_1(t) = -\lambda[f](t) m_1(t) 
			\end{align*}
			In the hypothesis that $\lambda[f]$ is independent of $f$ we can fix 
			\[
			\lambda[f](t)  = \dfrac{1}{t+2}\log^{-1}(t+2).
			\]
			 Then, the previous equation can be explicitly integrated  and we obtain  $m_1(t)=\frac{T_0}{\log(t+2)}$ for an initial datum $m_1(0)=T_0.$ This formula is precisely as suggested in \cite{P24}.  
		} 
	\end{remark} 

	\noindent {\bf Proof of Lemma \ref{lem1}}. 
	{We introduce $\mathcal I_\F(t) := \int_{\mathbb R^d}\F(x)(f^q(x,t)-f(x,t))dx$.}  Then, the  starting point is the equality \eqref{eq:HF} that, 
	in the case $k=1,$ simplifies to 
	\begin{align*}
		\dfrac{d}{dt}H(f|f^q_\F)(t) &= - I_H(f|f^q_\F)  -
		\frac{ \lambda[f]  }{m_1(t) }\int_{\mathbb R^d}\F(x) (f^q_\F(x,t) - f(x,t))dx. 
	\end{align*} 
	Since $I_H \geq 0,$ $\lambda, m_1(t) \geq0$ and by  assumption also $\int_{\mathbb R^d}\F(x) (f^q_\F(x,t) - f(x,t))dx\geq 0,$ 
	we obtain \[ H(f|f^q_\F)(t)\leq H(f|f^q_\F)(0)\] for all $t\geq0.$  Using the explicit form \eqref{lem:lambda} of $\lambda$ we furthermore have  	
	\begin{align}\label{lem:H}
		\dfrac{d}{dt}H(f|f^q_\F)(t) \leq  -  \alpha \frac{   \sqrt{H(0)} }{ \sqrt{2} \| \F \|_\infty }\int_{\mathbb R^d}\F(x) (f^q_\F(x,t) - f(x,t))dx.
	\end{align} 
	Due to the Czisar-Kullback-Pinsker inequality 
	\begin{align*}
		2H(f|f^q)(t)  \geq \|f-f^q_\F \|_{L^1}^2, 
	\end{align*}
	and hence due to the Cauchy-Schwarz inequality 
	\begin{align*}	 \sqrt{ 2 H(f|f^q)(t) }   \geq \|f-f^q_\F \|_{L^1}  \geq \frac{1}{ \| \F \|_\infty }  \int_{\mathbb R^d}\F(x) (f^q_\F(x,t) - f(x,t))dx.
	\end{align*}
	Since $H(f|f^q_\F)\geq0$, if $\mathcal I_\F(t)\ge0$, we get the monotone decay of the relative entropy 
	\[
	H(f|f^q_\F)(t) \leq H(f|f^q_\F)(0)
	\]
	 which enforces the inequality \eqref{lem:H}
	\begin{align}\label{lem:H2}
		\dfrac{d}{dt}H(f|f^q_\F)(t) \leq  -  \alpha \frac{  \sqrt{H(0)}  }{ \sqrt{2}  }  \sqrt{2 H(f|f^q_\F)(t)} \leq - \alpha H(f|f^q_\F)(t). 
	\end{align}  
	%%%%%%%%%%%%%%%%%%%%%%

	Summarizing, in the case $k=1,$ we propose the following closed-loop feedback depending  on the sign  of $\mathcal I_\F(t)$ 
	\begin{equation}
		\label{eq:lambda_kappa1}
		\lambda[f](t) = 
		\begin{cases}
			\alpha \dfrac{ m_1(t) \sqrt{H(0)} }{ \sqrt{2} \| \mathcal F\|_{\infty} } & \mathcal I_\F(t)\ge0 \\[5mm]
			\dfrac{1}{(t+2)\log(t+2)}& \mathcal I_\F(t)<0
		\end{cases}
	\end{equation}
		
	\noindent In the case of higher moments of $g,$ i.e., operators $\mathcal{D}$ of the type \eqref{operator D} with $k >1:$
	\begin{align*}
		\mathcal{D}[g](t) =  \int_{\mathbb R_+}T^k g(T,t) dT, 
	\end{align*}
	we proceed as follows.  Computing formally the evolution of the 
	moments $\mathcal{D}$ by integration \eqref{eq:FP_temp} leads to
	the dynamics \eqref{eq:Devo}. For the given system in the case $k>1,$  we do not expect to obtain a priori information on moments $m_k$ of the temperature $g$ for a general initial distribution $g(T,0) \in L^1(\mathbb R_+)$.  The condition \eqref{eq:lambda_cond_q} is required for the dissipation of the operator $\mathcal{D}.$ In the limit $\nu\to 0^+,$ we expect that $g\to g^q(T,t)$ given by the class of generalized Gamma densities defined in \eqref{eq:equilibrium} for any $p \in (0,\frac12)$. 
	
	From a kinetic point of view, this corresponds to the case when the temperature dynamics are faster than the particle dynamics. In more detail, in the limit $\nu \to 0^+$, the strong form of \eqref{eq:syst} for the probability density $f=f(x,t)$ reads 
	\begin{equation}\label{eq:fnu} 
		\partial_t f(x,t) = \nabla_x \cdot \left[ \nabla_x \F(x) f(x,t) + \mathcal D[g^q] \nabla_x f(x,t) \right],
	\end{equation}
	where $\mathcal{D}[g]$ has been replaced by $\mathcal{D}[g^q].$ 
	Furthermore, for the steady state $g^q$  the moments $m_k$ and $k>1$ are explicit \cite{Tosc}  

	\begin{equation}
		\label{eq:moment_q}
		m_k^q(t) = \int_{\mathbb R_+} T^k g^q(T,t)dT = \left( \dfrac{\sigma^2(1-p)}{\lambda[f](t)}\right)^{\frac{k}{2(1-p)}} \dfrac{\Gamma\left( \frac{1-2p+k}{2(1-p)} \right)}{\Gamma\left( \frac{1-2p}{2(1-p)}\right)}>0.
	\end{equation}
	Hence, the condition \eqref{eq:lambda_cond_q}, in the limit $\nu\to 0^+$, is given by 
	\begin{equation}\label{lambda-limit} 
		\begin{split}
			\lambda[f](t) &>\dfrac{\sigma^2(k-1)}{2} \dfrac{m_{k-2(p-1)}^q(t)}{m_k^q(t)}
			=: \gamma(\sigma,p,k). 
	\end{split}\end{equation}
	Note that $\sigma^2$ needs to be sufficiently small to guarantee that right--hand side $\gamma(\sigma,p,k)<1$ is less than one, such that $\lambda[f](t)$ can be also bounded from above by one. 
%%%%%%%%%
	As in Lemma \ref{lem1} we consider two cases. If $\int \F(x) (f^q_\F(x,t) - f(x,t) ) dx \geq 0$ (and $\nu$ sufficiently close to zero), we set 
	\begin{align}\label{l2}
		\lambda[f](t) = \alpha	\dfrac{  m^q_k(t)  \sqrt{H(0)} }{ \sqrt{2}\| \mathcal F\|_{\infty} }. 
	\end{align}
	Hence, the bounds on $\lambda$ can be fulfilled, e.g., for  $\alpha$ that fulfills the following  bounds
	\begin{align}\label{bound alpha}
		\frac{ \| \F \|_\infty \sigma^2 (k-1)   }{  \sqrt{H(0)} } \;  \dfrac{m_{k-2(p-1)}^q(t)}{ (m_k^q(t))^2}
		\leq  	\alpha \leq  \frac{  \sqrt{2} \| \F \|_\infty  }{   \sqrt{H(0)} m^q_k(t) }.
	\end{align}
	Note that for $\sigma$ sufficiently small, the interval for $\alpha$ is non--empty. Since $\gamma(\sigma,p,k) < \lambda[f](t) < 1,$ the bounds on $\alpha$ can be explicitly computed and only dependent on $\sigma, p,$ and $k,$ respectively. Similarly to the proof of Lemma \ref{lem1} we consider the decay
	of the relative entropy $H$ given by \eqref{eq:shannon}, but now for $f$ fulfilling equation \eqref{eq:fnu} and arbitrary initial data $f_0(x).$ The computation is analogous to the one in the previous section and yields similar to equation \eqref{eq:HF} the following equality
	\begin{equation}
		\begin{split}
			&\dfrac{d}{dt}H(f|f^q_\F)(t) = \\
			&\qquad - I_H(f|f^q_\F) + \dfrac{ k(k-1)\frac{\sigma^2}{2}m^q_{k-2(p-1)}(t) -k\lambda[f](t) m^q_k(t)  }{(m_k^q)^2(t)}\int_{\mathbb R^d}\F(x) (f^q_\F(x,t) - f(x,t))dx,
	\end{split}\end{equation}
	where the Gibb's distribution is now defined as follows
	\[
	f^q_\F = C(t) \exp\left\{ - \dfrac{ \F(x)}{\mathcal{D}[g^q](t)}\right\}, 
	\]
	being $C>0$ a normalization constant. Since $I_H \geq0$, under the assumption that $\int \F(x) (f^q_\F(x,t) - f(x,t) ) dx \geq 0$,  and the explicit form of $\lambda[f](\cdot)$ and the bounds on  $\alpha$ 
	we obtain 
	\begin{align*}
		\dfrac{d}{dt}H(f|f^q_\F)(t) & \leq   \dfrac{ k(k-1)\frac{\sigma^2}{2}m^q_{k-2(p-1)}(t) -k\lambda[f](t) \; m^q_k(t)  }{(m_k^q)^2(t)}\int_{\mathbb R^d}\F(x) (f^q_\F(x,t) - f(x,t))dx  \\
		&=  \left(\alpha \left( k \sqrt{H(0)}  / ( 2  \| \F\|_\infty )  \right) - k \alpha \sqrt{H(0)} /  ( \sqrt{2} \| \mathcal{F} \|_\infty ) \right) \int_{\mathbb R^d}\F(x) (f^q_\F(x,t) - f(x,t))dx \\
		&= - \alpha k  \left(\frac{1}{\sqrt2} -\frac12\right) \frac{ \sqrt{H(0)} }{ \| \F \|_\infty} \int_{\mathbb R^d}\F(x) (f^q_\F(x,t) - f(x,t))dx
	\end{align*}
	Proceeding as in equation \eqref{lem:H2} we obtain the decay in the case $k>1$, $\nu=0,$ and for $\alpha$ within the bounds \eqref{bound alpha} and $\lambda$ given by equation \eqref{l2} as
	\begin{align*}
		\dfrac{d}{dt}H(f|f^q_\F)(t) & \leq   - \alpha k  \left(\frac{1}{\sqrt2} -\frac12\right) H(f|f^q_\F)(t). 
	\end{align*}

	This is summarized in the following Lemma. 
	\begin{lemma}[Case $k>1$]\label{lem2}
		Consider the relative entropy $H$ defined in equation \eqref{eq:shannon} for a regular solution $f$ to equation \eqref{eq:fnu} with arbitrary initial data { $f_0 \geq0$ with mass $\int_{\mathbb R^d}f_0(x) dx=1.$ } Furthermore, let $\mathcal{D}$ be given by equation \eqref{operator D} for any $k\geq 2$. Assume  $g^q$ to be the  quasi--steady  to equation \eqref{eq:equilibrium}.  Furthermore, denote by $f^q_\F$ the quasi-equilibrium solution (depending on $\F$ and $g^q$) given by equation \eqref{eq:system_FP}. Then for any  $\alpha$ within the bounds \eqref{bound alpha} and 	for $H(0)=H(f_0, f^q_\F(\cdot,0) )$ the following statement holds:  
		
		\noindent If for $t\geq0,$  we have 
		$\int_{\mathbb R^d}\F(x) ( f^q_\F(x,t) - f(x,t))dx > 0$, then for 
		\begin{align}\label{lem:lambda2}
			\lambda[f](t) = \alpha \frac{ m_k^q(t) \sqrt{H(0)} }{ \sqrt{2} \| \F \|_\infty }
		\end{align}
		we have
		\begin{align*}
			\frac{d}{dt} H(f | f^q_\F) \leq - \alpha k \left(\frac{1}{\sqrt2} -\frac12\right)  H(f | f^q_\F).
		\end{align*}
	\end{lemma}
	
	Following the computation in Remark \ref{rem:lem1} we define a feedback in the case $\int_{\mathbb R^d} \F(x) (f^q(x,t)-f(x,t))dx <0$ as follows. We prescribe a  decay of the $k$th moment as 
	\[
	m^q_k(t )  = \dfrac{T_0}{\log(t+2)}, 
	\]
	which yields according to equation \eqref{eq:Devo} the following form for $\lambda:$ 
	\[
	\lambda[f](t) = \dfrac{\sigma^2(k-1)}{2} \dfrac{\Gamma\left( \frac{3-4p+k}{2(1-p)}\right)}{\Gamma\left( \frac{1-2p+k}{2(1-p)} \right)} + \dfrac{1}{k(t+2)\log(t+2)}>0,
	\]
	since 
	\[
	m_{k+2(1-p)}(t) = m_k(t)\dfrac{\Gamma\left( \frac{3-4p+k}{2(1-p)}\right)}{\Gamma\left( \frac{1-2p+k}{2(1-p)} \right)}
	\]
	
	Summarizing the case $k>1$ then yields the proposed feedback control 	
	\begin{equation}
		\label{eq:lambda_kappa>1}
		\lambda[f](t) = 
		\begin{cases}
			\alpha	\dfrac{  m^q_k(t)  \sqrt{H(0)} }{ \sqrt{2}\| \mathcal F\|_{\infty} }.  & \mathcal I_\F(t)\ge0 \\
			\dfrac{\sigma^2(k-1)}{2} \dfrac{\Gamma\left( \frac{3-4p+k}{2(1-p)}\right)}{\Gamma\left( \frac{1-2p+k}{2(1-p)} \right)} + \dfrac{1}{k(t+2)\log(t+2)}& \mathcal I_\F(t)<0, 
		\end{cases}
	\end{equation}
	
	Note that for all $k\geq1$  the estimate on $H$ holds only if the condition $\int_{\mathbb R^d}\F(x)(f^q(x,t)-f(x,t))dx\ge0 $ holds true for all $t\geq0.$ In the numerical results we however also observe a faster decay if the condition is only true on (small) time intervals. 
	
	\begin{remark}
		The $\lambda[f](t)$ in \eqref{eq:lambda_kappa>1} obtained in the case $\kappa>1$ for small $\nu\ll1$ is consistent with the case \eqref{eq:lambda_kappa1} in the limit  { $k \to 1^+$. }
	\end{remark}

	\section{Numerical results}\label{sec:numerics}
	In this section, we test the introduced entropic kinetic simulated annealing methods by comparing their performance with the existing version of such optimization method. In more detail,  we  define a direct simulation Monte Carlo approach for simulating the dynamics of a particle systems evolving through the kinetic non-Maxwellian model \eqref{eq:kinetic}, the {Direct Simulation Monte Carlo} (DSMC) method will build upon the classical acceptance-rejection technique for Boltzmann-type dynamics, see  \cite{pareschi_russo,PT} for an introduction. 
	
	{The algorithmic parameters are as follows: $p \in (0,\frac12)$ given by the integrability condition on the Gama density. The parameter $\alpha$ is non--negative and the upper bound is given by Lemma 3.1 and Lemma 3.3, respectively. The parameter $k$ controls the order of the moment of $g.$ Different choices for $k$ will be explored below. Further, 	} we discretize the time interval $[0,T_\textrm{max}]$ using 
	$n_\textrm{tot}$ time intervals $[t^n,t^{n+1}]$ of size $\Delta t >0$ such that $t^n = n\Delta t$, $n = 0,\dots,n_{tot}$ and we denote by $f^n(x), g^n(T)$ the approximations of $f(x,t^n)$ and $g(T,t^n)$ solution to \eqref{eq:kinetic} at time $t^n$, respectively.  Hence, for any $\nu>0$ the time discrete form of the entropic kinetic simulated annealing system \eqref{eq:kinetic} under the quasi-invariant scaling \eqref{eq:scaling} reads
	\begin{equation}\label{eq:DSMC}
		\begin{split}
			f^{n+1}(x) = \left(1-\dfrac{\Delta t}{\epsilon}\right)f^n(x) + \dfrac{\Delta t}{\epsilon}\mathcal H_\F^+(f^n,g^n) \\
			g^{n+1}(T) = \left(1-\dfrac{\Delta t}{\nu\epsilon}\right) g^n(T) + \dfrac{\Delta t}{\nu\epsilon }\mathcal J^+(g^n,f^n), 
		\end{split} 
	\end{equation}
	where $\mathcal H_\F^+(f^n,g^n)(x) = \mathcal J^+(f^n,g^n)(x) + f^n(x)$ and the gain operator of the temperature dynamics is defined as follows
	\[
	\mathcal J^+(g^n,f^n)(T) = \mathbb E_\eta \left[ \int_{\mathbb R_+} \dfrac{1}{|{}^\prime J_f|} g^n({}^\prime T)dT \right],
	\]
	{where ${}^\prime J_f$ has been defined in \eqref{eq:BG_strong}.} In the following,  we consider the case $\Delta t = \epsilon$. We report in Algorithm {1}  the DSMC method for the non-Maxwellian model \eqref{eq:kinetic}-\eqref{eq:DSMC}, see \cite{pareschi_russo} for an introduction. We highlight that the DSMC method for the second equation of \eqref{eq:DSMC} is Maxwellian whereas for the first equation, an acceptance-rejection routine should be incorporated. 
	
	\begin{algorithm}\label{alg_k1}
		\caption{Entropic kinetic simulated annealing (EntKSA) in the case $k = 1$. }\label{alg:entropic}
		\begin{algorithmic}
			\State Input parameters: $T_\textrm{max}$, $n_\textrm{ntot}$, $\epsilon>0$, $p\in (0,\frac{1}{2})$, $\theta \in (0,1) $
			\State Initialize $N>0$ trial points $X^{(0)}$, $T^{(0)}$
			\While{$t<n_{\textrm{tot}}$ }
			\State Generate $\xi \sim \mathcal N(0,1)$, $\eta \sim \mathcal U([-a,a])$ {from \eqref{eq:eta_bound_mod}}
			\State Compute $\lambda^{(t)}$ defined in \eqref{eq:lambda_kappa1} and $m_1^{(t)} = \frac{1}{N} \sum_{j=1}^n T_j^{(n)}$. 
			\State Compute $\tilde X = X^{(t)} + \sqrt{2\epsilon m_1^{(t)}} \xi$ 
			\State Evaluate $B^g_\F(X^{(t)} \to \tilde X)$ defined in \eqref{eq:kernel} with Gibbs measure \eqref{eq:BG}. 
			\State Compute $T^{(t+1)} = T^{(t)} - \lambda^{(t)}T^{(t)} + (T^{(t)})^{p} \chi(T^{(t)}\ge (1-p)\theta) \eta$
			\If{$\eta \le B_\F^g$}
			\State $X^{(t+1)} \gets \tilde X$
			\Else 
			\State $X^{(t+1)} \gets X^{(t)}$
			\EndIf
			\State $t \leftarrow t+1$
			\EndWhile
		\end{algorithmic}
	\end{algorithm}
	
	In the case $k>1$ we proposed to rely on the specific form of the quasi-equilibrium distribution for the temperature variable $g^q(T)$ in \eqref{eq:equilibrium}, that corresponds to a generalized gamma density 
	\begin{equation}
		\label{eq:gen_gamma}
		g^q(T,t) = \dfrac{\ell/a^d}{\Gamma(d/\ell)}T^{d-1}\exp\left\{ -(T/a)^d\right\},
	\end{equation}
	where
	\[
	d= 1-2p>0, \qquad \ell = 2(1-p)>0, \qquad a =\left( \dfrac{\lambda[f](t)}{\sigma^2(1-p)}\right)^{2(1-p)}>0. 
	\]
	The quasi-equilibrium, since $\lambda[f] = \lambda[f](t)$, characterizing the temperature dynamics is such that  \eqref{eq:moment_q} does hold. 
	
	\begin{algorithm}\label{alg_k>1}
		\caption{Entropic kinetic simulated annealing (EntKSA) in the case $k>1$. }
		\begin{algorithmic}
			\State Input parameters: $T_\textrm{max}$, $n_\textrm{ntot}$, $\epsilon>0$, $p\in (0,\frac{1}{2})$, $\theta \in (0,1) $
			\State Initialize $N>0$ trial points $X^{(0)}$, $T^{(0)}$
			\While{$t<n_{\textrm{tot}}$ }
			\State Generate $\xi \sim \mathcal N(0,1)$, $\eta \sim \mathcal U([-a,a])$  {from \eqref{eq:eta_bound_mod}}
			\State Compute $\lambda^{(t)}$ defined in \eqref{eq:lambda_kappa>1} and $m_k^{(t)}$ from \eqref{eq:Devo}. 
			\State Compute $\tilde X = X^{(t)} + \sqrt{2\epsilon m_k^{(t)}} \xi$ 
			\State Evaluate $B^g_\F(X^{(t)} \to \tilde X)$ defined in \eqref{eq:kernel} with Gibbs measure \eqref{eq:BG}. 
			\If{$\eta \le B_\F^g$}
			\State $X^{(t+1)} \gets \tilde X$
			\Else 
			\State $X^{(t+1)} \gets X^{(t)}$
			\EndIf
			\State $t \leftarrow t+1$
			\EndWhile
		\end{algorithmic}
	\end{algorithm}
	
	We compare the kinetic simulated annealing (KSA) and the entropic kinetic simulated annealing (EntKSA) approaches and validate their corresponding meanfield limit. This latter result may open new directions in both theoretical and applied fields as it highlights that gradient-based outcomes can be achieved through gradient-free methods which dissipate the entropy, in a suitable regime of parameters. 
	
	To this end, we will stick to a 1D case and we concentrate on a case where $\mathcal I_\F(0)>0$. In particular, we define the following nonconvex cost function
	\begin{equation}\label{eq:cost_function}
		\F(x) = \begin{cases}
			\cosh(\frac{x}{4})+3 & x \in \mathbb R\setminus [0,2] \\
			\cosh(\frac{x}{4}) - \cosh(x) + 3 & x\in [0,2]. 
		\end{cases}
	\end{equation}
	
	In addition to EntKSA we considered the solution obtained by solving the KSA proposed in \cite{P24} that is connected to the classical meanfield dynamics of the simulated annealing algorithm. {In all the subsequent tests we reconstructed the particle density over the integral $[-20,20]$ by means of $N_x = 501$ gridpoints. }
	
	\subsection{The case $k=1$}
	{In} the following, we will always assume $k = 1$ and $\nu = 1$ in \eqref{eq:Devo} defining the kinetic dynamics in \eqref{eq:kinetic}. In all the tests we will fix $N = 10^6$ and $\Delta t = \epsilon$, $p = 1/4$ and $\theta = 0.5$. As initial distribution  we consider
	\[
	\begin{split}
		f_0(x) = 
		\begin{cases}
			1 & x \in [1,2] \\
			0 & \textrm{elsewhere},
		\end{cases}
		\qquad
		g_0(T) = \delta(T-T_0),
	\end{split}
	\]
	where $T_0 = 2$. 
	
\begin{figure}
		\centering
		\includegraphics[scale = 0.35]{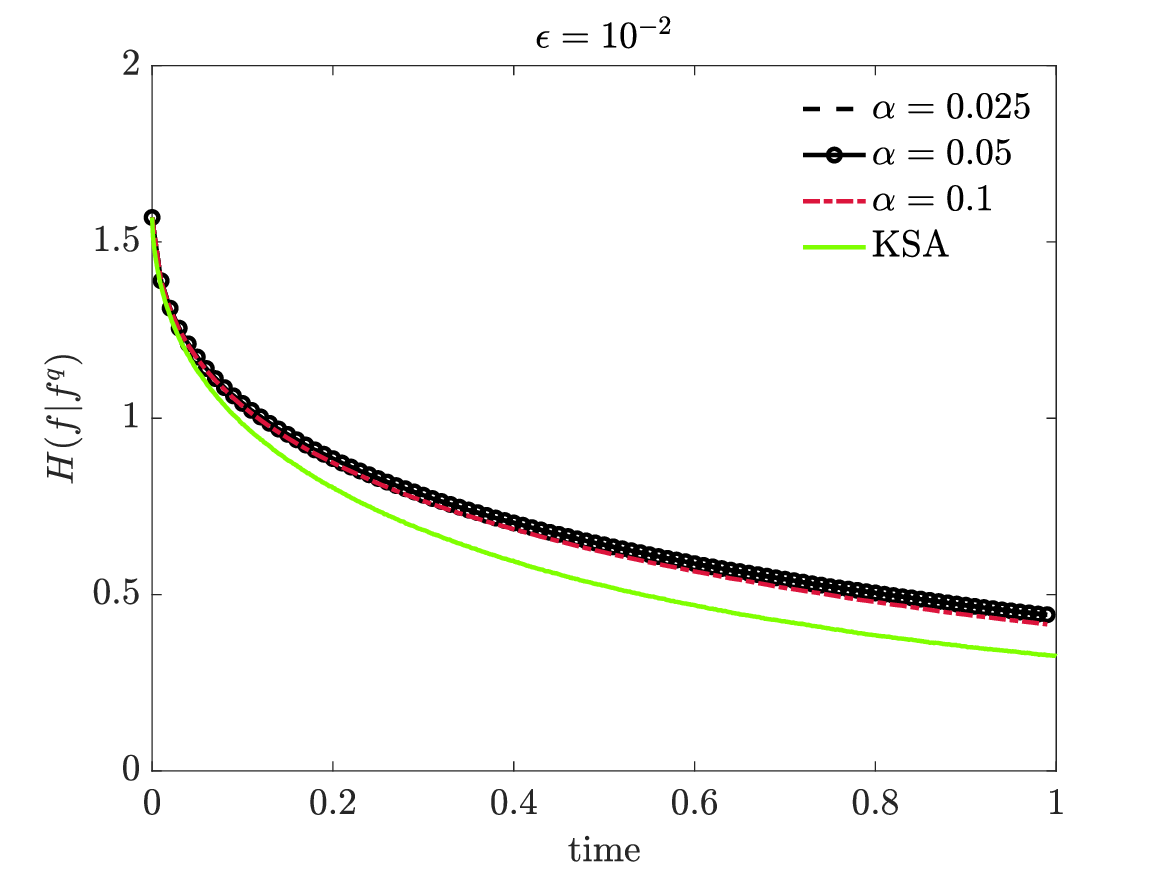}
		\includegraphics[scale = 0.35]{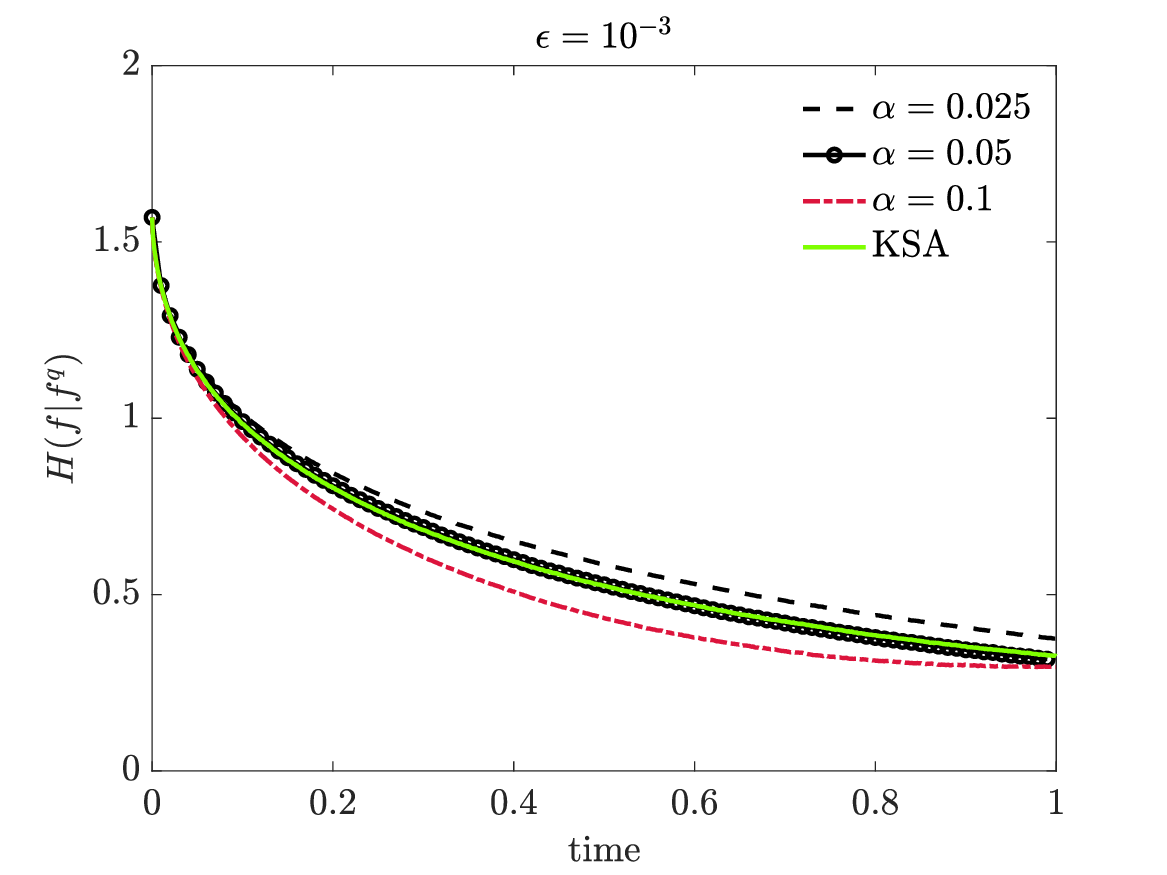} \\
		\includegraphics[scale = 0.35]{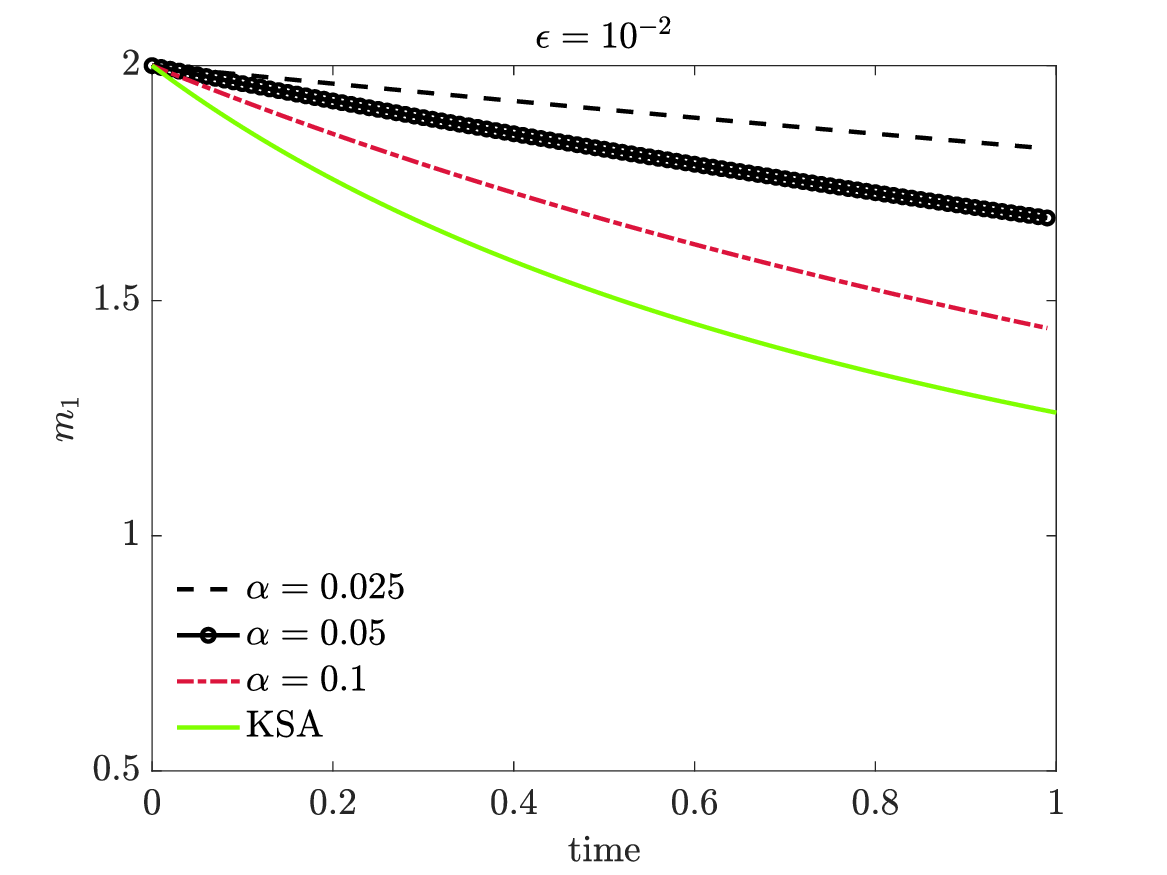} 
		\includegraphics[scale = 0.35]{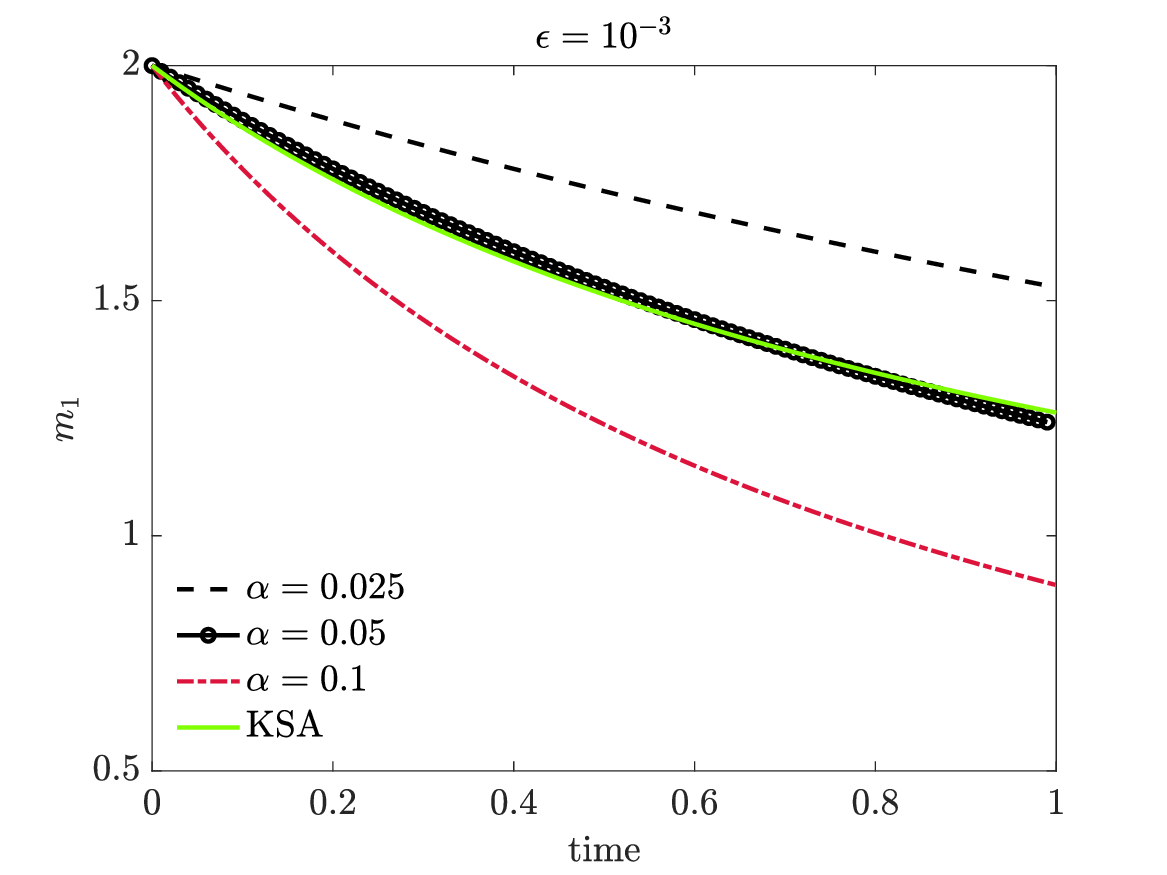} 
		\caption{Top row: evolution of the Shannon entropy $H(f|f^q)(t)$ for several values of $\alpha$ and for $\epsilon = 10^{-2}$ (left) and $\epsilon = 10^{-3}$ (right), we have reported in green approximated evolution of the KSA algorithm with $\epsilon = 10^{-4}$. Bottom row: evolution of the mean temperature of the system of particles $m_1(t) = \int_{\mathbb R_+}Tg(T,t)dT$ for $\epsilon = 10^{-2}$ (left) and $\epsilon = 10^{-3}$ (right). }
		\label{fig:ent_1}
	\end{figure}

	In Figure \ref{fig:ent_1} we present the evolution of the relative Shannon entropy defined by equation \eqref{eq:shannon} for two values of $\epsilon \in \{  10^{-3}, 10^{-2} \}$ and $\alpha \in \{  0.025,0.05,0.1 \}$ obtained by Algorithm {1}. As a reference entropy, we plot  the approximated evolution of the Shannon entropy of the KSA algorithm  for $\epsilon = 10^{-4}$. We  observe  for increasing values of $\alpha>0$ faster decay of the entropy (in the regime $\epsilon \ll 1$).  This result is further highlighted by the evolution of the mean temperature of the system of particles 
	\[
	m_1(t) = \int_{\mathbb R_+}Tg(T,t)dT
	\]
	The  evolution of $m_1$  is provided in the second row of Figure \ref{fig:ent_1}. We remark that higher values of $\alpha>0$ lead to faster decay of the mean temperature. 
	
\begin{figure}
		\centering
		\includegraphics[scale = 0.35]{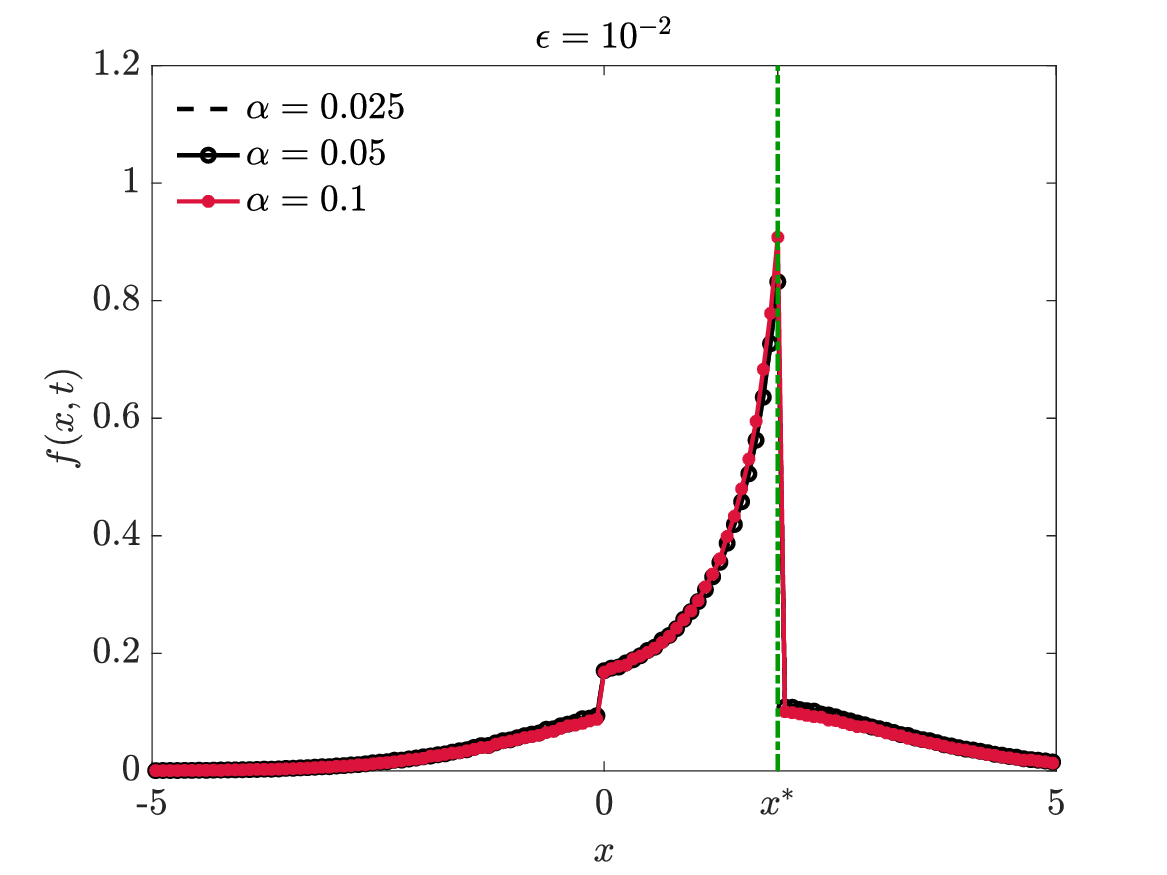}
		\includegraphics[scale = 0.35]{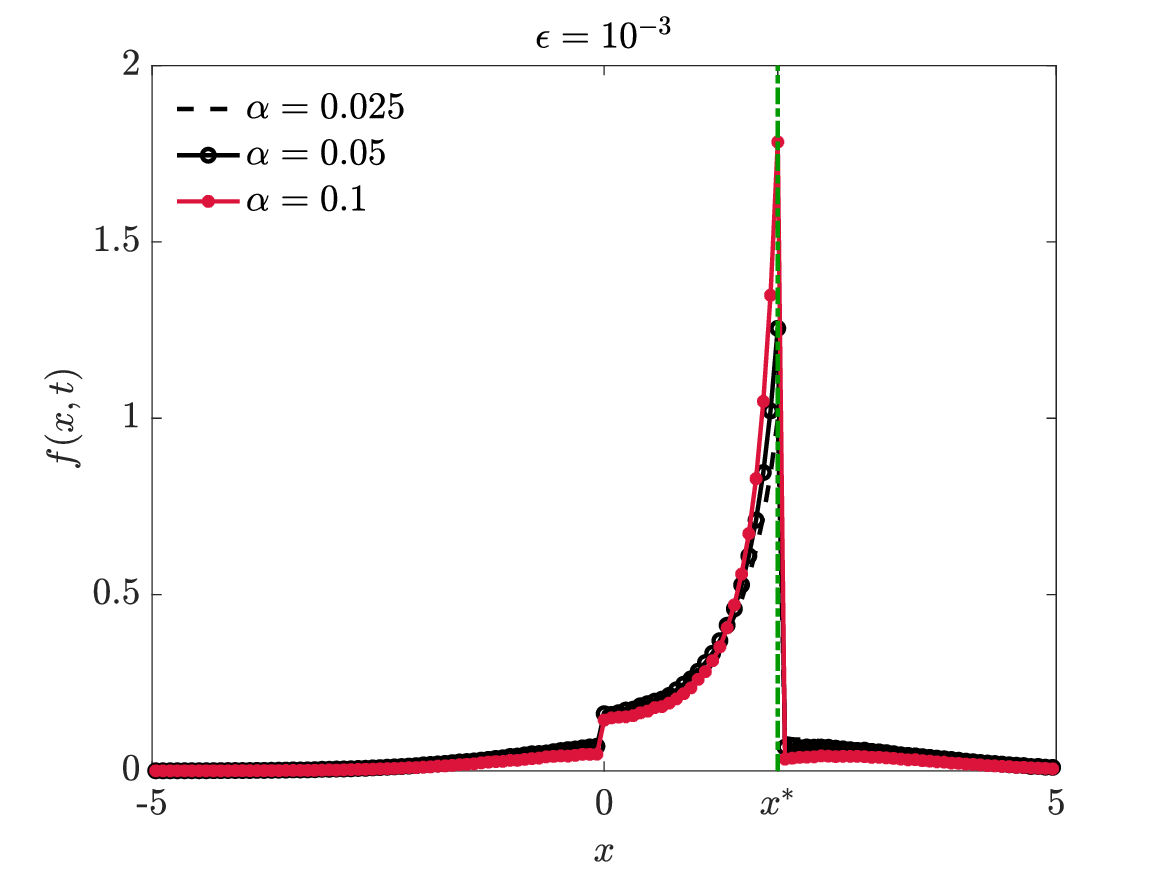}
		\caption{The distribution $f(x,t)$ at time $t = 1$ for several values of $\alpha = 0.025,0.05,0.1$ and for $\epsilon = 10^{-2}$ (left) or $\epsilon= 10^{-3}$ (right). In dotted green, we highlighted $x^*$ corresponding to the  exact minimum of the function $\F(x)$.  }
		\label{fig:dis_k1}
	\end{figure}

\begin{figure}
		\centering
		\includegraphics[scale = 0.225]{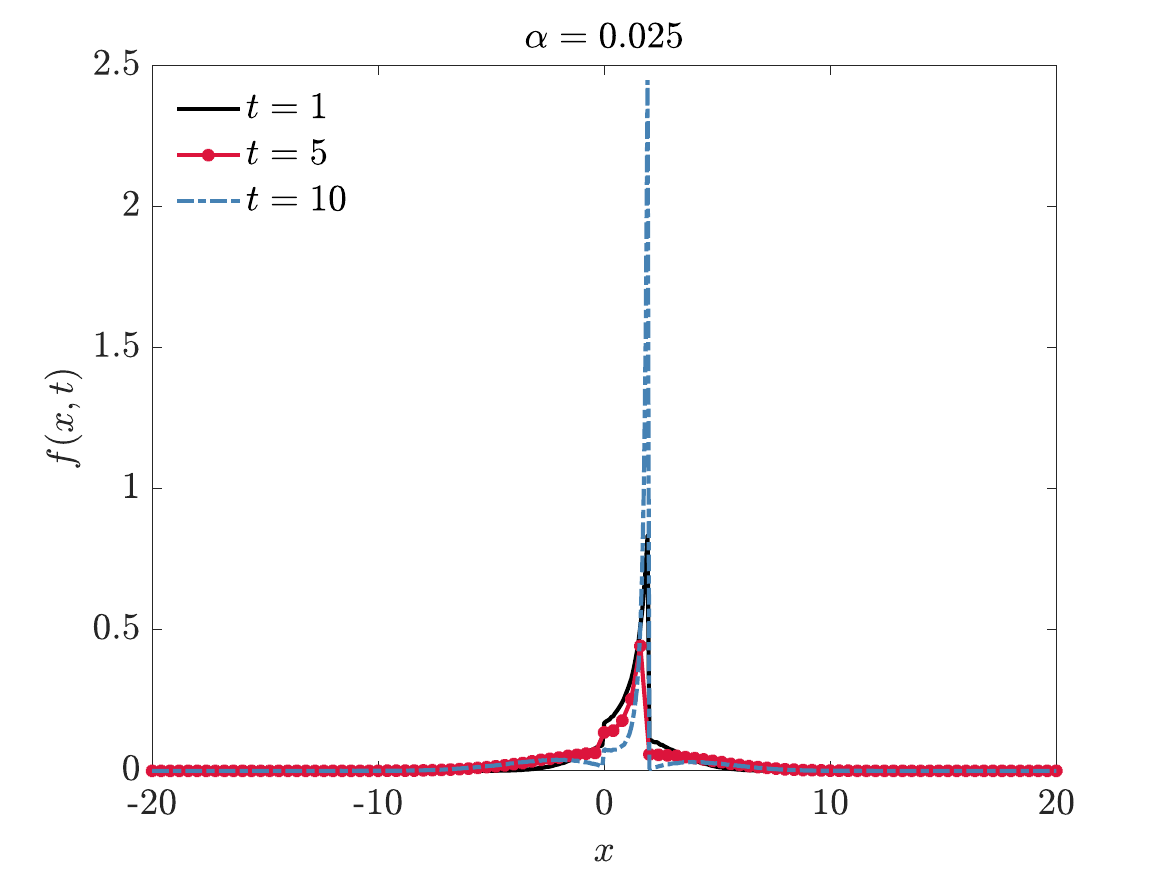}
		\includegraphics[scale = 0.225]{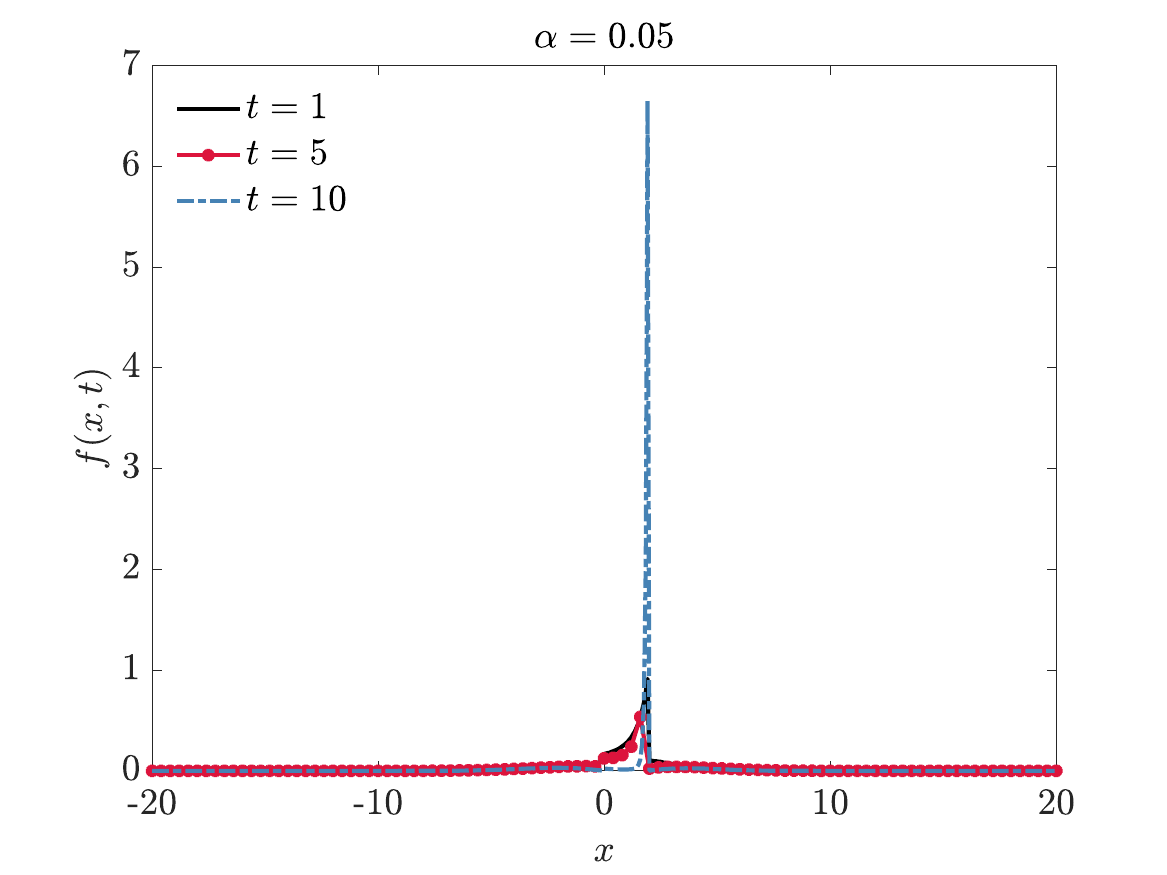}
		\includegraphics[scale = 0.225]{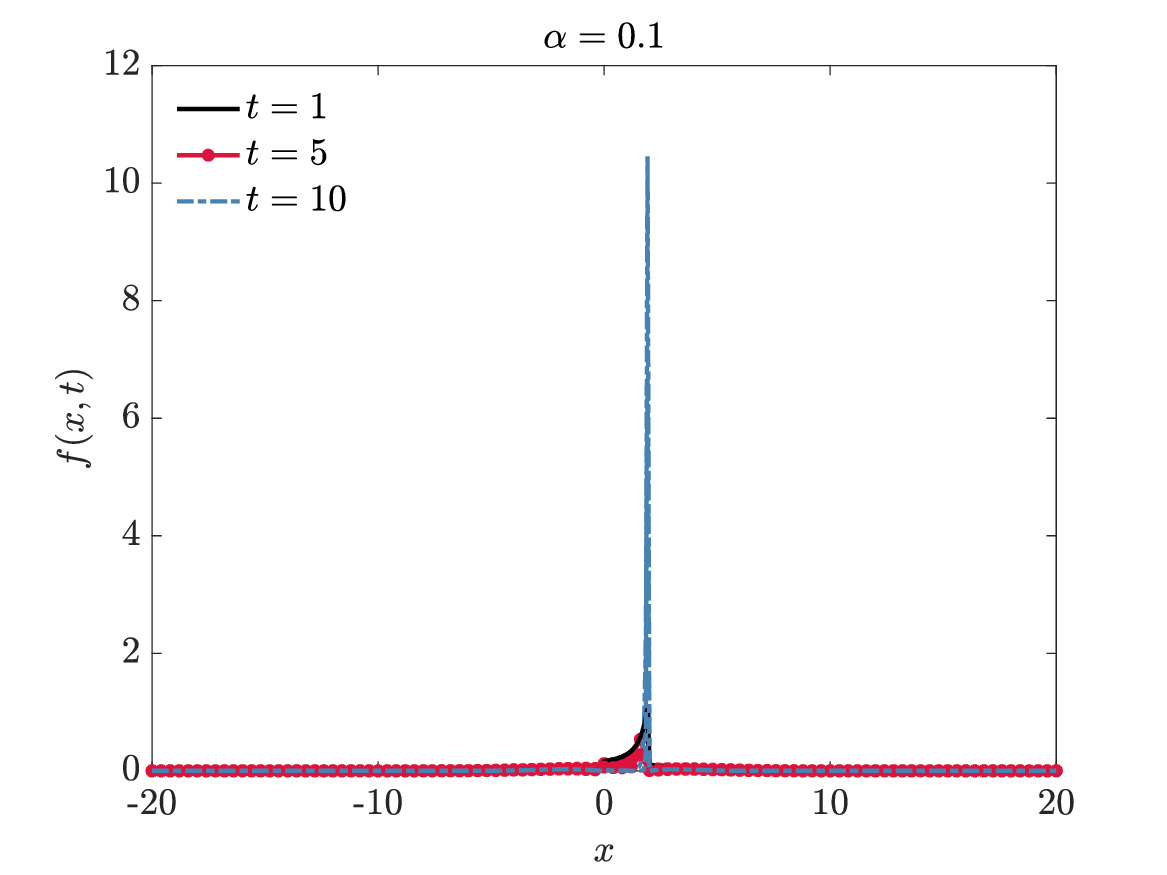} \\
		\includegraphics[scale = 0.225]{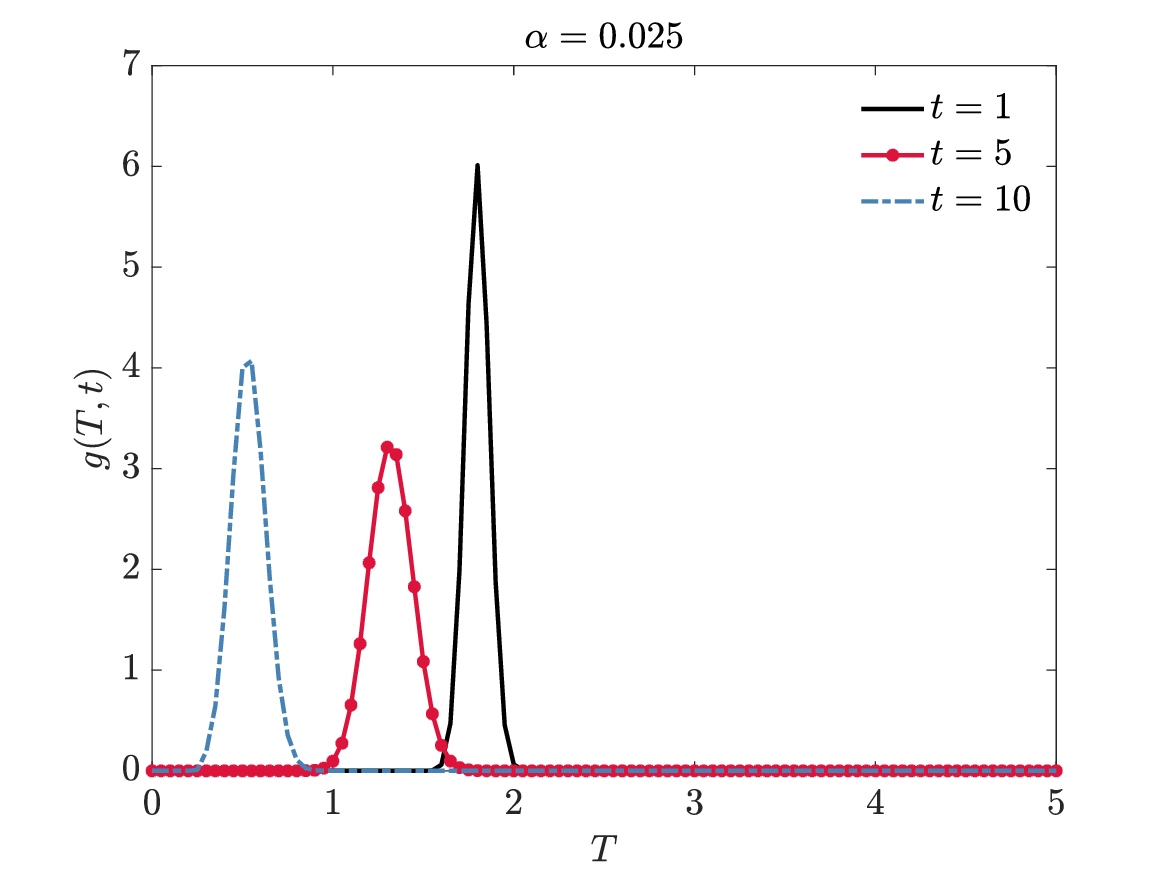}
		\includegraphics[scale = 0.225]{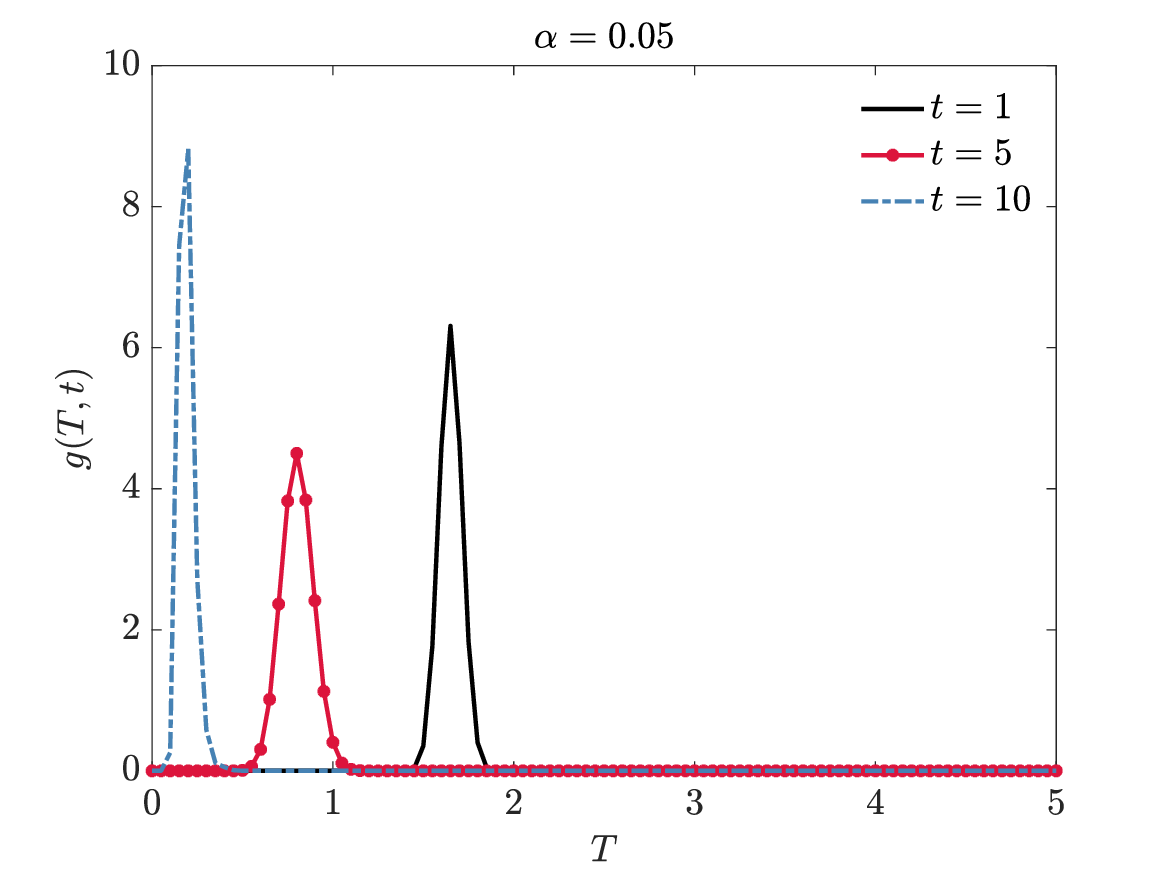}
		\includegraphics[scale = 0.225]{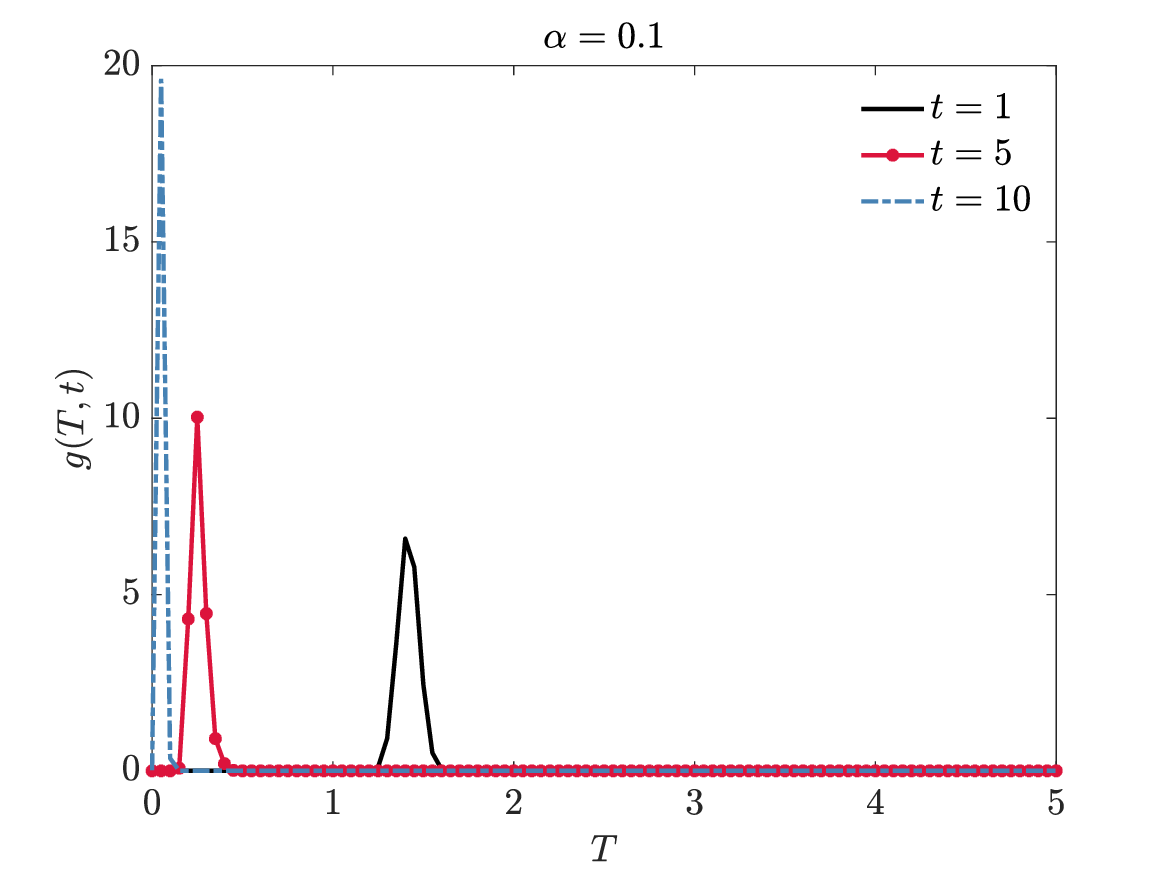}
		\caption{Evolution of the distributions $f(x,t)$, $g(T,t)$ at time $t = 1,5,10$ for several values of $\alpha = 0.025$ (left) $\alpha = 0.05$ (center), and $\alpha = 0.1$ (right)  for  a fixed $\epsilon= 10^{-3}$. In the top row, we report the evolution of $f(x,t)$, and in the bottom row the evolution of $g(T,t)$. We considered $N = 10^6$ particles both in space and temperature, $p= 1/4$, $\theta = 0.05$ and $\sigma^2 = 0.1$. }
		\label{fig:dis_fg}
	\end{figure}

\begin{figure}
		\centering
		\includegraphics[scale = 0.35]{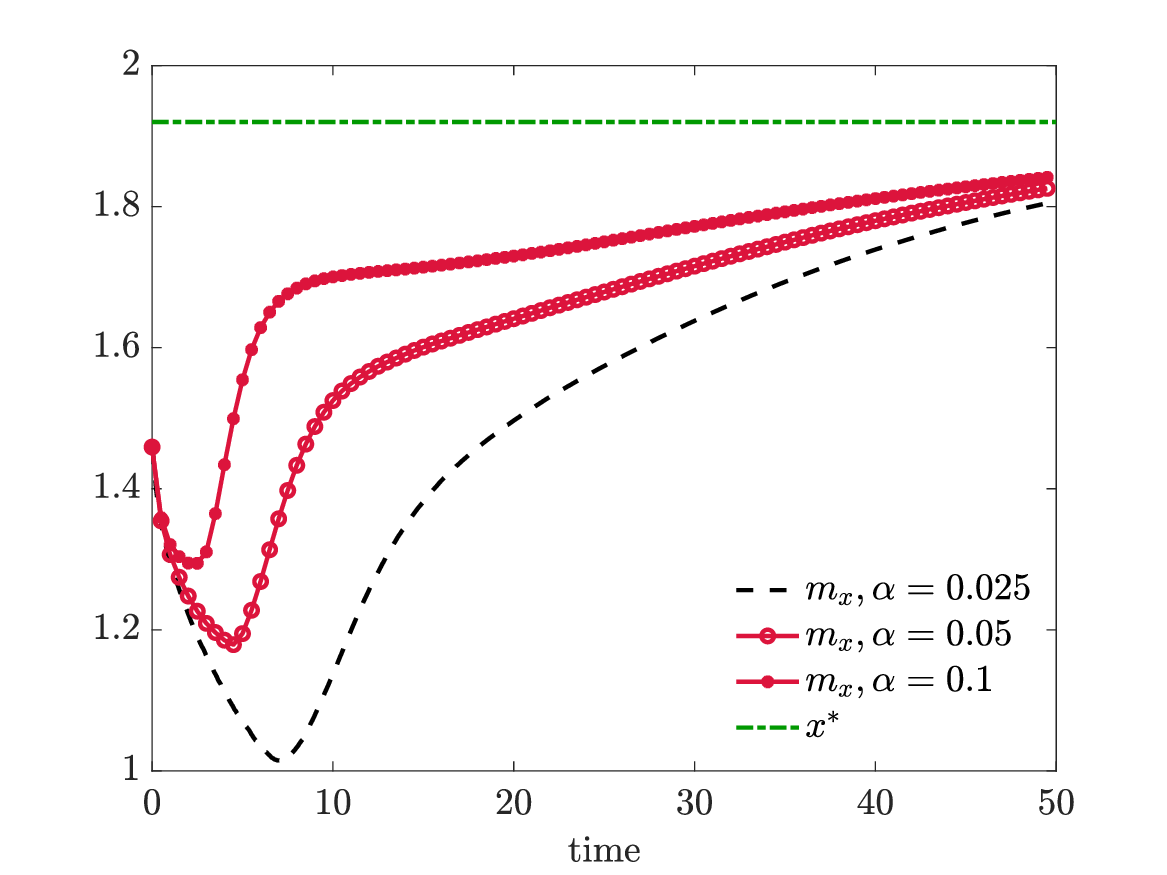}
		\includegraphics[scale = 0.35]{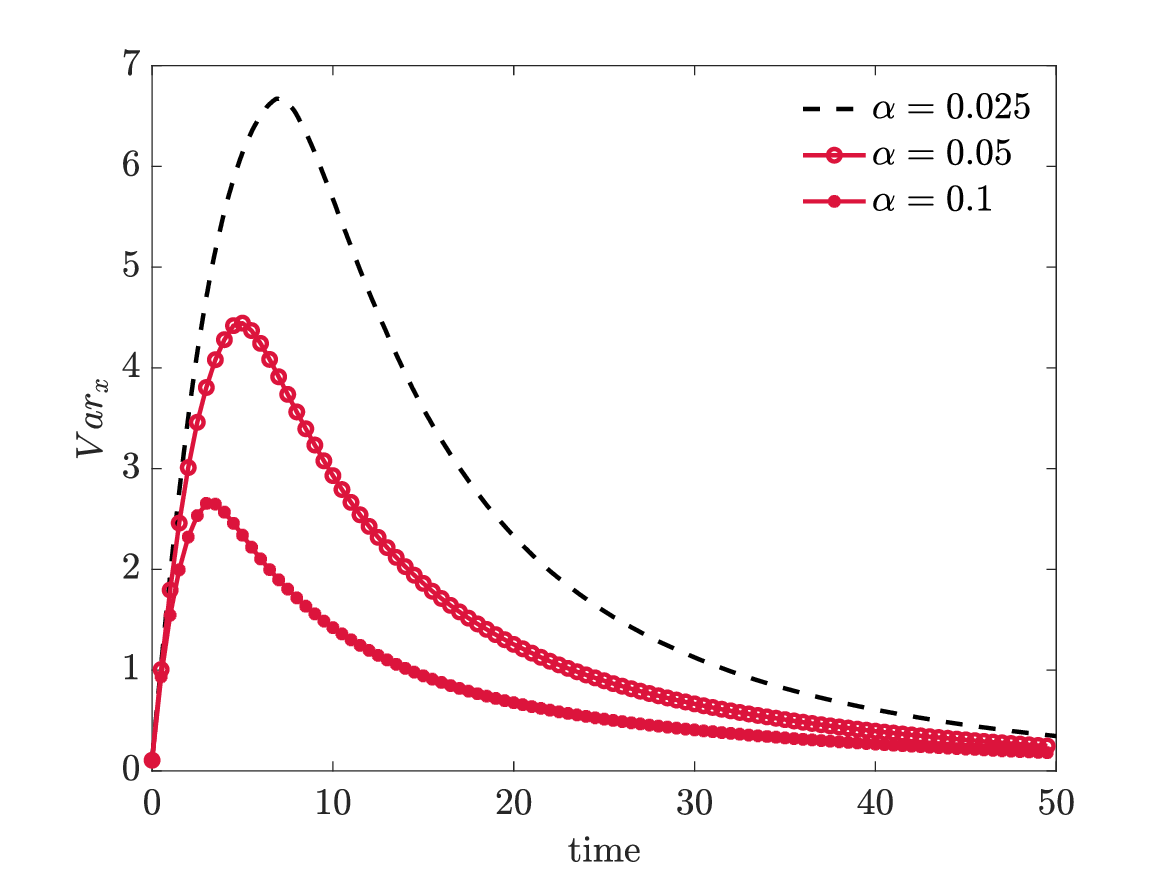}
		\caption{Evolution of the mean position $m_x = \int_{\mathbb R} x f(x,t)dx$ and of its variance $Var_x = \int_{\mathbb R} (x-m_x)^2 f(x,t)dx$ over the time interval $[0,50]$ and several values of the parameter $\alpha = 0.025,0.05,0.1$. }
		\label{fig:mx_varx}
	\end{figure}

	At the level of the kinetic distribution $f(x,t)$, we report in Figure \ref{fig:dis_k1} its approximation at time $t=1$ for several values of $\alpha>0$ and two values of $\epsilon = 10^{-3},10^{-2}$. {As expected a faster concentration toward $x^*$ is observed for increasing values of $\alpha>0$.} To highlight the convergence towards the correct minimum we reported with a dotted line the location of $x^* \in \mathbb R$. These trends can be further observed for longer times in Figure \ref{fig:dis_fg} where we consider the dynamics of $f(x,t)$ and $g(T,t)$ over the time interval $[0,T]$, $T = 10$, and in Figure \ref{fig:mx_varx} where we depict the evolution of the mean position and of its variance 
	\[
	m_x(t) = \int_{\mathbb R} xf(x,t)dx, \qquad \textrm{Var}_x(t) = \int_{\mathbb R} (x-m_x)^2f(x,t)dx. 
	\]

\begin{figure}
		\centering
		\includegraphics[scale = 0.35]{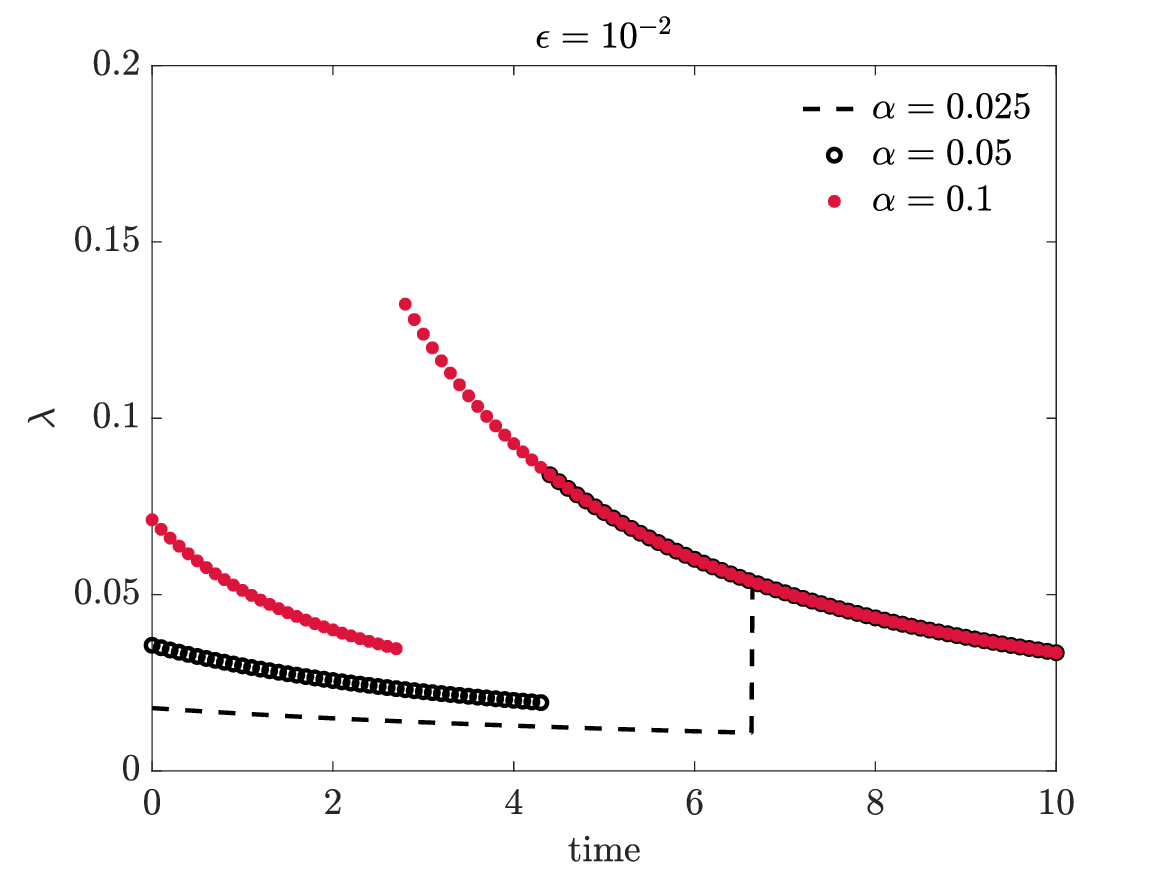}
		\includegraphics[scale = 0.35]{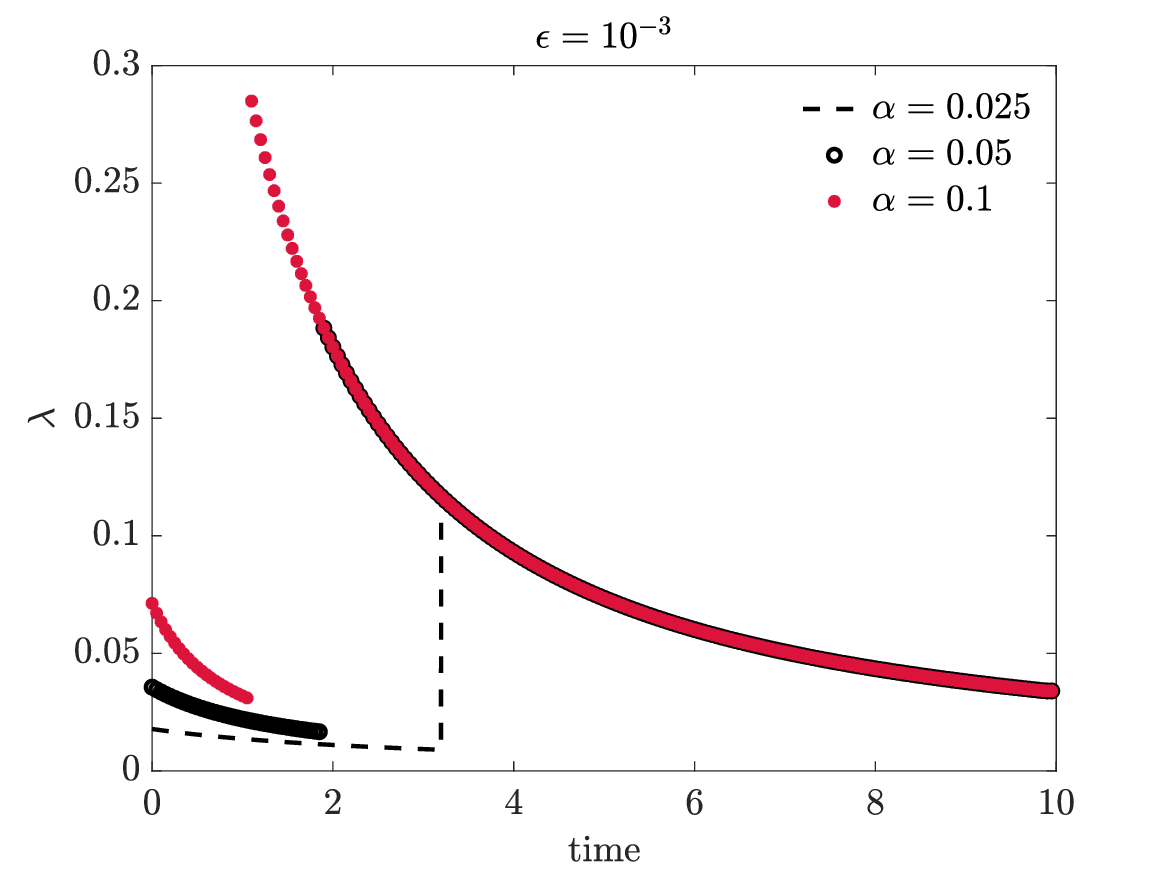} \\
		\includegraphics[scale = 0.35]{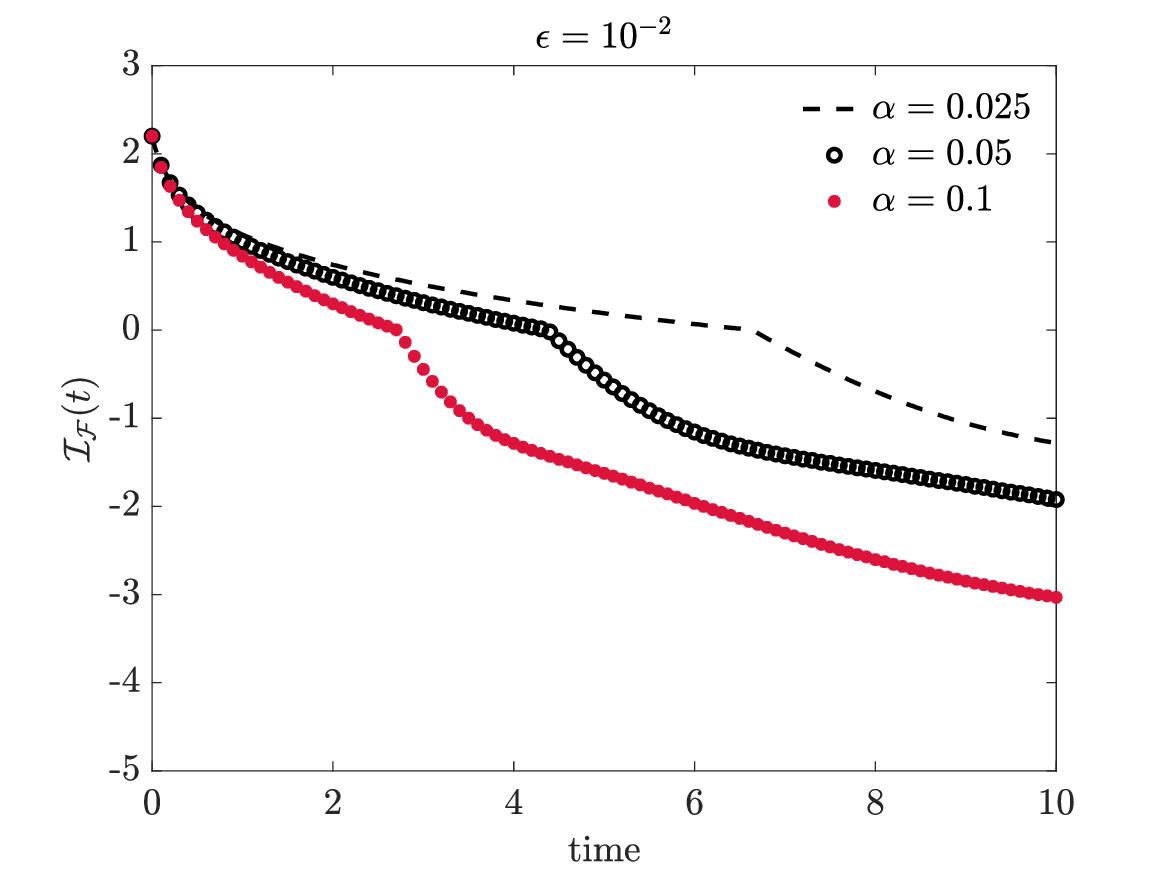}
		\includegraphics[scale = 0.35]{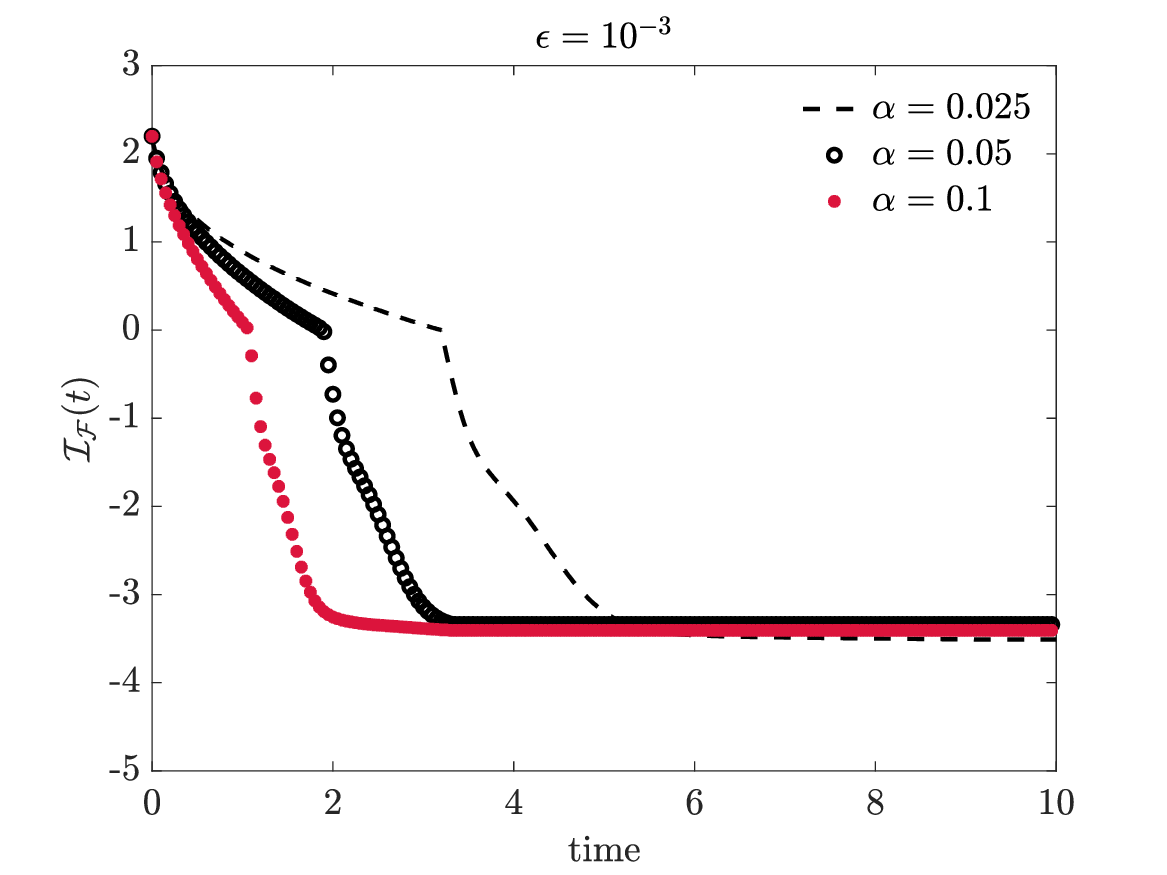}
		\caption{Top row: estimated values of $\lambda>0$ defined in \eqref{eq:lambda_kappa1} and several values of $\alpha = 0.025,0.05,0.1$. Bottom row: estimated values of $\mathcal I_\F(t) = \int_{\mathbb R} \F(x)(f^q(x,t)-f(x,t))dx$ for several values of $\alpha = 0.025,0.05,0.1$. We considered $\epsilon = 10^{-2}$ (left) and $\epsilon = 10^{-3}$ (right).  }
		\label{fig:lambda}
	\end{figure}

	In this numerical example, the obtained $\lambda>0$ has been computed by equation \eqref{eq:lambda_kappa1} and its trend is given in Figure \ref{fig:lambda}.  In particular, we observe that for large values of $\alpha$, the quantity  $\mathcal I_\F(t) = \int_{\mathbb{R}} \mathcal{F}(x)(f^q(x,t) - f(x,t))dx$ is negative for smaller times. Then, the feedback control is active and the convergence towards the local equilibrium is enhanced. The dynamics of $\mathcal I_\F(t)$ is reported in the second row of Figure \ref{fig:lambda}. 
	
	{To assess the impact of the particle number $N$ on the optimization method, we report the success rate of the EntSA method over $N_{\textrm{sim}} = 200$ simulations, computed as
\[
N_{\textrm{ave}} = \dfrac{\#\{j \in \{1,\dots,N\} : |x_j(T) - x^*| < \delta_{\textrm{threshold}}\}}{N},
\]
where $x^*$ is defined in \eqref{eq:optim} and $\delta_{\textrm{threshold}} = 0.25$, see \cite{CJLZ}.
We consider $N = 50, 100, 200$ and collect the results at final times $T = 50$ and $T = 100$ in Table~\ref{table:combined}. We observe that EntSA preserves good performance even for a relatively small number of particles. However, the analytical properties exploited by the method hold in the mean-field (PDE) regime, i.e.\ for large $N$. This behavior is consistent with the numerical results, which indicate that larger values of $\alpha > 0$ require an increasing number of particles.}

	{
	\begin{table}
\centering
\begin{tabular}{|c|c|c|c|c|}
\hline
    \multicolumn{2}{|c|}{} & \multicolumn{3}{c|}{$N_{\textrm{ave}}$} \\  \cline{3-5}
    \multicolumn{2}{|c|}{}   & $\alpha = 0.025$ & $\alpha = 0.05$ & $\alpha = 0.1$ \\
\hline
\multirow{3}{*}{$T = 50$}
& $N=50$  & 0.9486 & 0.9498 & 0.9328 \\ \cline{2-5}
& $N=100$ & 0.9493 & 0.9501 & 0.9485 \\ \cline{2-5}
& $N=200$ & 0.9548 & 0.9532 & 0.9527 \\
\hline
\multirow{3}{*}{$T = 100$}
& $N=50$  & 0.9825 & 0.9695 & 0.9376 \\ \cline{2-5}
& $N=100$ & 0.9940 & 0.9908 & 0.9656 \\ \cline{2-5}
& $N=200$ & 0.9960 & 0.9954 & 0.9858 \\
\hline
\end{tabular}
\caption{Values of $N_{\textrm{ave}}$ for different $N$, $\alpha>0$, and final times $T$.}
\label{table:combined}
\end{table}
}
	
	\subsection{The case $k>1$ }
	In this section, we focus on the case where the temperature dynamics is considered to evolve faster compared with the one of $f(x,t)$, i.e. $\nu \to 0^+$. This allows to obtain conditions on $\lambda[f](t)>0$ by the quasi-equilibrium state of the temperature dynamics in the form of a generalized gamma density \eqref{eq:gen_gamma}. We recall the results in Lemma \ref{lem2} for the explicit form of the parameter $\lambda[f](t)$ for a general $k>1$, $t\ge0$ as in  \eqref{eq:lambda_kappa>1}. In the following, we will consider the nonconvex cost function defined in \eqref{eq:cost_function} and the parameters $\sigma^2 = 1/10$, $p = 1/4$, $\theta = 1/2$. Hence, the resulting kinetic equation is approximated through {Algorithm 2}. 
	
	We considered as initial distribution 
	\begin{equation}
		\label{eq:f0T0_k>1}
		f_0(x) = \begin{cases}
			1 & x \in [1,2], \\
			0 & \textrm{elsewhere}
		\end{cases}
	\end{equation}
	whereas we computed the initial temperature such that 
	$m_k(0) = \frac{2}{\log(2)}$. 
	
\begin{figure}
		\centering
		\includegraphics[scale = 0.225]{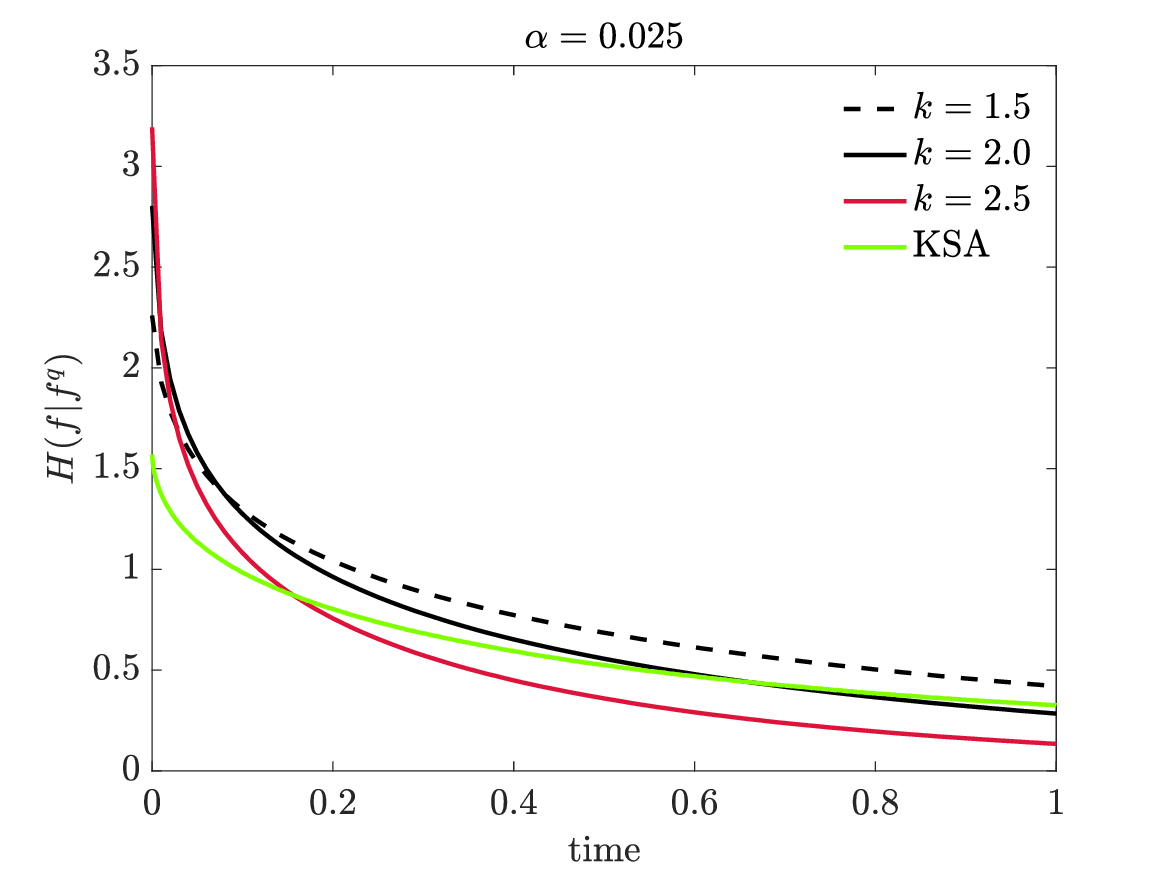}
		\includegraphics[scale = 0.225]{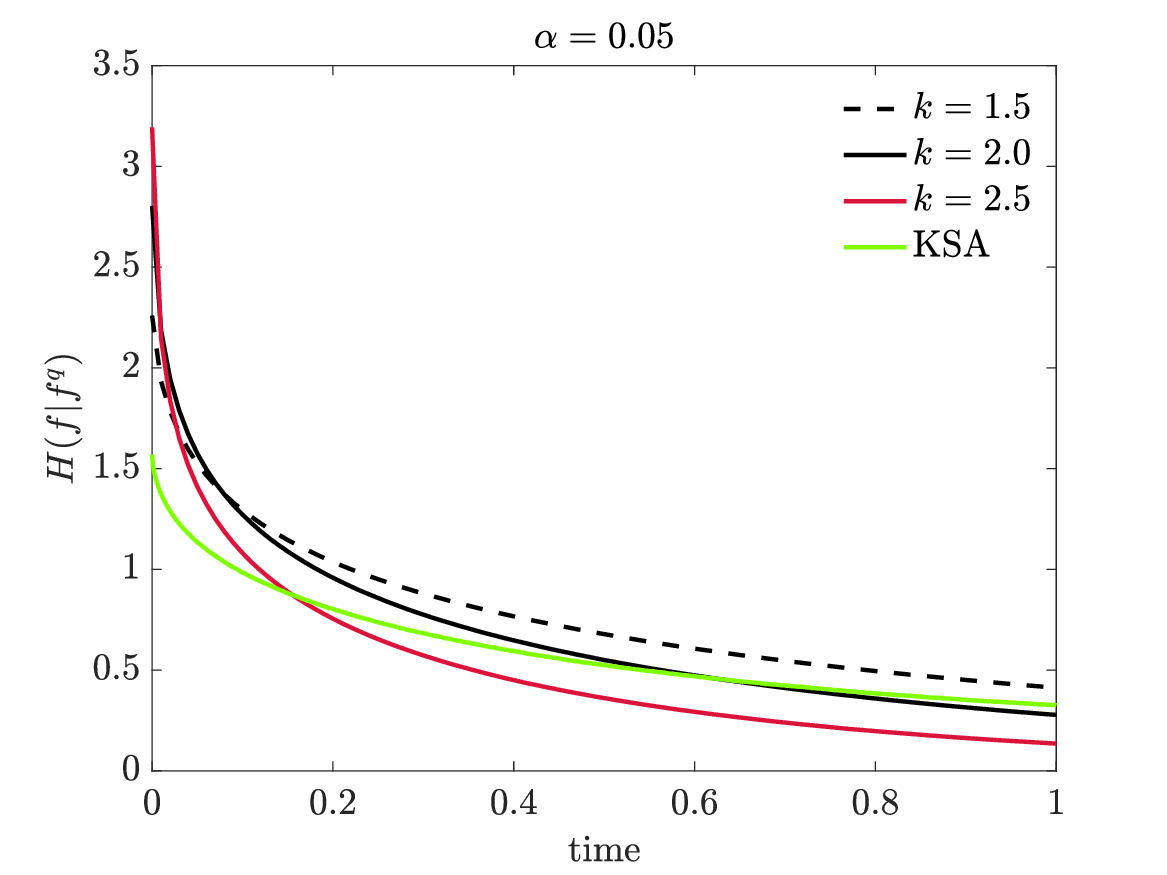}
		\includegraphics[scale = 0.225]{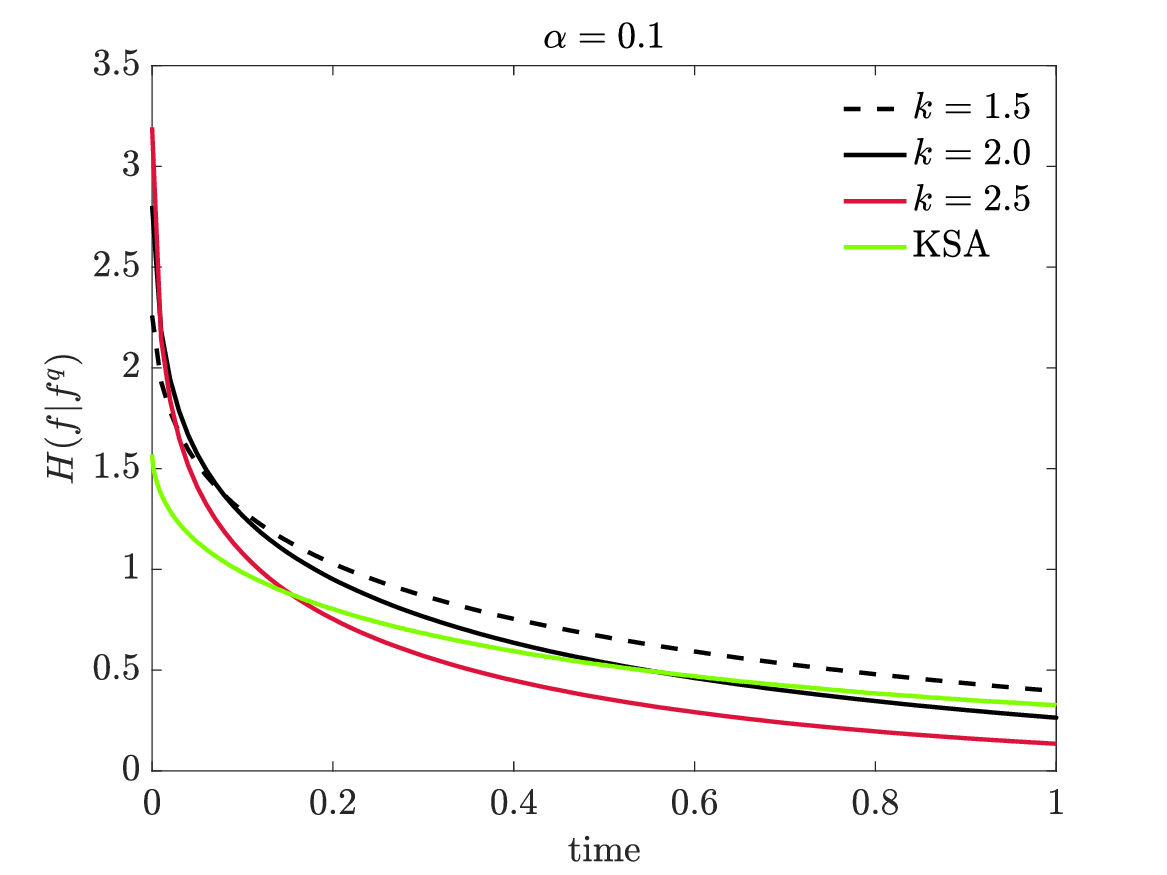}\\
		\includegraphics[scale = 0.225]{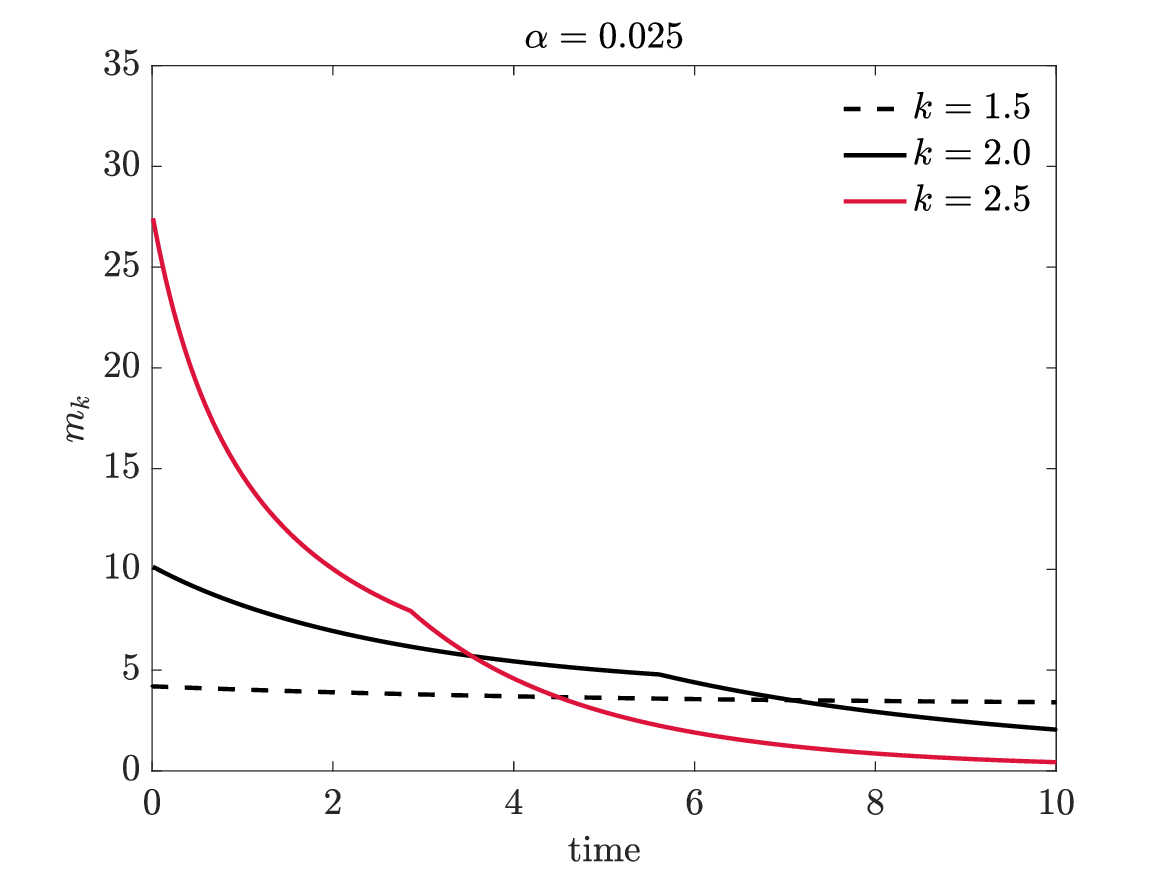} 
		\includegraphics[scale = 0.225]{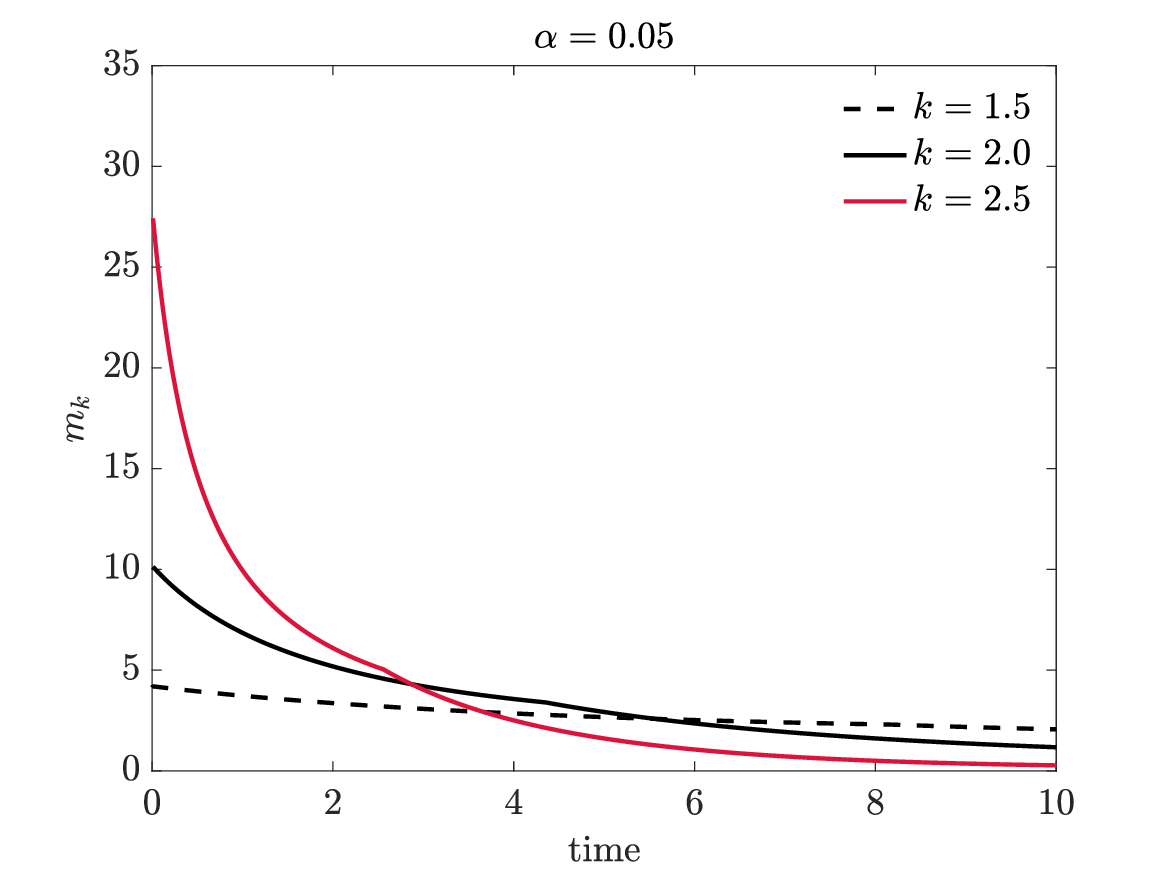} 
		\includegraphics[scale = 0.225]{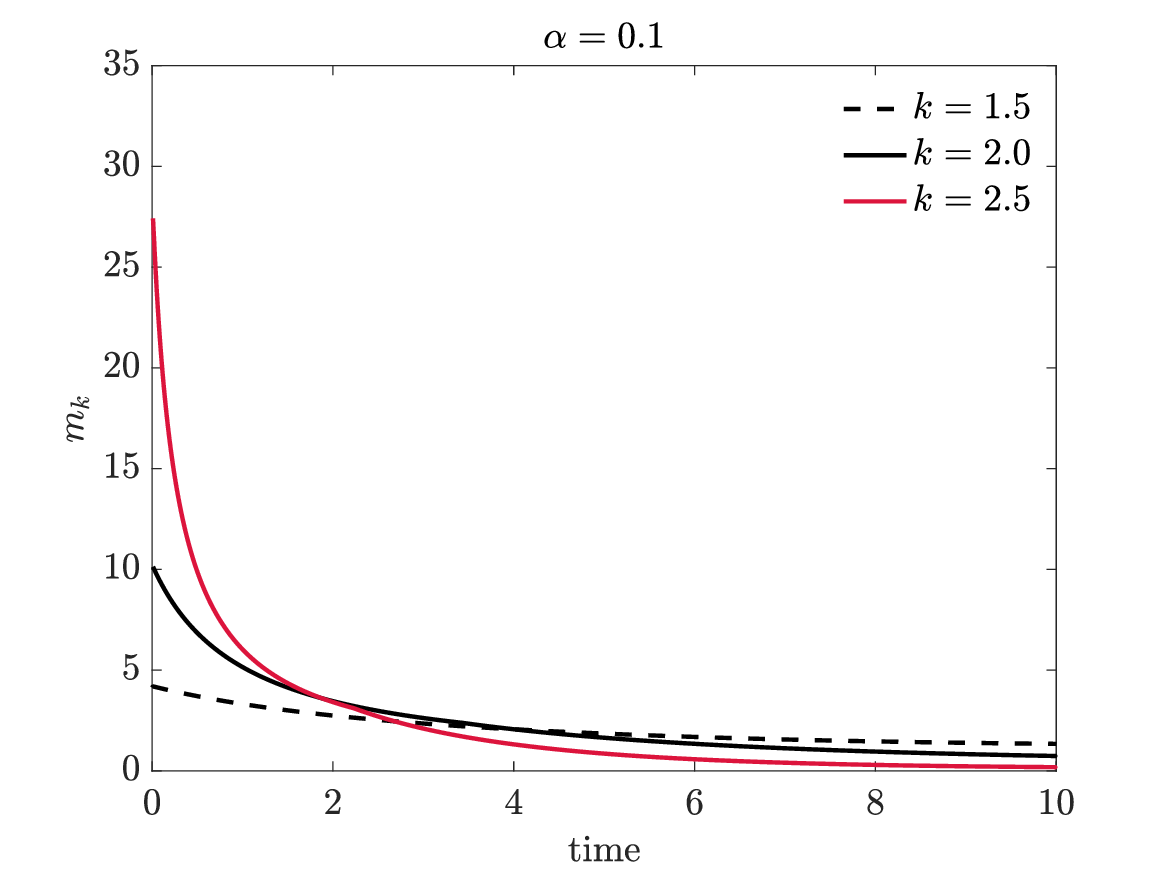} 
		\caption{Top row: evolution of the Shannon entropy $H(f|f^q)(t)$ for several values of $\alpha$ and for $k  \in\{ 1.5,2.0,2.5\}$ (left) . Bottom row: evolution of the mean temperature of the system of particles $m_k(t) = \int_{\mathbb R_+}Tg^q(T,t)dT$. }
		\label{fig:ent_k>1}
	\end{figure}

	In Figure \ref{fig:ent_k>1} we show the evolution of the Shannon entropy for the resulting EntKSA with quasi-equilibrium temperature dynamics in the cases {$k   \in \{ 1.5,2.0,2.5 \}$} and for several values of $\alpha \in \{  0.025,0.05,0.1  \}$ and under  the quasi-invariant scaling, $\epsilon = 10^{-2}$. As before, as reference entropy, we also plotted the evolution of the Shannon entropy of the KSA algorithm obtained with $\epsilon = 10^{-4}$. We  observe  even in the quasi-equilibrium case the dissipation of the entropy. Furthermore, the evolution of the moment $m_k(t)$ over the time interval $[0,10]$ is shown and it  decaying over time. Its decay remains coupled to the dynamics of $f(x,t)$ through $\lambda[f](t)$. 
	
\begin{figure}
		\centering
		\includegraphics[scale = 0.35]{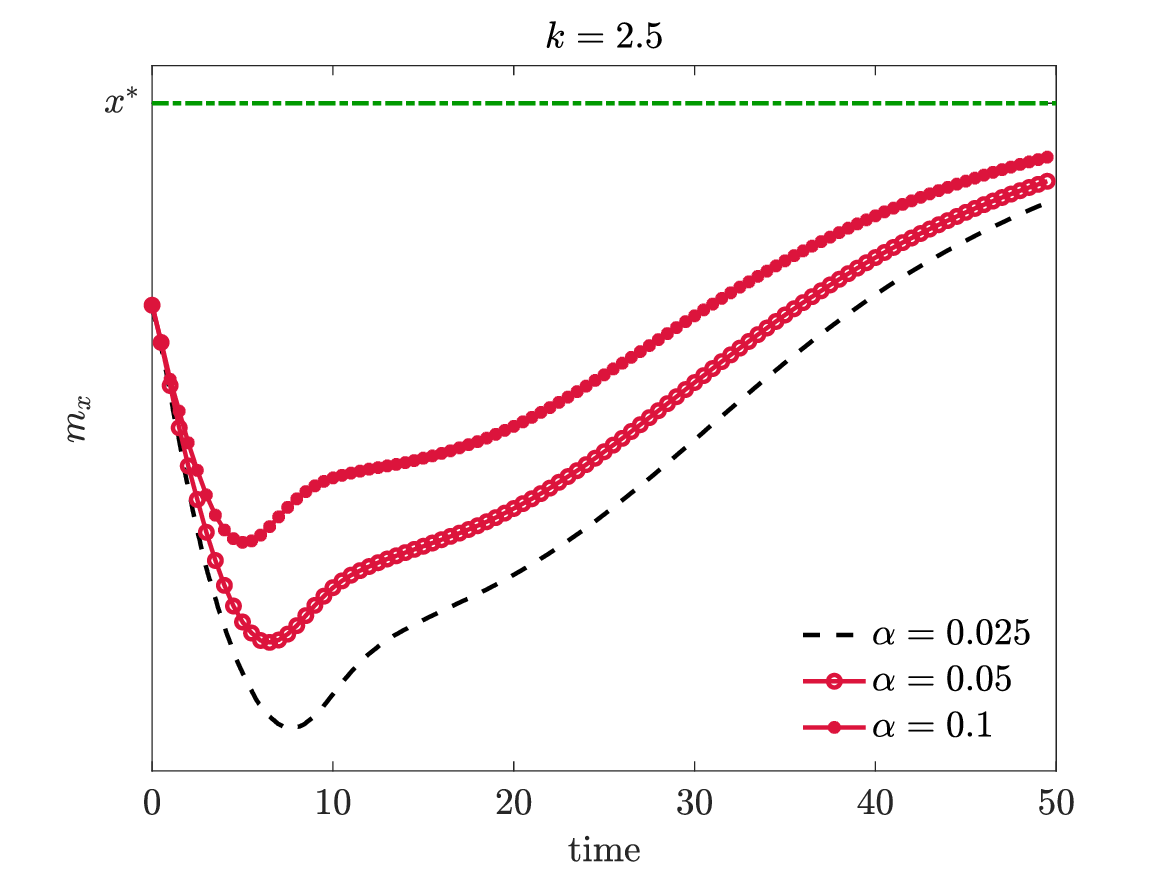}
		\includegraphics[scale = 0.35]{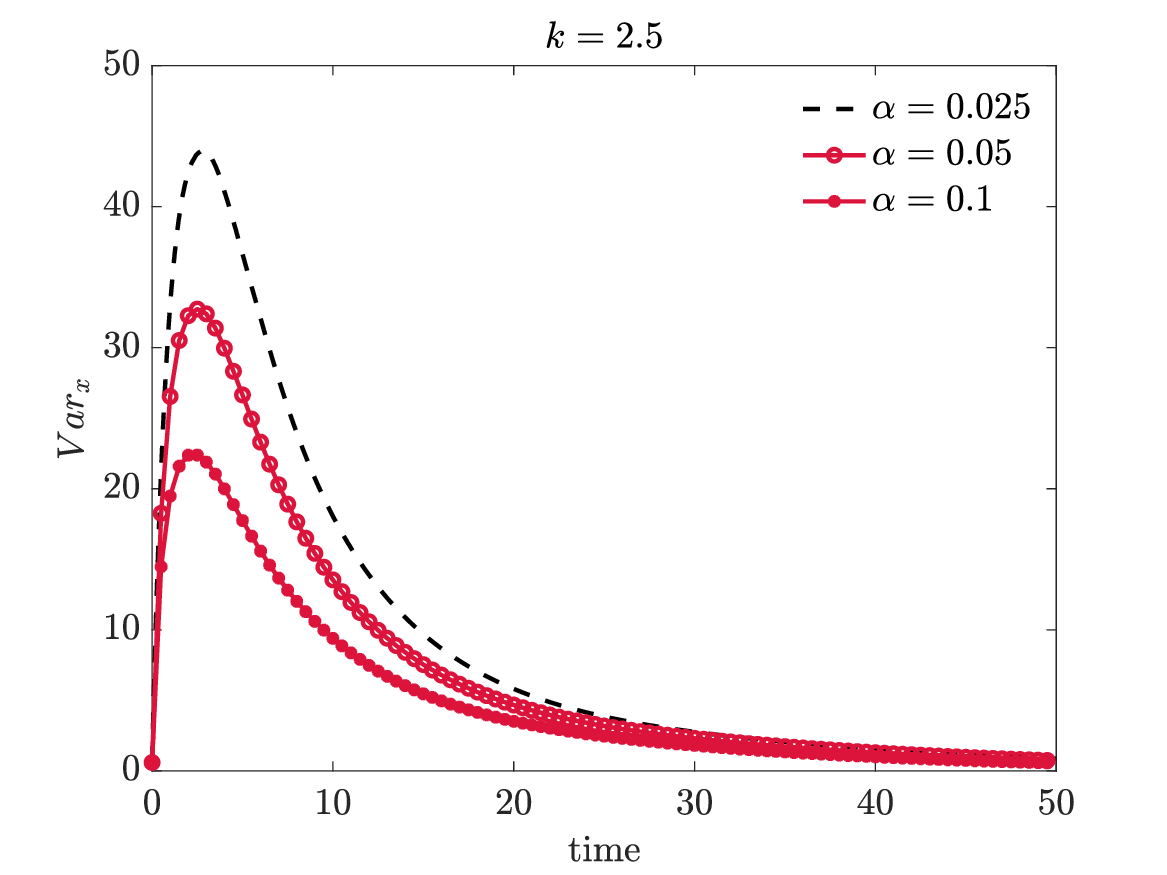}
		\caption{Evolution of the mean position $m_x =\int_{\mathbb  R}xf(x,t)dx$ and of its variance $Var_x =\int_{\mathbb R}(x-m_x)^2f(x,t)dx$ over the time interval $[0,50]$ and several values of the parameter $\alpha= 0.025,0.05,0.1$ in the case $k = 2.5$. }
		\label{fig:mx_vx_k}
	\end{figure}

	In Figure \ref{fig:mx_vx_k} we report the evolution of mean and variance of the kinetic density $f(x,t)$ over the time interval $[0,50]$ from which we may observe how it approaches the global minimum of $\F(x)$ while asymptotically reducing its variance. To highlight this behavior we report in Figure \ref{fig:ftime_k>1} the reconstructed kinetic density of particles $f(x,t)$ at times $t = 1,5,10$ and for several choices of $\kappa  \in \{ 1.5,2.0,2.5 \}$ and $\alpha \in \{  0.025,0.05,0.1$ \}, respectively. We  observe that  the distribution has its maximum in the global minimum of the cost as expected. 
	
	\begin{figure}
		\centering
		\includegraphics[scale = 0.225]{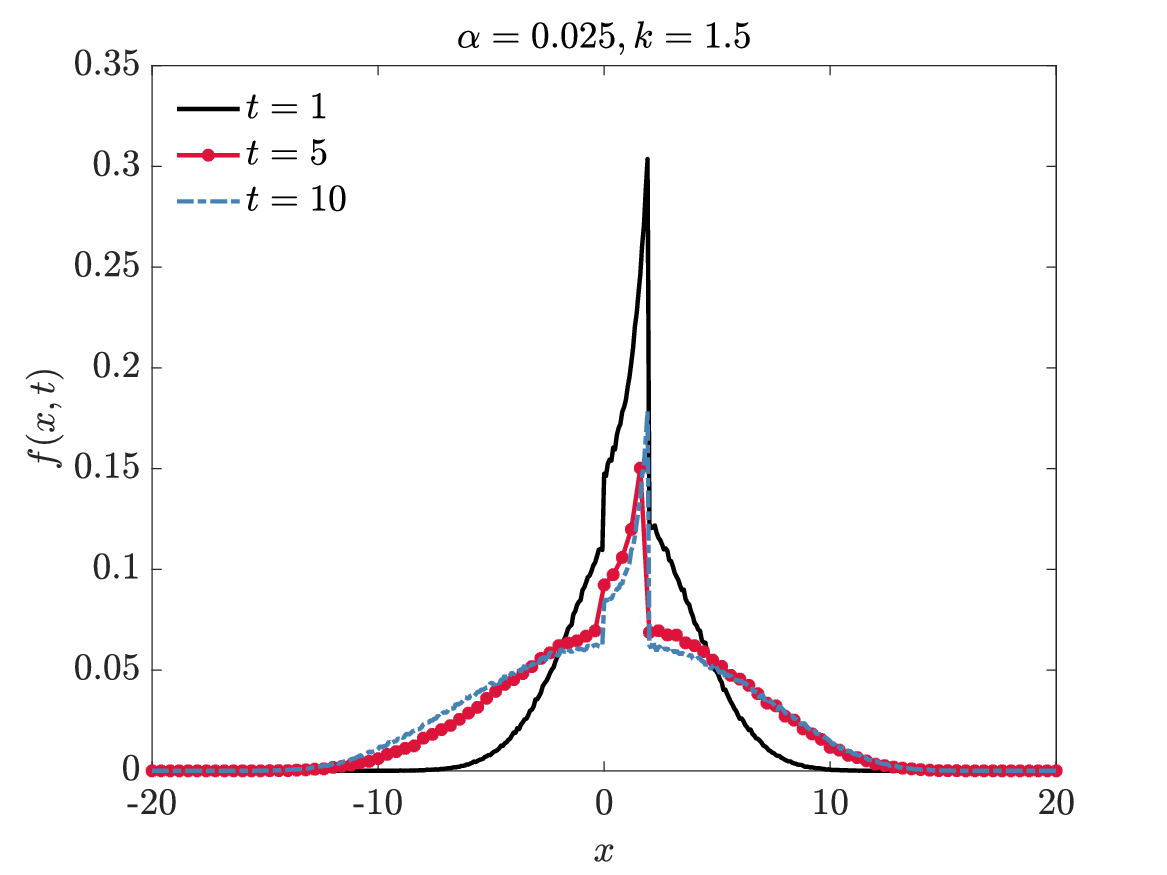}
		\includegraphics[scale = 0.225]{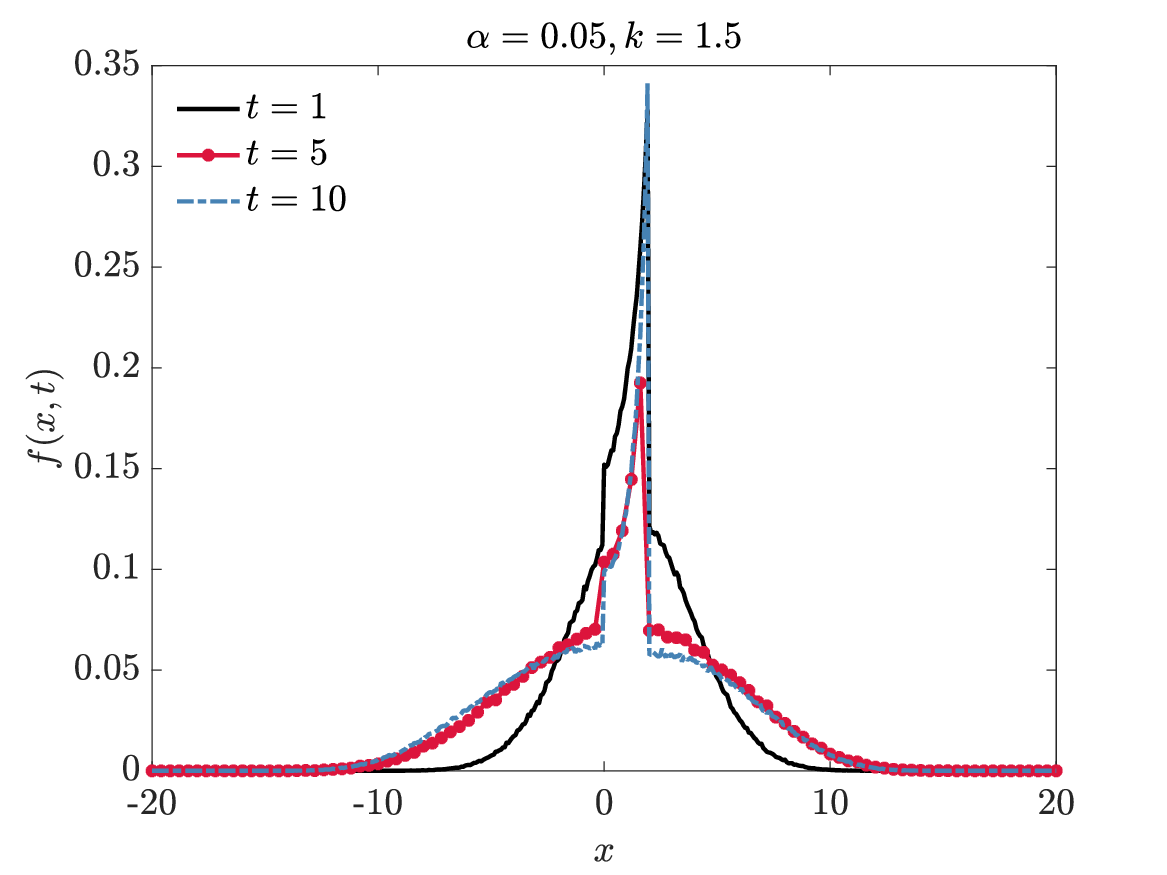}
		\includegraphics[scale = 0.225]{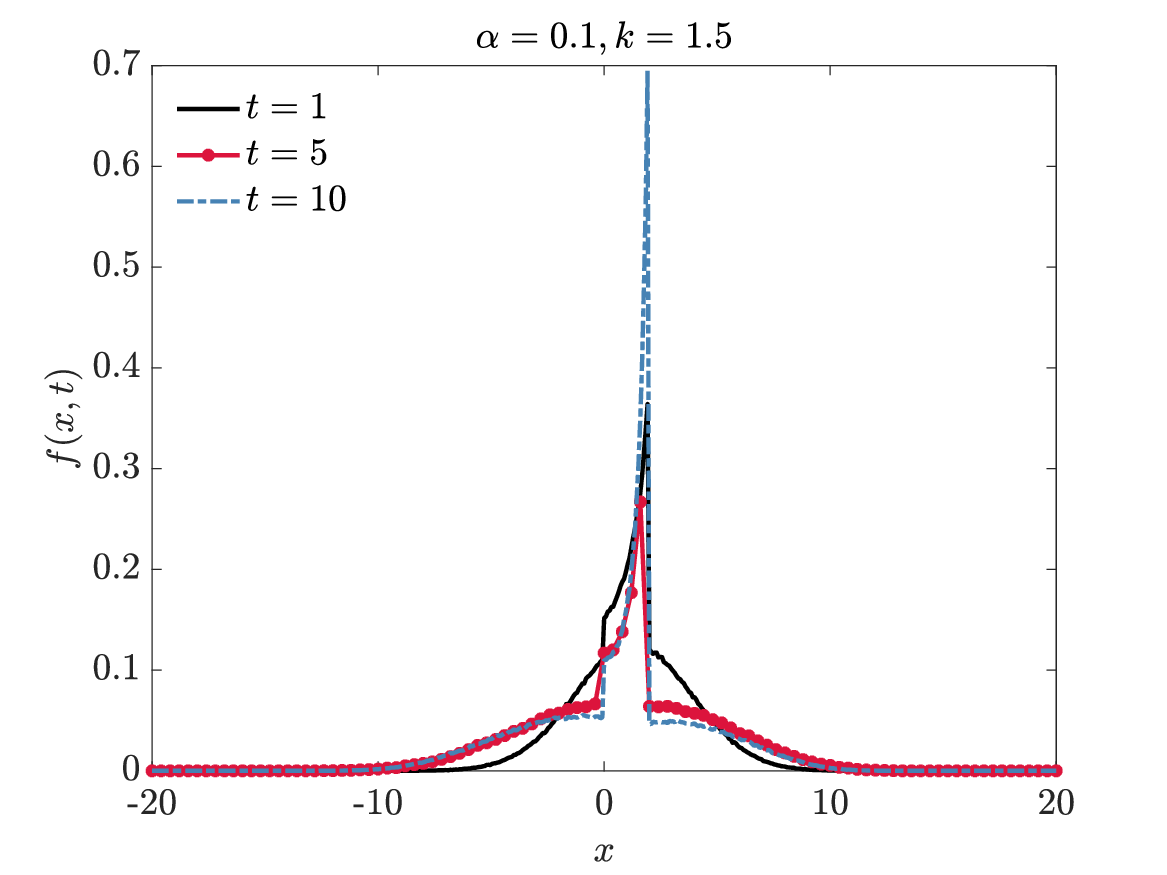}\\
		\includegraphics[scale = 0.225]{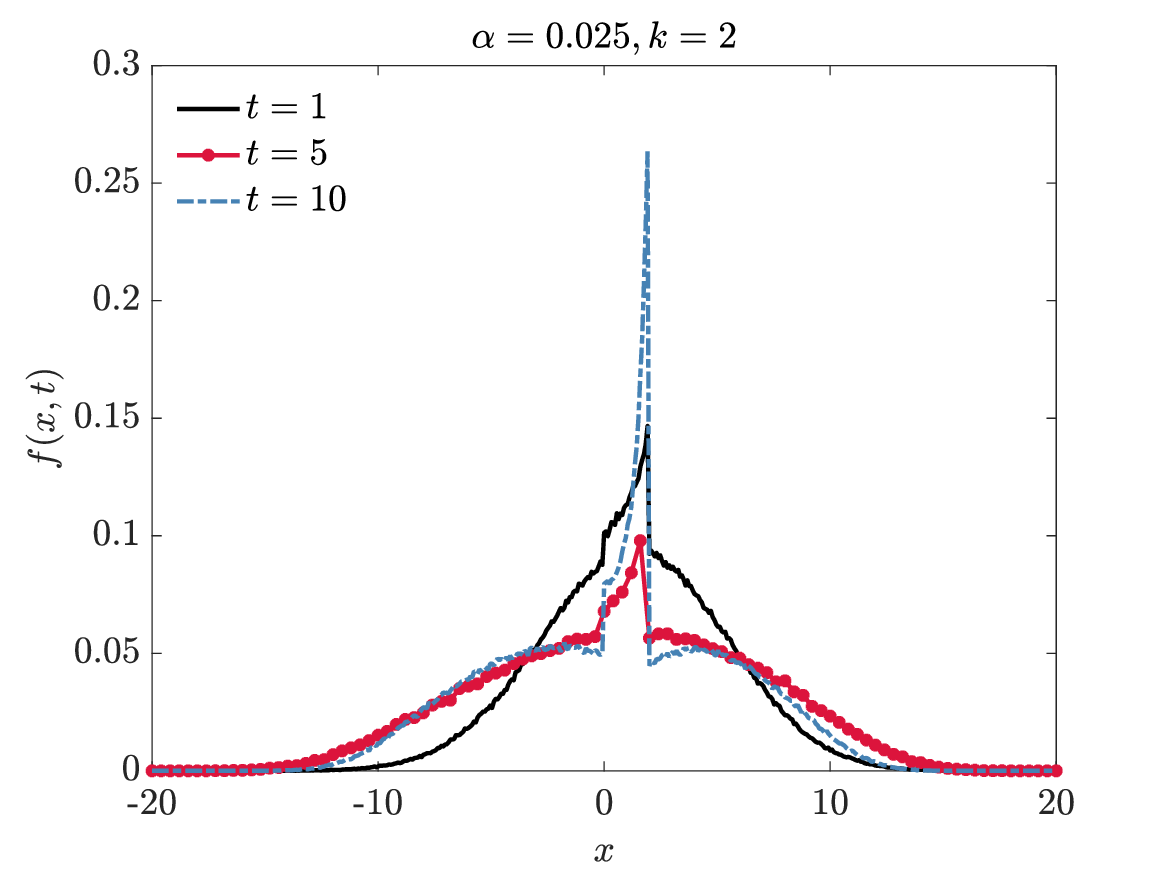}
		\includegraphics[scale = 0.225]{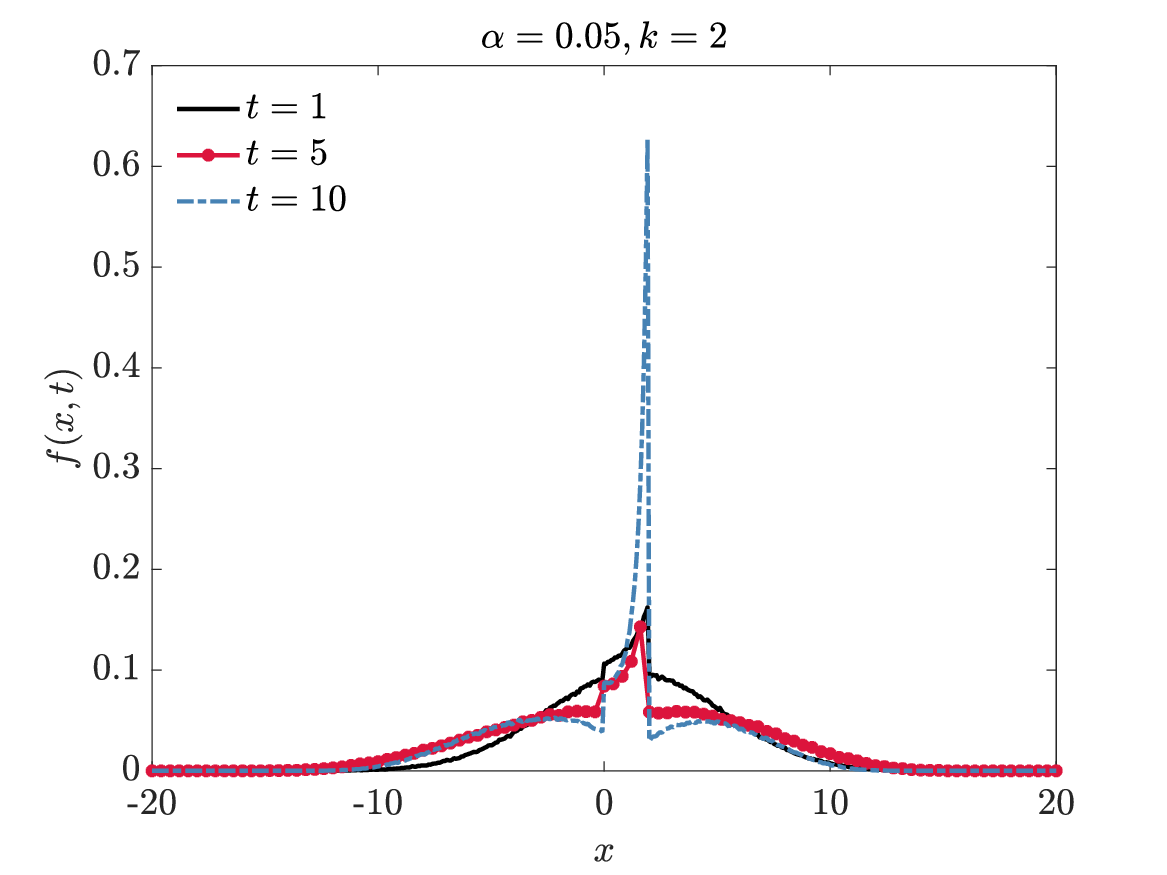}
		\includegraphics[scale = 0.225]{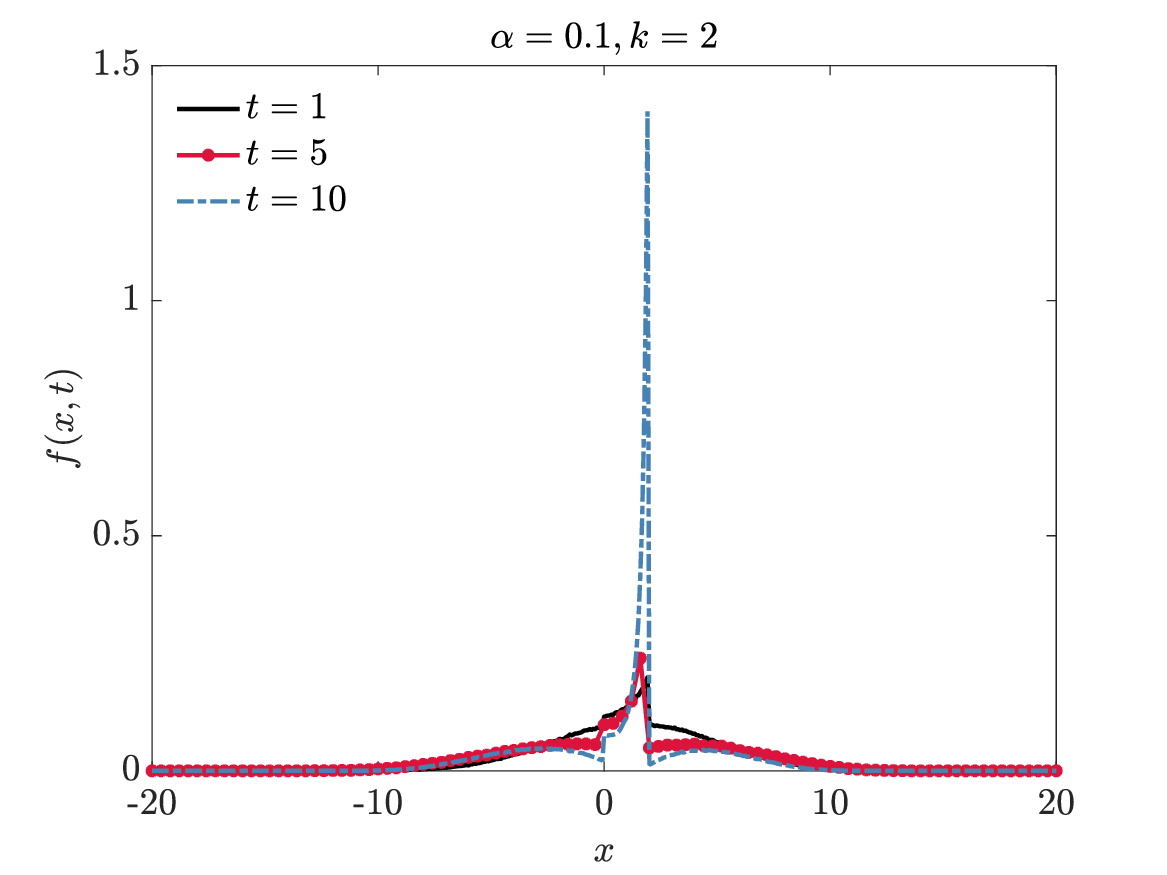}\\
		\includegraphics[scale = 0.225]{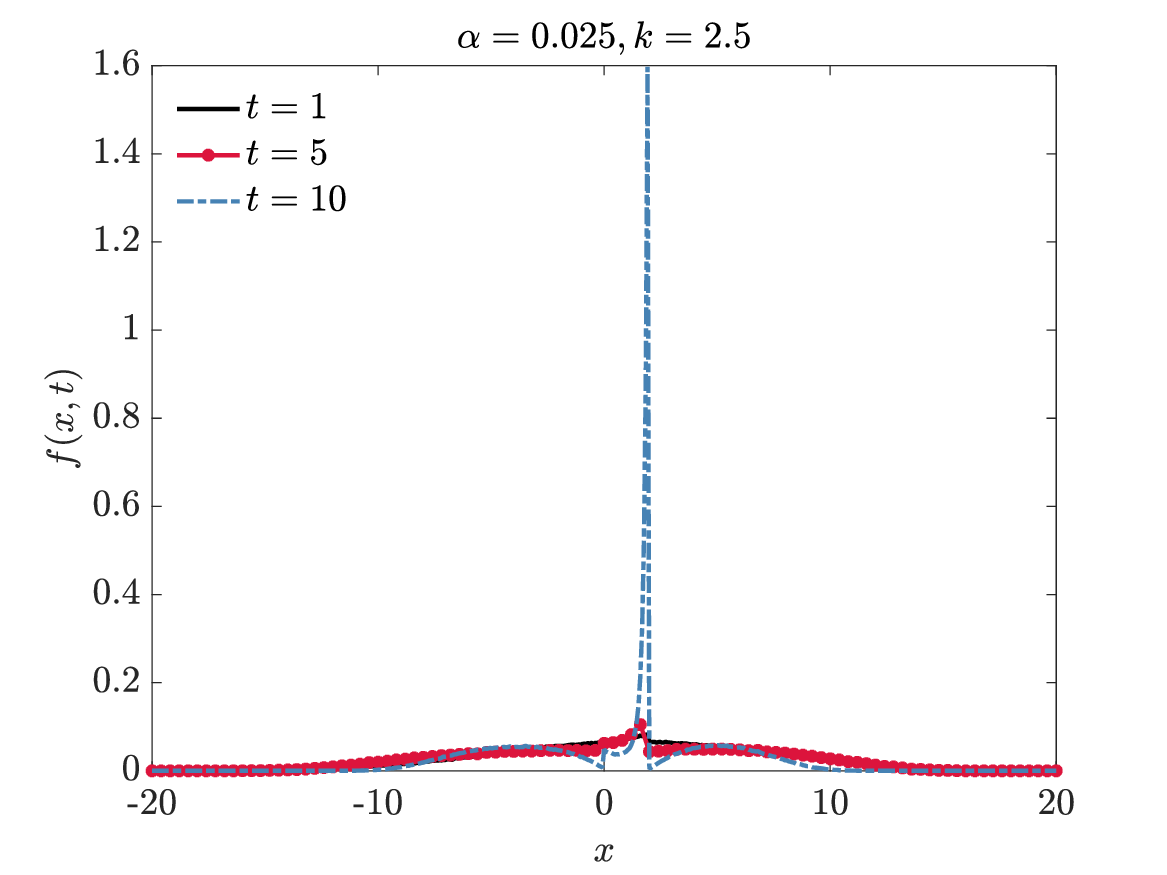}
		\includegraphics[scale = 0.225]{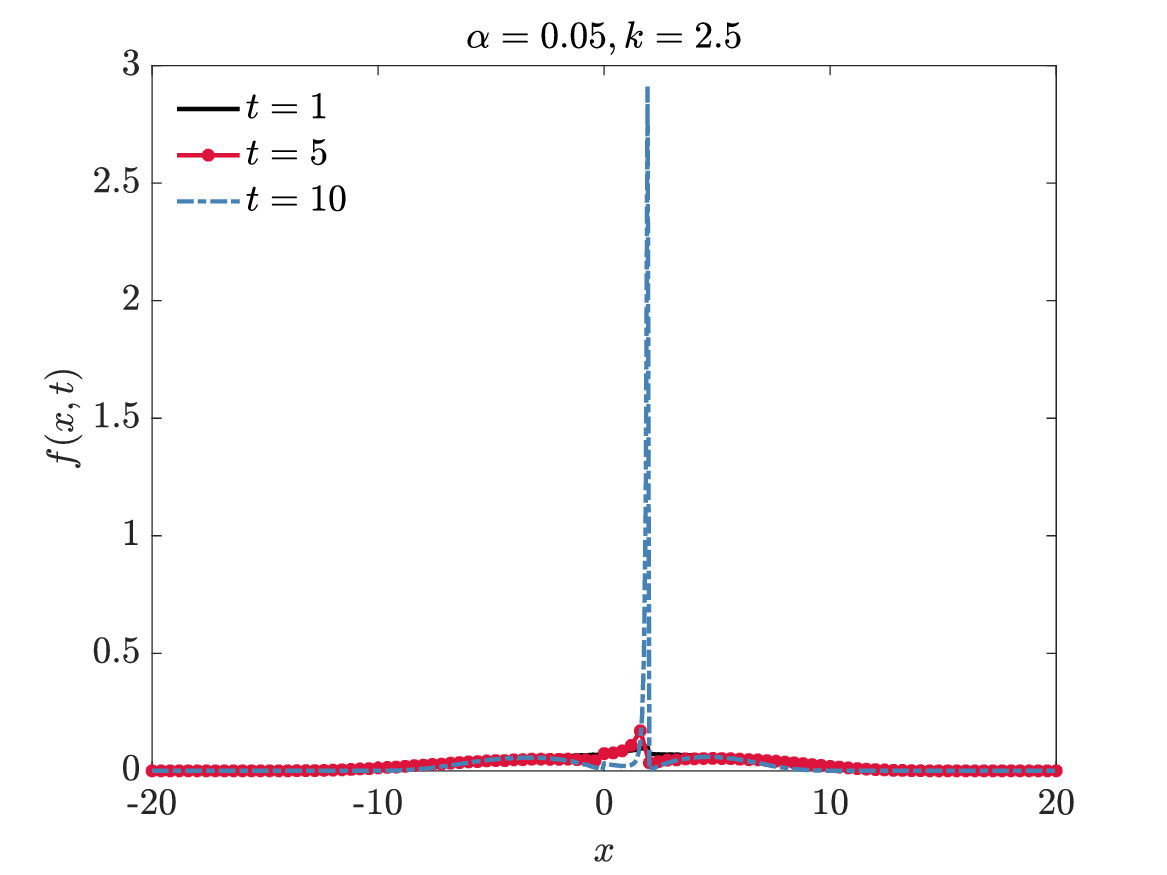}
		\includegraphics[scale = 0.225]{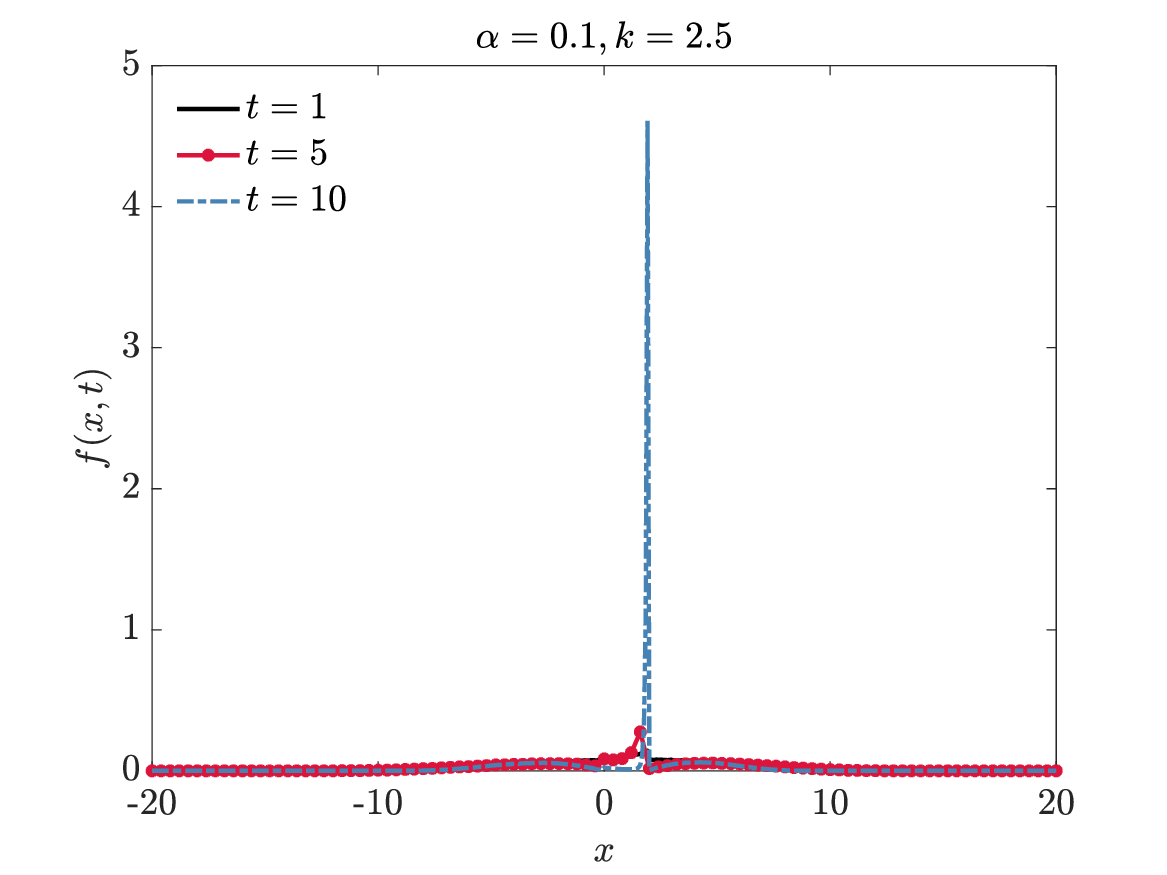}
		\caption{Evolution of the distribution $f(x,t)$ at time $t = 1,5,10$ for several values of $\alpha = 0.025$ (left), $\alpha = 0.05$ (center), $\alpha = 0.1$ (right), and $k = 1.5$ (first row), $k = 2.0$ (second row), $k = 2.5$ (third row). In all the tests we fixed $N = 10^6$, $p = 1/4$, $\theta = 0.5$ and $\sigma = 0.1$, the definition of $\lambda[f](\cdot)$ is \eqref{eq:lambda_kappa>1} and the initial state is \eqref{eq:f0T0_k>1}.  }
		\label{fig:ftime_k>1}
	\end{figure}
	
	\begin{figure}
		\centering
		\includegraphics[scale = 0.225]{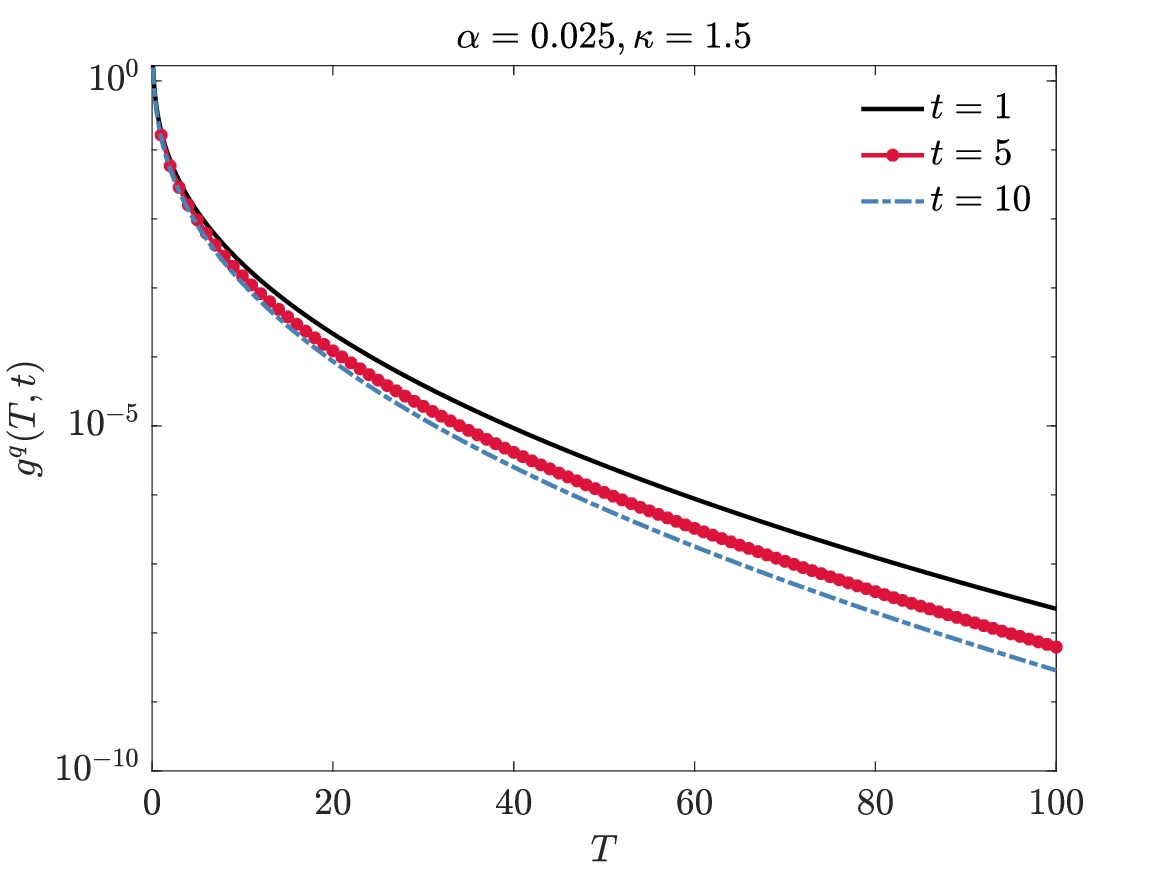}
		\includegraphics[scale = 0.225]{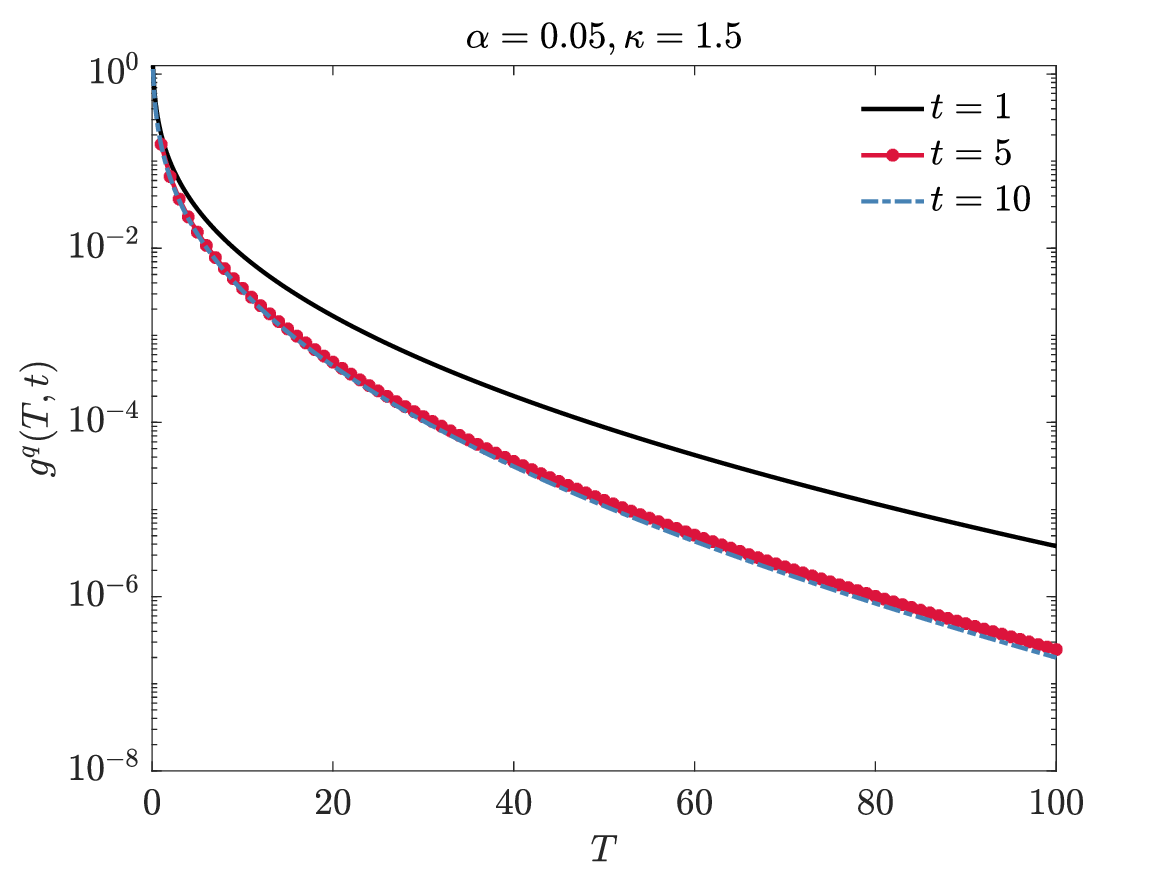}
		\includegraphics[scale = 0.225]{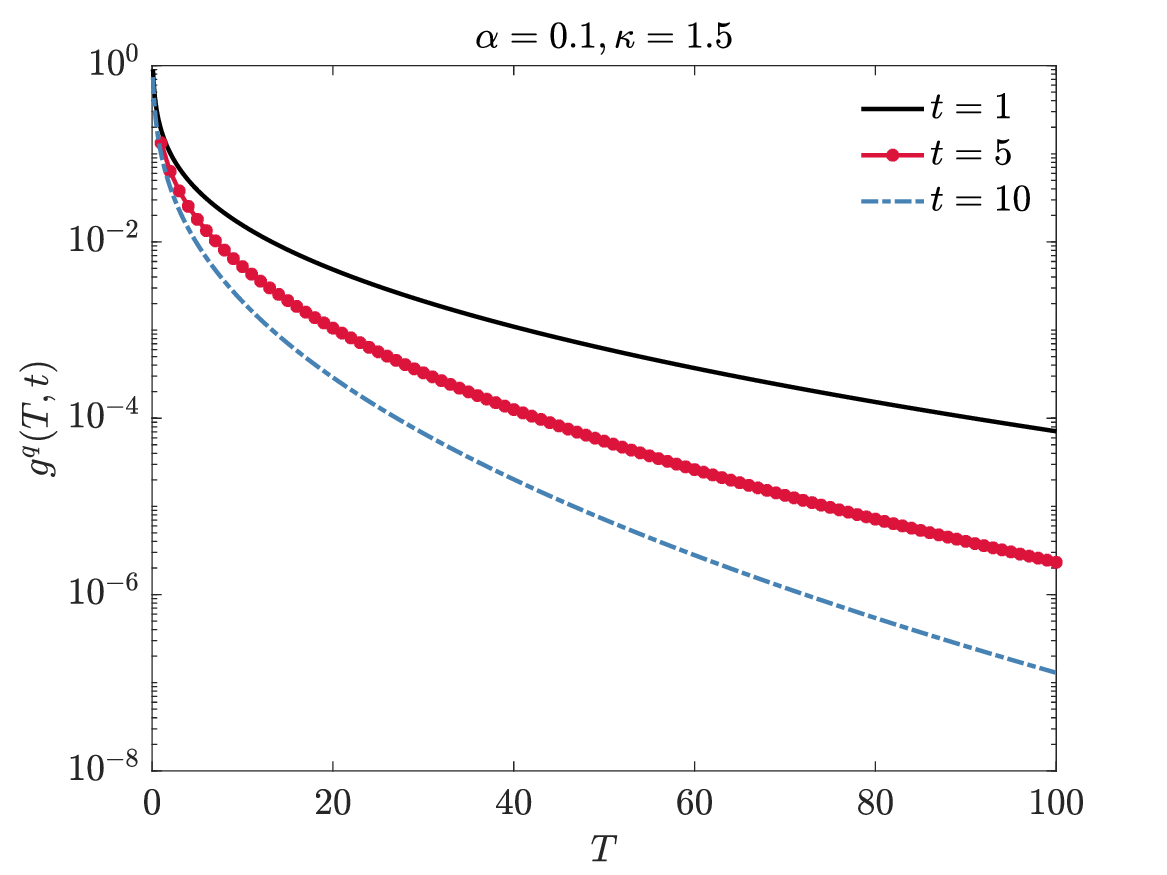}\\
		\includegraphics[scale = 0.225]{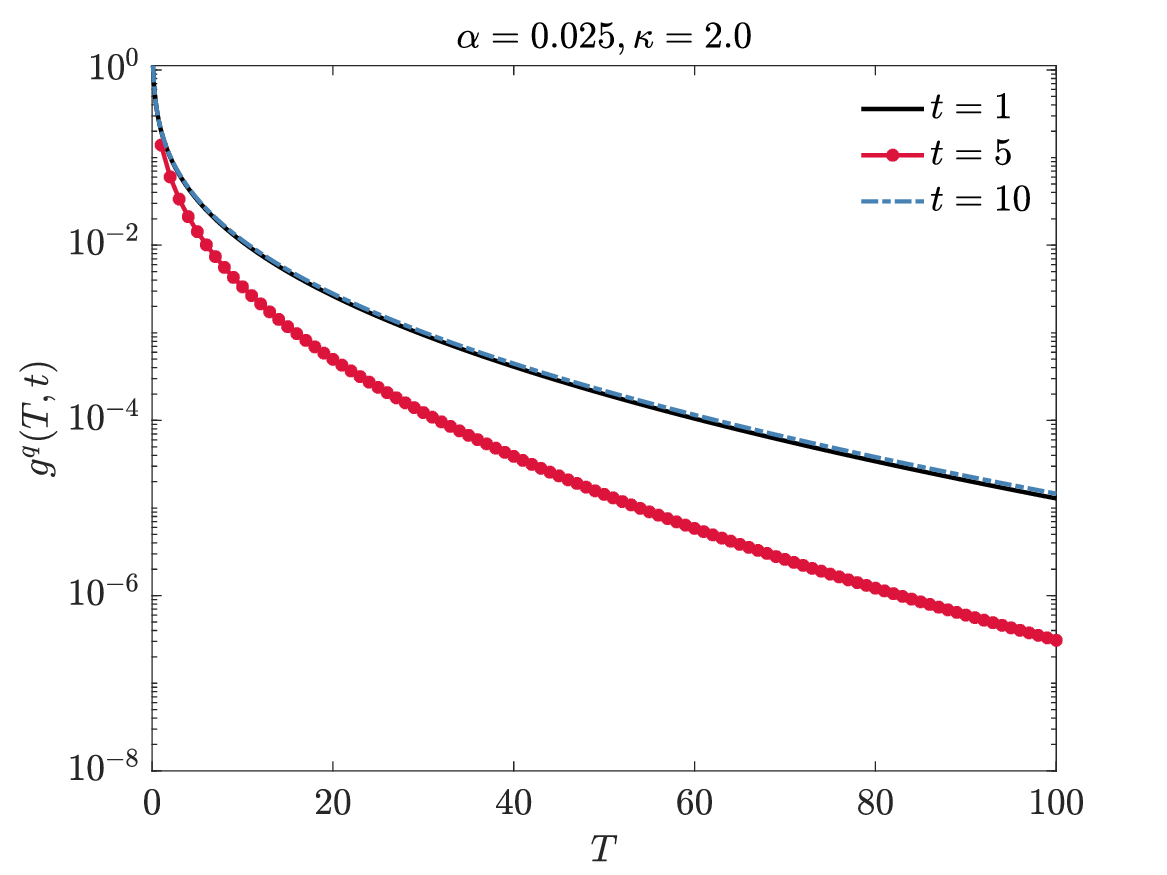}
		\includegraphics[scale = 0.225]{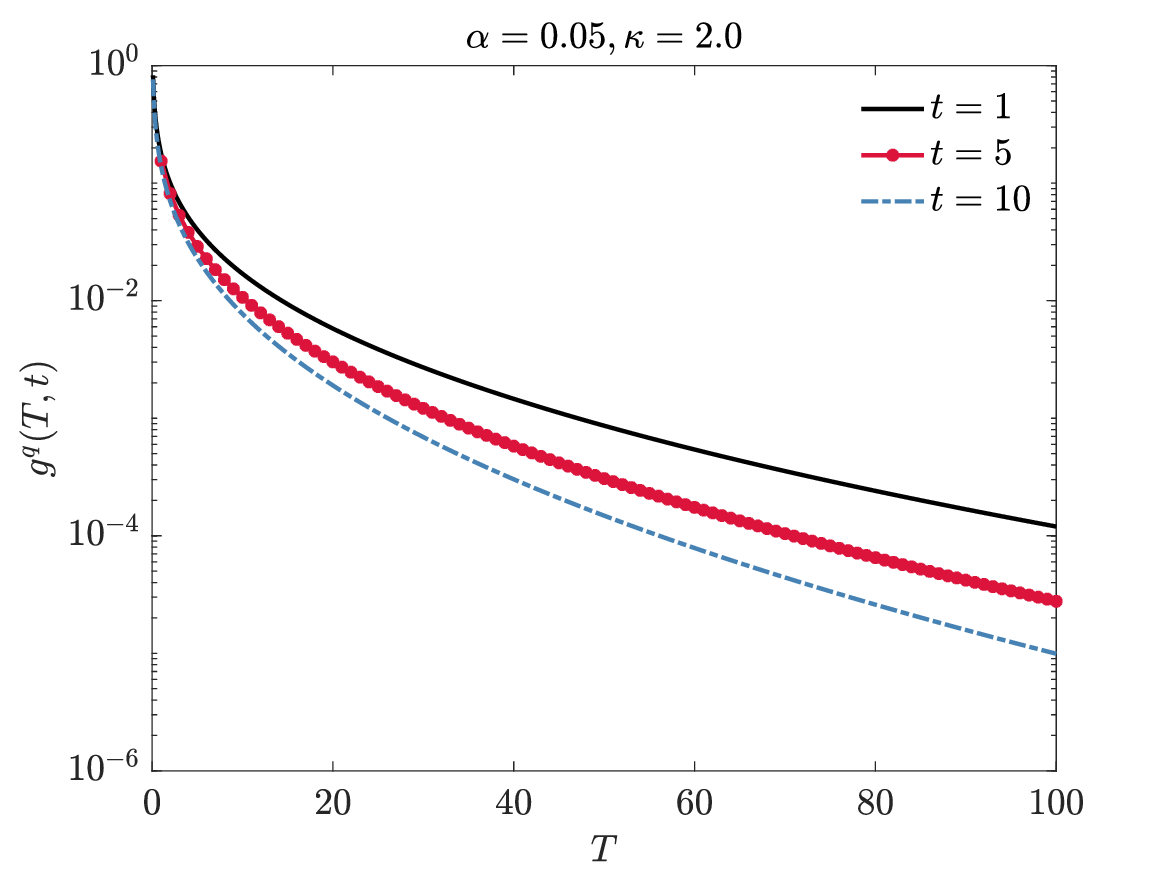}
		\includegraphics[scale = 0.225]{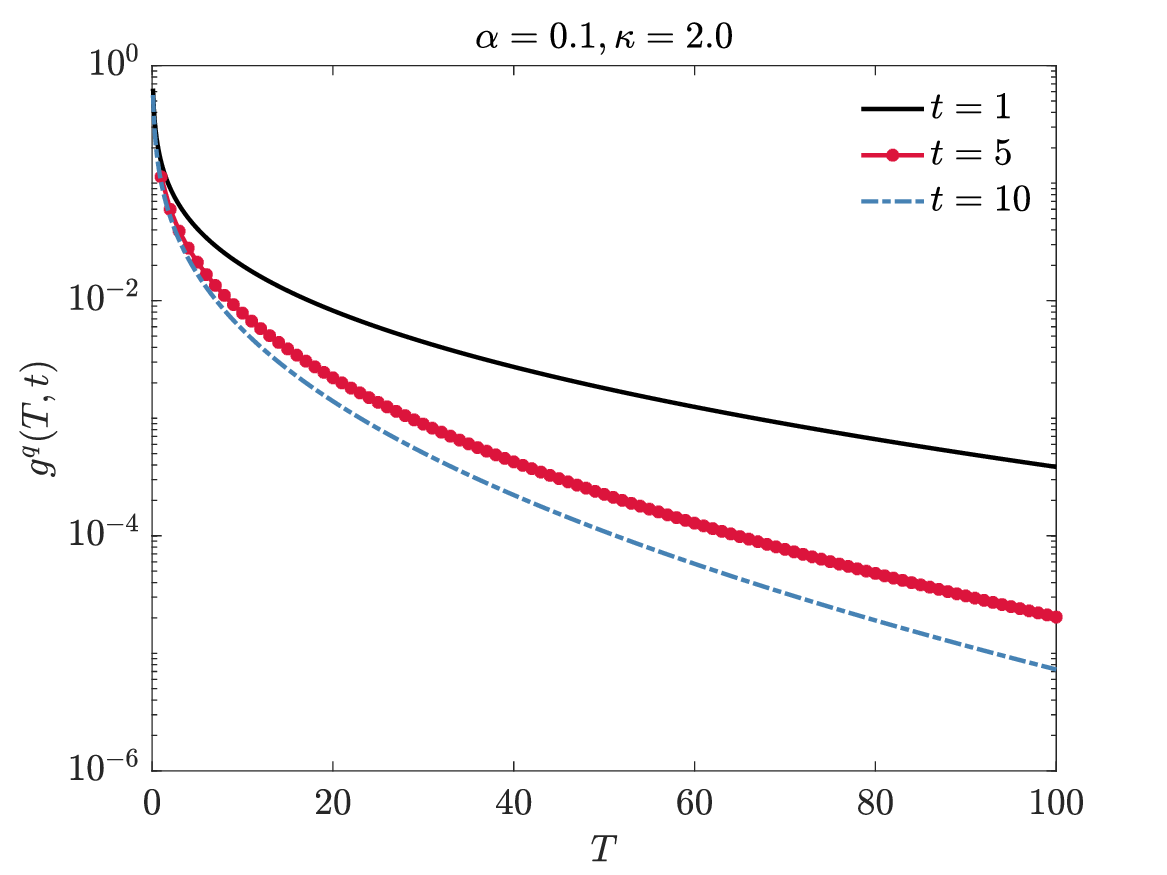}\\
		\includegraphics[scale = 0.225]{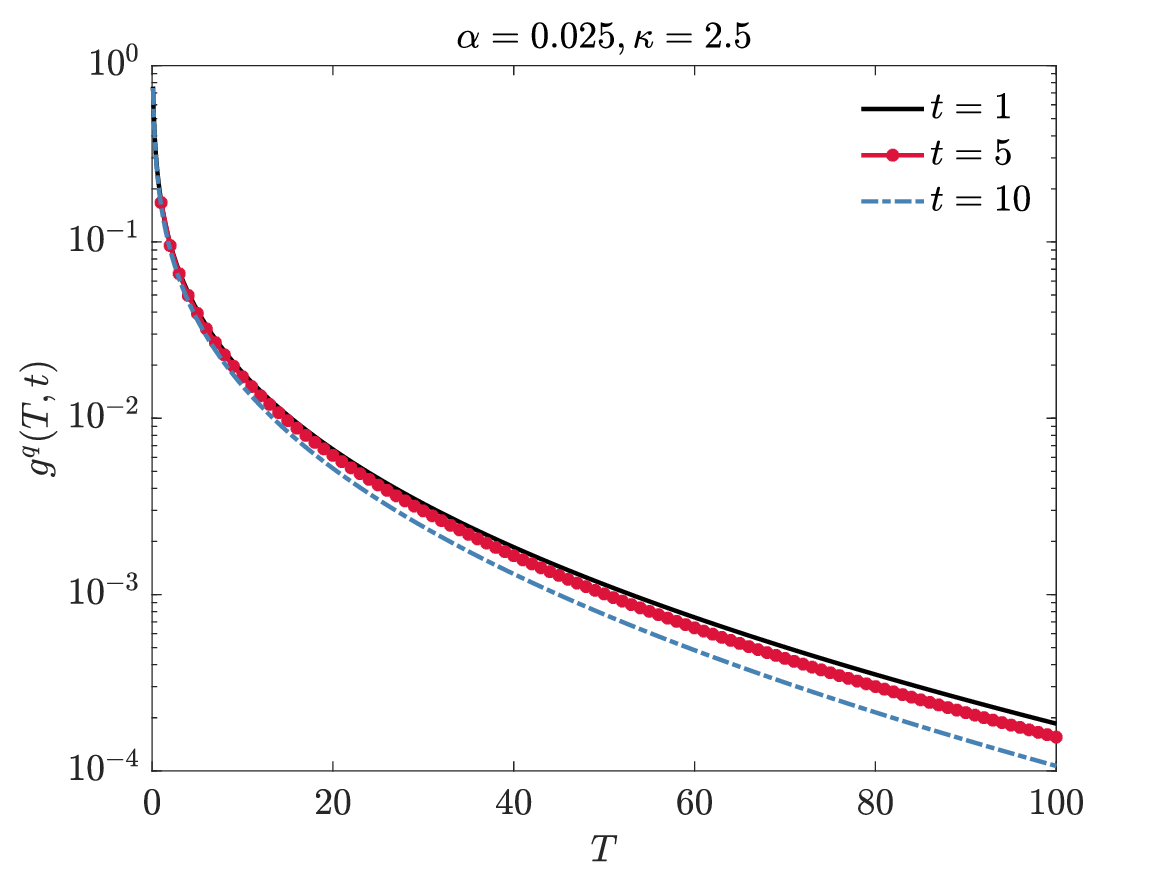}
		\includegraphics[scale = 0.225]{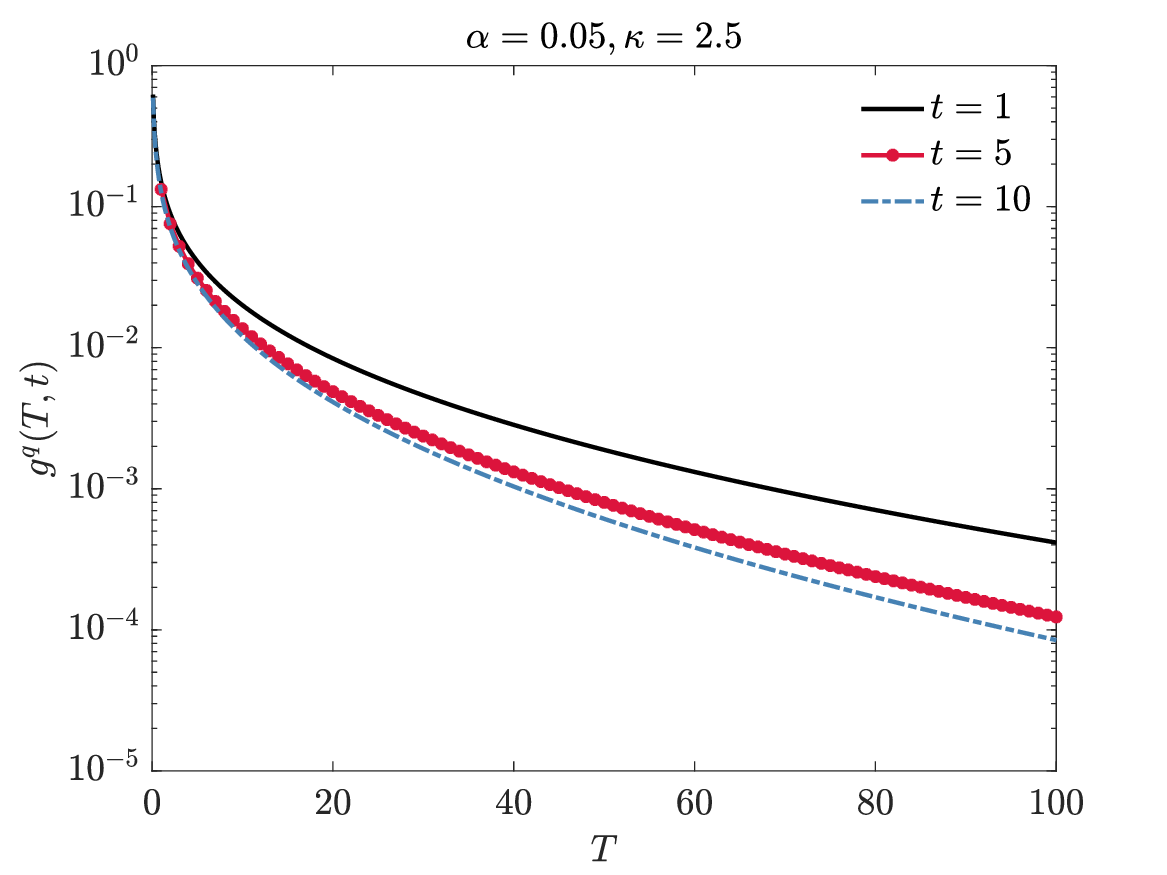}
		\includegraphics[scale = 0.225]{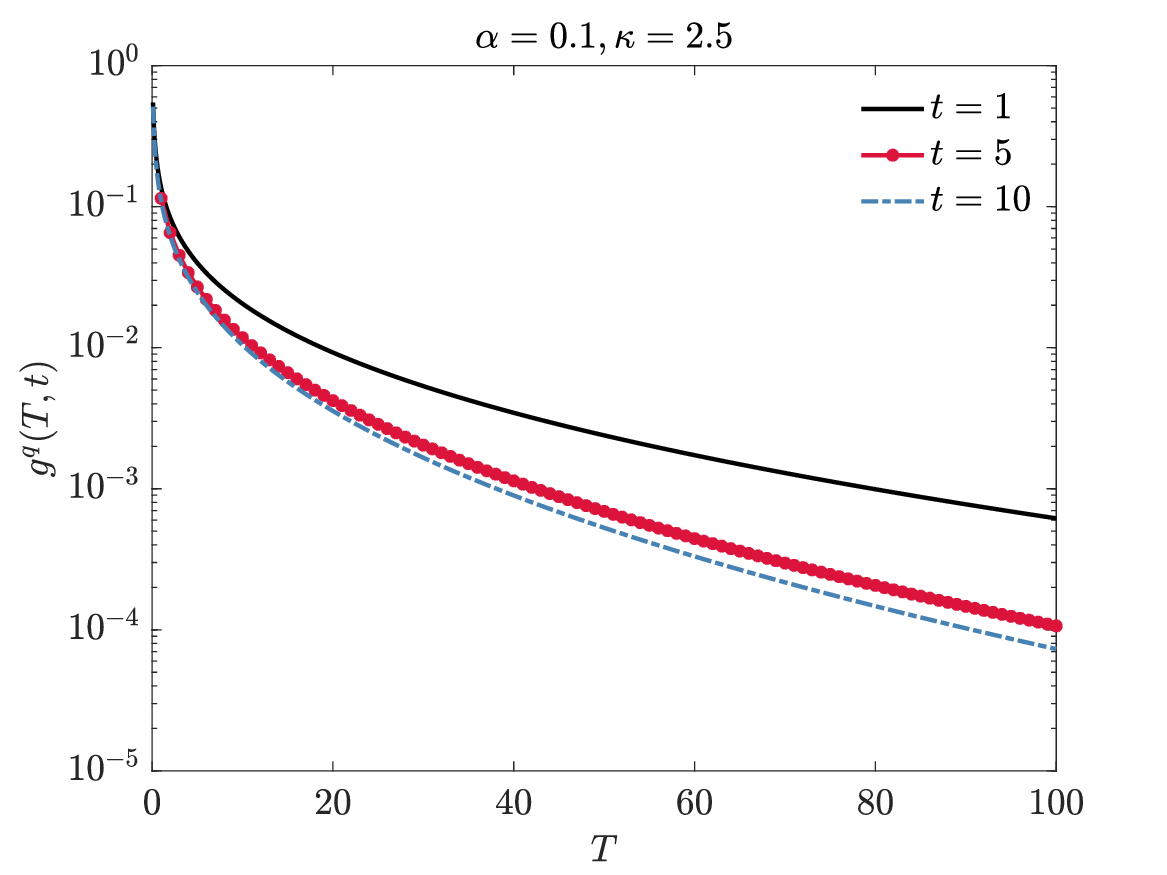}
		\caption{Evolution of the quasi-equilibrium distribution of the temperature $g^q(T,t)$ at time $t = 1,5,10$ for several values of $\alpha = 0.025$ (left), $\alpha = 0.05$ (center), $\alpha = 0.1$ (right), and $k = 1.5$ (first row), $k = 2.0$ (second row), $k = 2.5$ (third row). We computed generalized gamma distribution defined in \eqref{eq:gen_gamma} where the definition of $\lambda[f](\cdot)$ is \eqref{eq:lambda_kappa>1}.  }
		\label{fig:gq_k>1}
	\end{figure}

	\begin{figure}
		\centering
		\includegraphics[scale = 0.225]{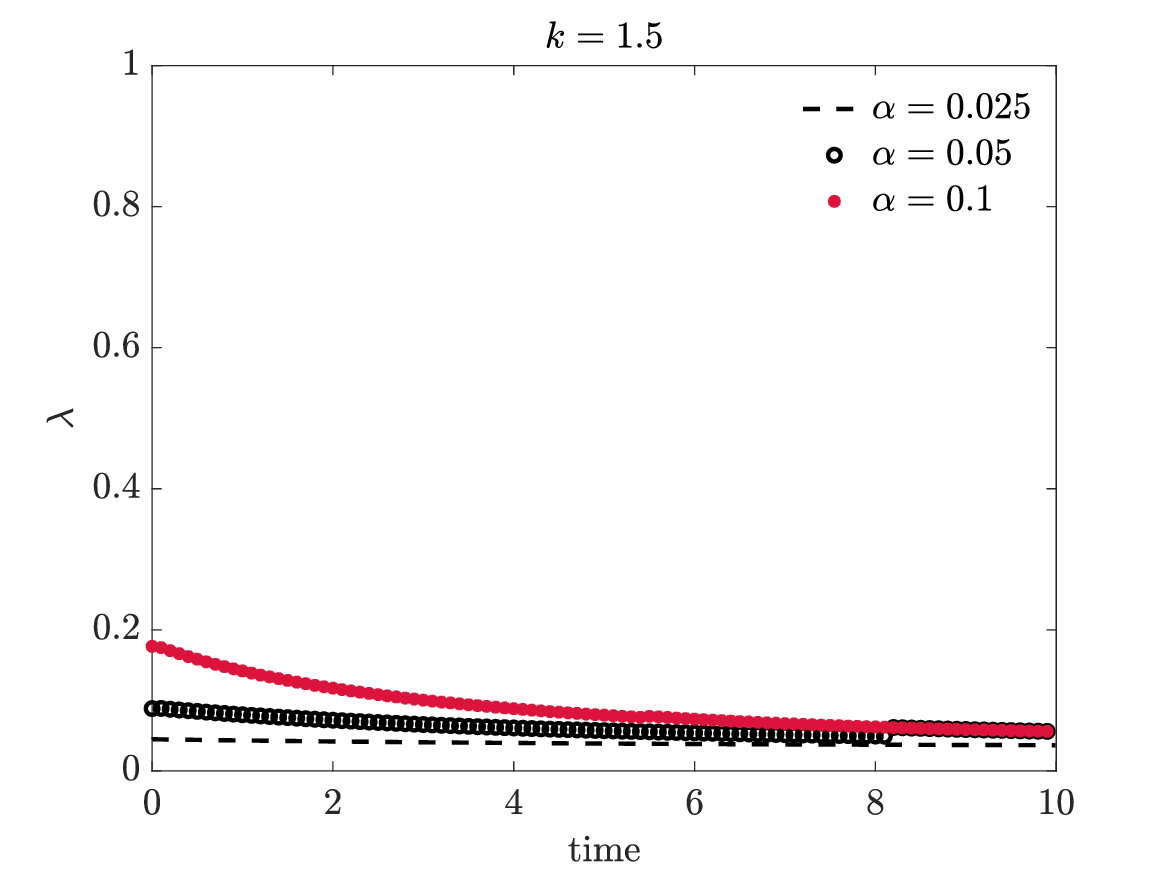}
		\includegraphics[scale = 0.225]{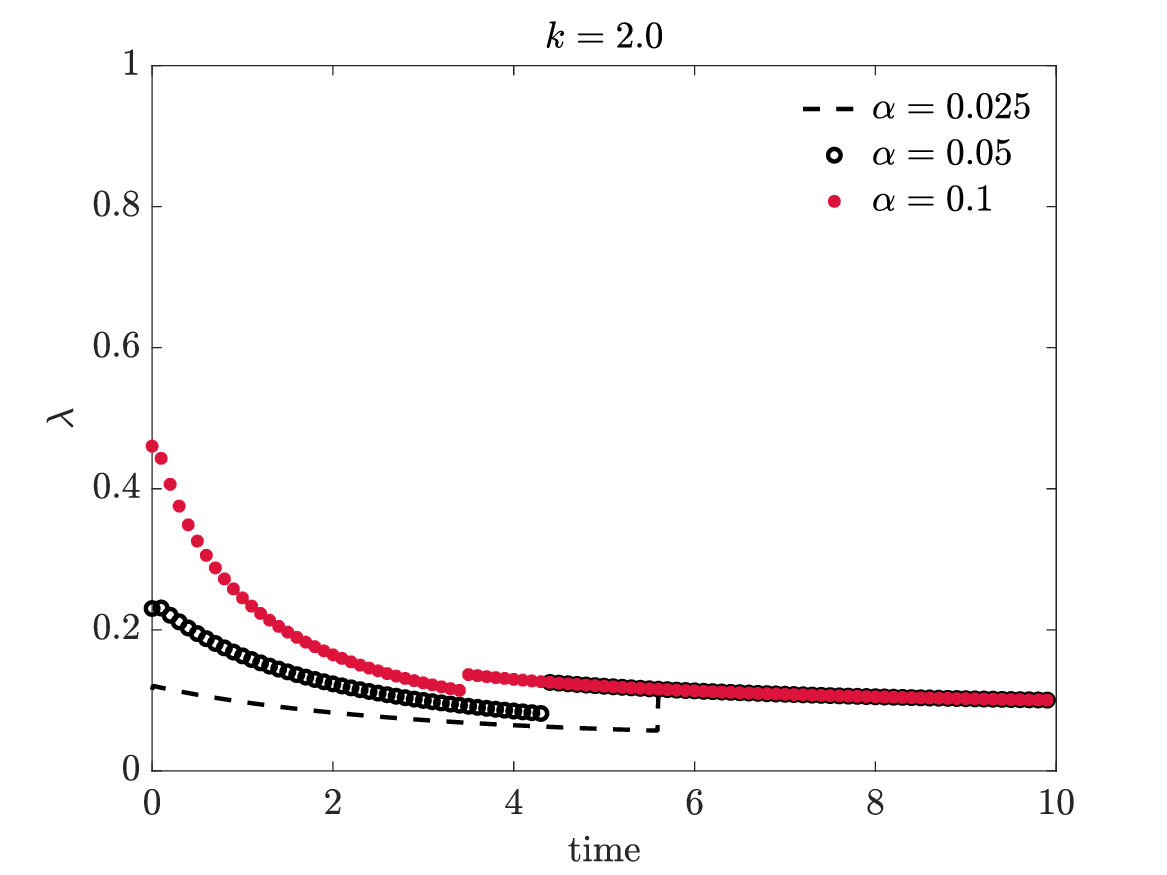}
		\includegraphics[scale = 0.225]{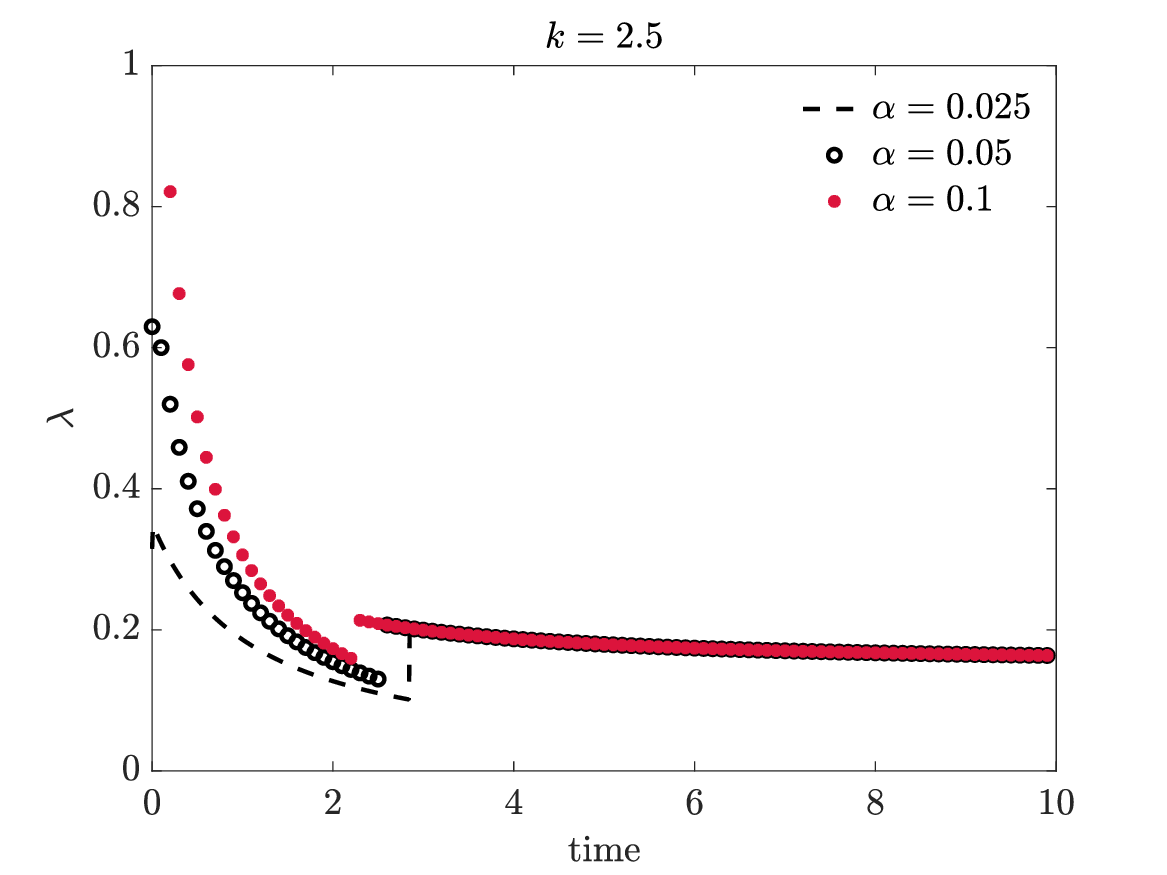} \\
		\includegraphics[scale = 0.225]{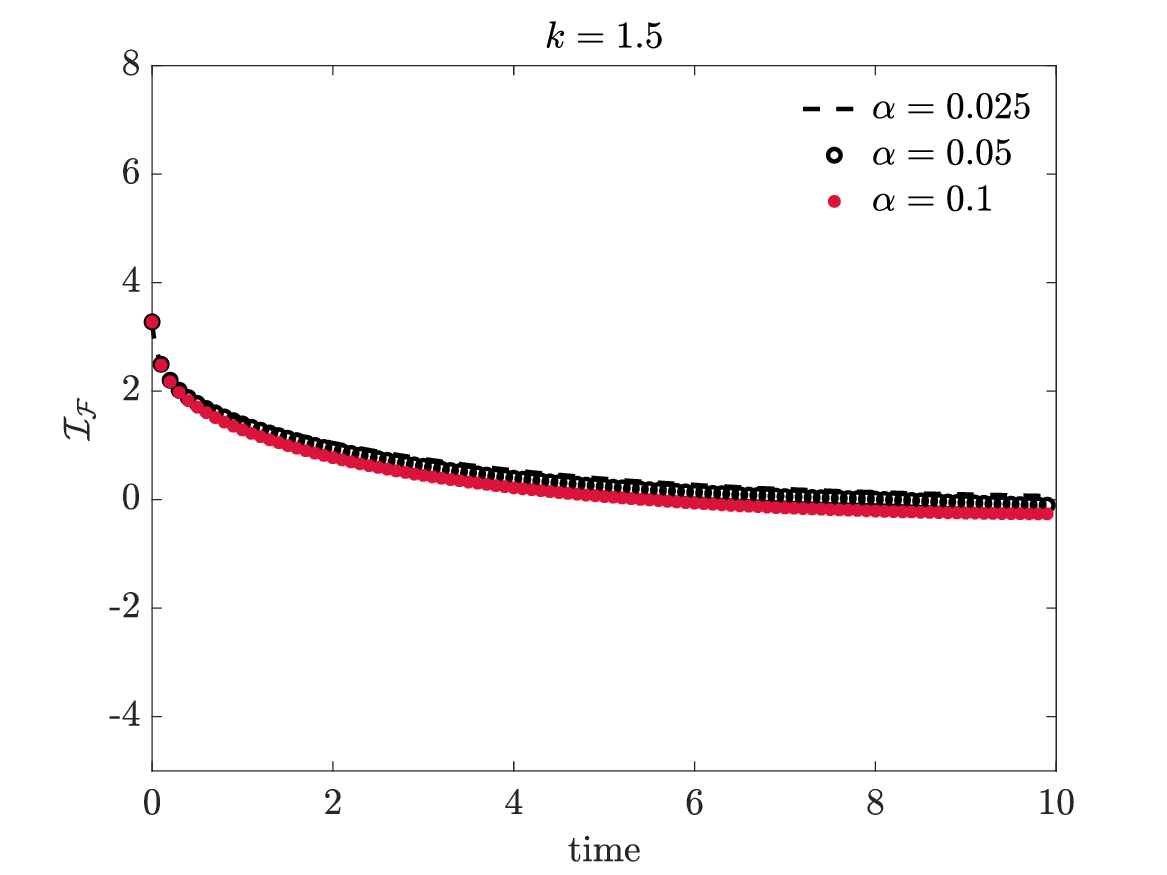}
		\includegraphics[scale = 0.225]{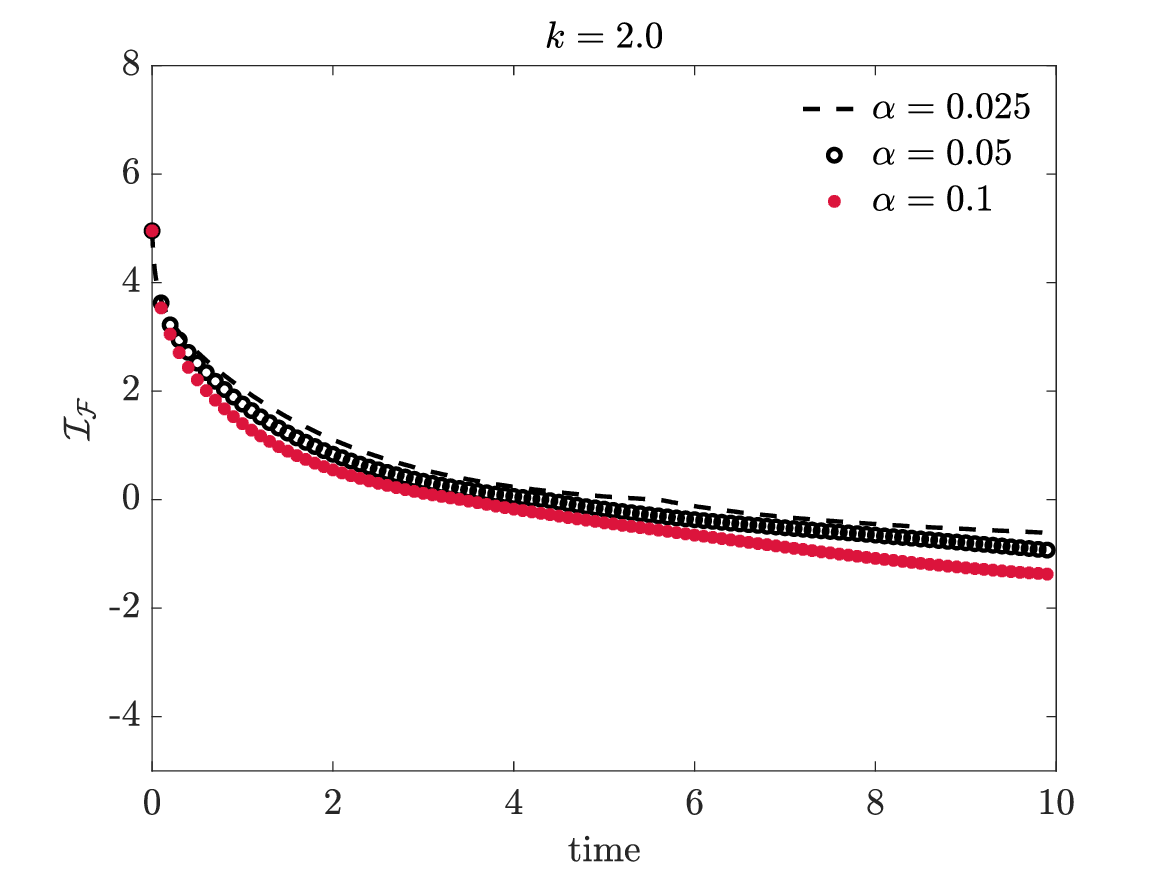}
		\includegraphics[scale = 0.225]{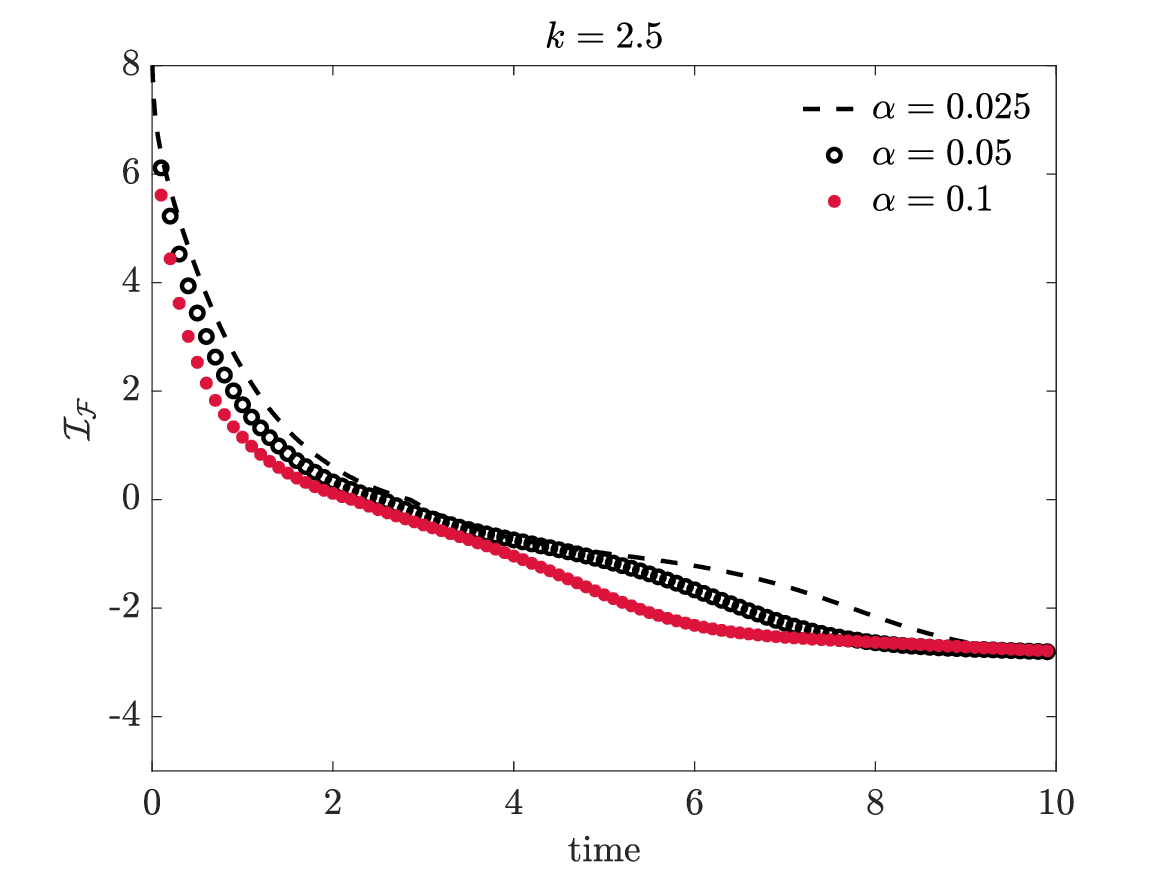} 
		\caption{ Top row: estimated values of $\lambda>0$ as in \eqref{eq:lambda_kappa>1} for several values of $\alpha = 0.025,0.05,0.1$ and $\kappa = 1.5$ (left), $\kappa = 2.0$ (center), $\kappa = 2.5$ (right). Bottom row: estimated values of $\mathcal I_\F(t)$. We considered $N = 10^6$ particles both in space and temperature, $p= 1/4$, $\theta = 0.05$ and $\sigma^2 = 0.1$. The initial condition is reported in \eqref{eq:f0T0_k>1}. }
		\label{fig:ftime_k>1}
	\end{figure}
	
	\section*{Conclusion}
In this paper, we propose an extended version of the simulated annealing method in which the temperature is dynamically controlled through a state-dependent feedback mechanism. Specifically, we introduce a feedback law that, under suitable assumptions, guarantees an exponential decay of the system's entropy. This approach leads to an efficient exploration of the energy landscape by dynamically adapting the cooling schedule in response to the current state of the system.

The analysis is carried out by reformulating the underlying particle dynamics in terms of binary interactions, which allows us to derive a Boltzmann-type kinetic equation for the evolution of the probability density. Entropy decay estimates are rigorously established in the quasi-invariant scaling, highlighting the effectiveness of the proposed control strategy. {The control strategy exploits the decay of the temperature. This includes the use of the self-similar generalised gamma profiles. }

The numerical approximation of the resulting coupled Boltzmann dynamics has been investigated through a direct simulation Monte Carlo method for the update of the space and temperature variables. Numerical experiments confirm the theoretical findings, demonstrating the expected exponential decay of entropy obtained through a  prototypical nonconvex functional. Furthermore, an extensive numerical investigation is conducted to assess the role of the feedback parameter, and a comparative analysis with the classical simulated annealing method is presented. The results indicate a significant improvement in convergence behavior under the proposed controlled framework.

Future developments will focus on extending the analysis to settings with more general temperature dynamics, with the goal of understanding how these influence the entropy decay behavior of the algorithm.	

{
\section*{Data Availability Statement}
The codes developed for this study are publicly available on Zenodo. \footnote{https://doi.org/10.5281/zenodo.18173007}
}

	\section*{Acknowledgements}
	M.Z. acknowledges partial support of GNFM group of INdAM (National Institute of High Mathematics) and of PRIN2022PNRR project No.P2022Z7ZAJ, European Union - NextGenerationEU and by ICSC - Centro Nazionale di Ricerca in High Performance Computing, Big Data and Quantum Computing, funded by European Union - NextGenerationEU. MH acknowledges funding 
	by the  Deutsche Forschungsgemeinschaft (DFG, German Research Foundation) for the financial support through 320021702/GRK2326 and through 333849990/GRK2379
	(IRTG Hierarchical and Hybrid Approaches in Modern Inverse Problems). MH  received funding from the European Union’s Horizon Europe research and innovation programme under the Marie Sklodowska-Curie Doctoral Network Datahyking (Grant No. 101072546).

	\bibliographystyle{siam}
	\bibliography{refs}

@unpublished{PFZ_25,
	author = {Pareschi, L. and Franceschi, J. and Zanella, M. },
	note = {Preprint arXiv:2506.09001},
	title = {Superlinear drift in consensus-based optimization with condensation phenomena},
	year = {2025}}

@article{ingber89,
	author = {Ingber, L.},
	journal = {Math. Comput. Model.},
	number = {8},
	pages = {967--973},
	title = {Very fast simulated re-annealing},
	volume = {12},
	year = {1989}}

@article{ingber93,
	author = {Ingber, L.},
	journal = {Math. Comput. Model.},
	number = {11},
	pages = {29--57},
	title = {Simulated annealing: {P}ractice versus theory},
	volume = {18},
	year = {1993}}

@article{parisi92,
	author = {Marinari, E. and Parisi, G.},
	journal = {Europhys. Lett.},
	number = {6},
	title = {Simulated Tempering: {A} New {M}onte {C}arlo Scheme},
	volume = {19},
	year = {1992}}

@article{BMTZ_25,
	author = {Bondesan, A. and Menale, M. and Toscani, G. and Zanella, M.},
	journal = {Nonlinearity},
	pages = {075026},
	title = {Lotka-{V}olterra-type kinetic equations for interacting species},
	volume = {38},
	year = {2025}}

@article{FPTZ_17,
	author = {Furioli, G. and Pulvirenti, A. and Terraneo, E. and Toscani, G.},
	journal = {Math. Mod. Meth. Appl. Sci.},
	number = {01},
	pages = {115--158},
	title = {Fokker-{P}lanck equations in the modeling of socio-economic phenomena},
	volume = {27},
	year = {2017}}

@article{PZ18,
	author = {Pareschi, Lorenzo and Zanella, Mattia},
	journal = {J. Sci. Comput.},
	number = {3},
	pages = {1575--1600},
	title = {Structure preserving schemes for nonlinear {F}okker-{P}lanck equations and applications},
	volume = {74},
	year = {2018}}

@book{grippo2023introduction,
	author = {Grippo, Luigi and Sciandrone, Marco},
	publisher = {Springer Nature},
	title = {Introduction to Methods for Nonlinear Optimization},
	volume = {152},
	year = {2023}}

@book{pareschi2013interacting,
	author = {Pareschi, Lorenzo and Toscani, Giuseppe},
	publisher = {OUP Oxford},
	title = {Interacting multiagent systems: kinetic equations and Monte Carlo methods},
	year = {2013}}

@article{toscani98,
	author = {Toscani, G.},
	doi = {10.1051/m2an/1998320607631},
	fjournal = {RAIRO Mod\'elisation Math\'ematique et Analyse Num\'erique},
	issn = {0764-583X},
	journal = {RAIRO Mod\'el. Math. Anal. Num\'er.},
	mrclass = {82C40 (35Q99 45K05)},
	mrnumber = {1652617},
	mrreviewer = {Lorenzo\ Pareschi},
	number = {6},
	pages = {763--772},
	title = {The grazing collisions asymptotics of the non-cut-off {K}ac equation},
	url = {https://doi.org/10.1051/m2an/1998320607631},
	volume = {32},
	year = {1998}}

@article{bobylev00,
	author = {Bobylev, A. V. and Nanbu, K.},
	doi = {10.1103/PhysRevE.61.4576},
	issue = {4},
	journal = {Phys. Rev. E},
	month = {Apr},
	numpages = {0},
	pages = {4576--4586},
	publisher = {American Physical Society},
	title = {Theory of collision algorithms for gases and plasmas based on the {B}oltzmann equation and the {L}andau-{F}okker-{P}lanck equation},
	url = {https://link.aps.org/doi/10.1103/PhysRevE.61.4576},
	volume = {61},
	year = {2000}}

@article{villani98,
	author = {Villani, C.},
	doi = {10.1007/s002050050106},
	fjournal = {Archive for Rational Mechanics and Analysis},
	issn = {0003-9527},
	journal = {Arch. Rational Mech. Anal.},
	mrclass = {82C40 (35D05 35Q99 45K05 76P05)},
	mrnumber = {1650006},
	mrreviewer = {Alexei\ Heintz},
	number = {3},
	pages = {273--307},
	title = {On a new class of weak solutions to the spatially homogeneous {B}oltzmann and {L}andau equations},
	url = {https://doi.org/10.1007/s002050050106},
	volume = {143},
	year = {1998}}

@article{TW_2020,
	author = {Claudia Totzeck and Marie-Therese Wolfram},
	doi = {10.3934/mbe.2020320},
	issn = {1551-0018},
	journal = {Math. Biosci. Eng.},
	number = {5},
	pages = {6026-6044},
	title = {Consensus-based global optimization with personal best},
	url = {https://www.aimspress.com/article/doi/10.3934/mbe.2020320},
	volume = {17},
	year = {2020}}

@inproceedings{nitanda2022convex,
	author = {Nitanda, Atsushi and Wu, Denny and Suzuki, Taiji},
	booktitle = {International Conference on Artificial Intelligence and Statistics},
	organization = {PMLR},
	pages = {9741--9757},
	title = {Convex analysis of the mean field {L}angevin dynamics},
	year = {2022}}

@article{stacy62,
	author = {Stacy, E. W.},
	journal = {Ann. Math. Statist.},
	number = {3},
	pages = {1187-1192},
	title = {A Generalization of the Gamma Distribution},
	volume = {33},
	year = {1962}}

@article{Albi_23,
	author = {Albi, Giacomo and Ferrarese, Federica and Totzeck, Claudia},
	doi = {10.1142/S0218202523500641},
	journal = {Math. Mod. Meth. Appl. Sci.},
	number = {14},
	pages = {2905-2933},
	title = {Kinetic-based optimization enhanced by genetic dynamics},
	url = {https://doi.org/10.1142/S0218202523500641},
	volume = {33},
	year = {2023}}

@article{Ha_CBO20,
	author = {Ha, Seung-Yeal and Jin, Shi and Kim, Doheon},
	doi = {10.1142/S0218202520500463},
	journal = {Math. Mod. Meth. Appl. Sci.},
	number = {12},
	pages = {2417-2444},
	title = {Convergence of a first-order consensus-based global optimization algorithm},
	url = {https://doi.org/10.1142/S0218202520500463},
	volume = {30},
	year = {2020}}

@article{OTTO2000361,
	author = {F. Otto and C. Villani},
	doi = {https://doi.org/10.1006/jfan.1999.3557},
	issn = {0022-1236},
	journal = {J. Funct. Anal.},
	number = {2},
	pages = {361-400},
	title = {Generalization of an Inequality by {T}alagrand and Links with the Logarithmic {S}obolev Inequality},
	url = {https://www.sciencedirect.com/science/article/pii/S0022123699935577},
	volume = {173},
	year = {2000}}

@article{FPTZ22,
	author = {G. Furioli and A. Pulvirenti and E. Terraneo and G. Toscani},
	doi = {10.1007/s00032-022-00352-3},
	journal = {Milan J. Math.},
	pages = {177--208},
	title = {One-Dimensional {F}okker-{P}lanck Equations and Functional Inequalities for Heavy Tailed Densities},
	volume = 90,
	year = 2022}

@article{FPTZ19,
	author = {G. Furioli and A. Pulvirenti and E. Terraneo and G. Toscani},
	doi = {10.1016/J.ANIHPC.2019.07.005},
	journal = {Ann. Inst. H. Poincar{\'e} Anal. Non Lin{\'e}aire},
	pages = {2065--2082},
	title = {{W}right-{F}isher--type equations for opinion formation, large time behavior and weighted logarithmic-{S}obolev inequalities},
	volume = 36,
	year = 2019}

@article{ATZ_23,
	author = {Ferdinando Auricchio and Giuseppe Toscani and Mattia Zanella},
	doi = {https://doi.org/10.1016/j.aml.2023.108746},
	issn = {0893-9659},
	journal = {Appl. Math. Lett.},
	pages = {108746},
	title = {Trends to equilibrium for a nonlocal {F}okker-{P}lanck equation},
	url = {https://www.sciencedirect.com/science/article/pii/S0893965923001787},
	volume = {145},
	year = {2023}}

@article{pareschi_russo,
	author = {Pareschi, L. and Russo, G.},
	doi = {10.1051/proc:2001004},
	journal = {ESAIM: Proc.},
	pages = {35-75},
	title = {An introduction to {M}onte {C}arlo method for the {B}oltzmann equation},
	url = {https://doi.org/10.1051/proc:2001004},
	volume = 10,
	year = 2001}

@article{CT07,
	author = {Carrillo, J. A. and Toscani, G.},
	fjournal = {Rivista di Matematica della Universit\`a{} di Parma. Serie 7},
	issn = {0035-6298},
	journal = {Riv. Mat. Univ. Parma (7)},
	mrclass = {82C40 (35B40 35F20 60B10 82-01)},
	mrnumber = {2355628},
	mrreviewer = {Carlo\ Cercignani},
	pages = {75--198},
	title = {Contractive probability metrics and asymptotic behavior of dissipative kinetic equations},
	volume = {6},
	year = {2007}}

@article{Toscani99,
	author = {Toscani, G.},
	doi = {10.1090/qam/1704435},
	fjournal = {Quarterly of Applied Mathematics},
	issn = {0033-569X,1552-4485},
	journal = {Quart. Appl. Math.},
	mrclass = {82C31 (35Q99)},
	mrnumber = {1704435},
	mrreviewer = {Carlo\ Cercignani},
	number = {3},
	pages = {521--541},
	title = {Entropy production and the rate of convergence to equilibrium for the {F}okker-{P}lanck equation},
	url = {https://doi.org/10.1090/qam/1704435},
	volume = {57},
	year = {1999}}

@article{CCTT,
	author = {Carrillo, Jos\'e{} A. and Choi, Young-Pil and Totzeck, Claudia and Tse, Oliver},
	doi = {10.1142/S0218202518500276},
	fjournal = {Mathematical Models and Methods in Applied Sciences},
	issn = {0218-2025,1793-6314},
	journal = {Math. Models Methods Appl. Sci.},
	mrclass = {90C26 (35Q84 37N40 60H10 90C59)},
	mrnumber = {3804923},
	number = {6},
	pages = {1037--1066},
	title = {An analytical framework for consensus-based global optimization method},
	url = {https://doi.org/10.1142/S0218202518500276},
	volume = {28},
	year = {2018}}

@article{blum2003metaheuristics,
	author = {Blum, Christian and Roli, Andrea},
	journal = {ACM computing surveys (CSUR)},
	number = {3},
	pages = {268--308},
	publisher = {Acm New York, NY, USA},
	title = {Metaheuristics in combinatorial optimization: Overview and conceptual comparison},
	volume = {35},
	year = {2003}}

@article{MR3629153,
	author = {Bellomo, Nicola and Ha, Seung-Yeal},
	doi = {10.1142/S0218202517500154},
	fjournal = {Mathematical Models and Methods in Applied Sciences},
	issn = {0218-2025,1793-6314},
	journal = {Math. Models Methods Appl. Sci.},
	mrclass = {92D50 (35Q20 91D10)},
	mrnumber = {3629153},
	number = {4},
	pages = {745--770},
	title = {A quest toward a mathematical theory of the dynamics of swarms},
	url = {https://doi.org/10.1142/S0218202517500154},
	volume = {27},
	year = {2017}}

@article{MR3352763,
	author = {Bisi, Marzia and Ca\~nizo, Jos\'e{} A. and Lods, Bertrand},
	doi = {10.1016/j.jfa.2015.05.002},
	fjournal = {Journal of Functional Analysis},
	issn = {0022-1236,1096-0783},
	journal = {J. Funct. Anal.},
	mrclass = {82C40 (47G20)},
	mrnumber = {3352763},
	number = {4},
	pages = {1028--1069},
	title = {Entropy dissipation estimates for the linear {B}oltzmann operator},
	url = {https://doi.org/10.1016/j.jfa.2015.05.002},
	volume = {269},
	year = {2015}}

@article{MR1165528,
	author = {Desvillettes, L.},
	doi = {10.1080/00411459208203923},
	fjournal = {Transport Theory and Statistical Physics},
	issn = {0041-1450,1532-2424},
	journal = {Transport Theory Statist. Phys.},
	mrclass = {82C40 (76P05)},
	mrnumber = {1165528},
	mrreviewer = {Reinhard\ Illner},
	number = {3},
	pages = {259--276},
	title = {On asymptotics of the {B}oltzmann equation when the collisions become grazing},
	url = {https://doi.org/10.1080/00411459208203923},
	volume = {21},
	year = {1992}}

@article{chak2023generalized,
	author = {Chak, Martin and Kantas, Nikolas and Pavliotis, Grigorios A},
	journal = {SIAM/ASA J. Uncert. Quantif.},
	number = {1},
	pages = {139--167},
	publisher = {SIAM},
	title = {On the generalized {L}angevin equation for simulated annealing},
	volume = {11},
	year = {2023}}

@article{chizat2022mean,
	author = {Chizat, L{\'e}na{\"\i}c},
	journal = {arXiv preprint arXiv:2202.01009},
	title = {Mean-field {L}angevin dynamics: Exponential convergence and annealing},
	year = {2022}}

@article{MR942621,
	author = {Hajek, Bruce},
	doi = {10.1287/moor.13.2.311},
	fjournal = {Mathematics of Operations Research},
	issn = {0364-765X,1526-5471},
	journal = {Math. Oper. Res.},
	mrclass = {90C30},
	mrnumber = {942621},
	mrreviewer = {Th.\ M.\ Liebling},
	number = {2},
	pages = {311--329},
	title = {Cooling schedules for optimal annealing},
	url = {https://doi.org/10.1287/moor.13.2.311},
	volume = {13},
	year = {1988}}

@article{MR1188544,
	author = {B\'elisle, Claude J. P.},
	doi = {10.2307/3214721},
	fjournal = {Journal of Applied Probability},
	issn = {0021-9002,1475-6072},
	journal = {J. Appl. Probab.},
	mrclass = {65C05 (65K10 90C30)},
	mrnumber = {1188544},
	mrreviewer = {R.\ Shonkwiler},
	number = {4},
	pages = {885--895},
	title = {Convergence theorems for a class of simulated annealing algorithms on {${\mathbb{R}}^d$}},
	url = {https://doi.org/10.2307/3214721},
	volume = {29},
	year = {1992}}

@book{MR983115,
	author = {Aarts, Emile and Korst, Jan},
	isbn = {0-471-92146-7},
	mrclass = {90C27 (68Q99 92A09)},
	mrnumber = {983115},
	mrreviewer = {Th.\ M.\ Liebling},
	note = {A stochastic approach to combinatorial optimization and neural computing},
	pages = {xii+272},
	publisher = {John Wiley \& Sons, Ltd., Chichester},
	series = {Wiley-Interscience Series in Discrete Mathematics and Optimization},
	title = {Simulated annealing and {B}oltzmann machines},
	year = {1989}}

@article{PTTM,
	author = {Pinnau, Ren{\'e} and Totzeck, Claudia and Tse, Oliver and Martin, Stephan},
	journal = {Math. Mod. Meth. Appl. Sci.},
	number = {01},
	pages = {183--204},
	publisher = {World Scientific},
	title = {A consensus-based model for global optimization and its mean-field limit},
	volume = {27},
	year = {2017}}

@article{P24,
	author = {Pareschi, Lorenzo},
	doi = {10.1142/S0218202524500428},
	journal = {Math. Mod. Meth. Appl. Sci.},
	number = {12},
	pages = {2191-2216},
	title = {Optimization by linear kinetic equations and mean-field {L}angevin dynamics},
	volume = {34},
	year = {2024}}

@book{PT,
	author = {Pareschi, L. and Toscani, G.},
	publisher = {OUP Oxford},
	title = {Interacting {M}ultiagent {S}ystems: {K}inetic {E}quations and {M}onte {C}arlo {M}ethods},
	year = {2013}}

@incollection{Tot,
	author = {Totzeck, Claudia},
	booktitle = {Active particles. {V}ol. 3. {A}dvances in theory, models, and applications},
	doi = {10.1007/978-3-030-93302-9\_6},
	isbn = {978-3-030-93301-2; 978-3-030-93302-9},
	mrclass = {93A16 (93D50 93E03)},
	mrnumber = {4433532},
	pages = {201--226},
	publisher = {Birkh\"auser/Springer, Cham},
	series = {Model. Simul. Sci. Eng. Technol.},
	title = {Trends in consensus-based optimization},
	url = {https://doi.org/10.1007/978-3-030-93302-9_6},
	year = {[2022] \copyright 2022}}

@article{Tosc,
	author = {Toscani, Giuseppe},
	doi = {10.1007/s11587-019-00471-x},
	fjournal = {Ricerche di Matematica. A Journal of Pure and Applied Mathematics},
	issn = {0035-5038,1827-3491},
	journal = {Ric. Mat.},
	mrclass = {35Q84 (37A30 39B62 60J60 82B21)},
	mrnumber = {4272007},
	number = {1},
	pages = {35--50},
	title = {Entropy-type inequalities for generalized {G}amma densities},
	url = {https://doi.org/10.1007/s11587-019-00471-x},
	volume = {70},
	year = {2021}}

@book{MR2244940,
	author = {Nocedal, Jorge and Wright, Stephen J.},
	edition = {Second},
	isbn = {978-0387-30303-1; 0-387-30303-0},
	mrclass = {90-01 (49Mxx 65K05 90-02 90C30)},
	mrnumber = {2244940},
	pages = {xxii+664},
	publisher = {Springer, New York},
	series = {Springer Series in Operations Research and Financial Engineering},
	title = {Numerical {O}ptimization},
	year = {2006}}

@book{MR1678201,
	author = {Kelley, C. T.},
	doi = {10.1137/1.9781611970920},
	isbn = {0-89871-433-8},
	mrclass = {90-02 (65K05 90-01 90C30)},
	mrnumber = {1678201},
	mrreviewer = {William\ W.\ Hager},
	pages = {xvi+180},
	publisher = {Society for Industrial and Applied Mathematics (SIAM), Philadelphia, PA},
	series = {Frontiers in Applied Mathematics},
	title = {Iterative methods for optimization},
	url = {https://doi.org/10.1137/1.9781611970920},
	volume = {18},
	year = {1999}}

@article{KGV,
	author = {Kirkpatrick, S. and Gelatt, Jr., C. D. and Vecchi, M. P.},
	doi = {10.1126/science.220.4598.671},
	fjournal = {American Association for the Advancement of Science. Science},
	issn = {0036-8075,1095-9203},
	journal = {Science},
	mrclass = {90C27 (82A05)},
	mrnumber = {702485},
	number = {4598},
	pages = {671--680},
	title = {Optimization by simulated annealing},
	url = {https://doi.org/10.1126/science.220.4598.671},
	volume = {220},
	year = {1983}}

@article{GH,
	author = {Geman, Stuart and Hwang, Chii-Ruey},
	doi = {10.1137/0324060},
	fjournal = {SIAM J. Contr. Optim.},
	issn = {0363-0129},
	journal = {SIAM J. Control Optim.},
	mrclass = {49D10 (60J70)},
	mrnumber = {854068},
	mrreviewer = {Pierre-Louis\ Lions},
	number = {5},
	pages = {1031--1043},
	title = {Diffusions for global optimization},
	url = {https://doi.org/10.1137/0324060},
	volume = {24},
	year = {1986}}

@article{CJLZ,
	author = {Carrillo, Jos\'e{} A. and Jin, Shi and Li, Lei and Zhu, Yuhua},
	doi = {10.1051/cocv/2020046},
	fjournal = {ESAIM. Control, Optimisation and Calculus of Variations},
	issn = {1292-8119,1262-3377},
	journal = {ESAIM Control Optim. Calc. Var.},
	mrclass = {60H35 (65K10 68T05 70F10)},
	mrnumber = {4222159},
	pages = {Paper No. S5, 22},
	title = {A consensus-based global optimization method for high dimensional machine learning problems},
	url = {https://doi.org/10.1051/cocv/2020046},
	volume = {27},
	year = {2021}}

@article{MR4793478,
	author = {Fornasier, Massimo and Klock, Timo and Riedl, Konstantin},
	doi = {10.1137/22M1527805},
	fjournal = {SIAM Journal on Optimization},
	issn = {1052-6234,1095-7189},
	journal = {SIAM J. Optim.},
	mrclass = {65K10 (35Q84 35Q90 90C26 90C56)},
	mrnumber = {4793478},
	mrreviewer = {Aurea\ Mart\'inez},
	number = {3},
	pages = {2973--3004},
	title = {Consensus-based optimization methods converge globally},
	url = {https://doi.org/10.1137/22M1527805},
	volume = {34},
	year = {2024}}

@misc{borghi2024kineticmodelsoptimizationunified,
	archiveprefix = {arXiv},
	author = {Giacomo Borghi and Michael Herty and Lorenzo Pareschi},
	eprint = {2410.10369},
	primaryclass = {math.OC},
	title = {Kinetic models for optimization: a unified mathematical framework for metaheuristics},
	url = {https://arxiv.org/abs/2410.10369},
	year = {2024}}

\end{document}